\documentclass{article}
\usepackage[table]{xcolor}
\usepackage{amssymb,amsmath,amsthm,tikz,multirow,enumitem,nccrules,geometry,graphicx}
\usetikzlibrary{calc, arrows, arrows.meta}
\title{Tiling of Hyperbolic Surface by Maximal Number of Tiles}
\author{Chunlin Li, Erxiao Wang\thanks{Corresponding author (wang.eric@zjnu.edu.cn).  Research was supported by National Natural Science Foundation of China NSFC-RGC 12361161603 and Key Projects of Zhejiang Natural Science Foundation LZ22A010003.}, Jie Wu,
Zhejiang Normal University \\
Min Yan\thanks{Research was supported by NSFC-RGC Joint Research Scheme N-HKUST607/23 and Hong Kong RGC General Research Fund 16310925.}, 
Hong Kong University of Science and Technology}

\usepackage[hidelinks, hyperindex]{hyperref}

\newcommand{\sub}{\subset}
\newcommand{\pa}{\partial}

\newcommand{\mc}{\mathcal}
\newcommand{\bb}{\mathbb}

\newcommand{\hbullet}{\tikz \fill (0.08,0) arc (0:180:0.08);}
\newcommand*\circled[1]{\tikz[baseline=(char.base)]{
		\node[shape=circle,draw,inner sep=0.5pt] (char) {#1};}}

\colorlet{Dbrown}{brown!50!black!70}
\colorlet{Gold}{orange!75!yellow!75}

\makeatletter
\newcommand{\drawhypgeodesic}[5][]{
  
  \pgfmathsetmacro{\xone}{#2}
  \pgfmathsetmacro{\yone}{#3}
  \pgfmathsetmacro{\xtwo}{#4}
  \pgfmathsetmacro{\ytwo}{#5}
  \pgfmathsetmacro{\denom}{\xone*\ytwo - \xtwo*\yone}

  \ifdim\denom pt < 0.00001pt
    \ifdim\denom pt > -0.00001pt
      \draw[#1] (#2,#3) -- (#4,#5);
      \def\isLine{1}
    \else
      \def\isLine{0}
    \fi
  \else
    \def\isLine{0}
  \fi
  
  \ifnum\isLine=0

    \pgfmathsetmacro{\mone}{1 + \xone*\xone + \yone*\yone}
    \pgfmathsetmacro{\mtwo}{1 + \xtwo*\xtwo + \ytwo*\ytwo}
    \pgfmathsetmacro{\dd}{(\yone*\mtwo - \ytwo*\mone)/\denom}
    \pgfmathsetmacro{\ee}{(\xtwo*\mone - \xone*\mtwo)/\denom}
    \pgfmathsetmacro{\cx}{-\dd/2}
    \pgfmathsetmacro{\cy}{-\ee/2}
    \pgfmathsetmacro{\rr}{sqrt(\cx*\cx + \cy*\cy - 1)}
    \pgfmathsetmacro{\thetaone}{atan2(\yone-\cy, \xone-\cx)}
    \pgfmathsetmacro{\thetatwo}{atan2(\ytwo-\cy, \xtwo-\cx)}   
    \pgfmathsetmacro{\dtheta}{\thetatwo - \thetaone}
    \pgfmathsetmacro{\dtheta}{\dtheta > 180 ? \dtheta - 360 : \dtheta}
    \pgfmathsetmacro{\dtheta}{\dtheta <= -180 ? \dtheta + 360 : \dtheta}
    
    \draw[#1] (#2, #3) arc (\thetaone : \thetaone + \dtheta : \rr);
  \fi
}
\makeatother

\newtheorem{theorem}{Theorem}
\newtheorem{lemma}[theorem]{Lemma}

\newtheorem{proposition}[theorem]{Proposition}
\newtheorem*{theorem*}{Theorem}

\theoremstyle{definition}
\newtheorem{definition}{Definition}
\newtheorem*{definition*}{Definition}
\newtheorem*{case*}{Case}
\newtheorem*{subcase*}{Subcase}

\theoremstyle{remark}

\numberwithin{equation}{section}

\begin{document}

\maketitle

\begin{abstract}
Tilings of a closed hyperbolic surface by congruent $n$-gons with $n\ge 7$ is a finite problem. One extreme of the finite problem is that the number of tiles is maximal. The problem is related to edge-to-edge tilings of the Poincar\'e disk ${\bb H}^2$ by congruent $n$-gons, such that all vertices have degree 3 and the vertex sets  satisfy certain balance condition. For the special case $n=7$, we show that the heptagon in any such tiling belongs to one of the four families, and tilings of ${\bb H}^2$ by generic heptagons in each family is unique. As a consequence, for any closed hyperbolic surface of Euler characteristic $\chi$, we may classify all tilings of the surface by maximal number $-6\chi$ of congruent generic heptagons in the four families. 


{\em 2020 MR Subject Classification}: 52C20, 05B45, 51M10, 51M20.
\end{abstract}

\section{Overview}
\label{intro}

Consider tilings of closed hyperbolic surfaces by polygons congruent to a {\em prototile} $n$-gon $P$, with $n\ge 7$. Here hyperbolic surfaces have constant curvature $-1$, and edges of $P$ are geodesics. Then the Euler characteristics of surfaces $\chi<0$. We show in Section \ref{finite} that for fixed $\chi$, under the natural assumption of at least three tiles at each vertex, both $n$ and the number of tiles $f$ have only finite number of choices. For example, we have $n\le 18$ for $\chi=-2$ (the surface is $2{\bb T}^2$ or $4{\bb P}^2$). Moreover, we have $f\le-6\chi$ for $n=7$. Therefore the problem is finite.

The first step in studying the finite problem is the extreme cases. In Lemma \ref{bound}, we prove that the number of tiles $f$ equals the maximal possible number $\frac{-6\chi}{n-6}$ if and only if the tiling is edge-to-edge, and all vertices have degree 3. The goal of this paper is the classification for the case $n=7$, i.e., tilings of hyperbolic surfaces by the maximal number $f=-6\chi$ of congruent heptagons. The following result describes the heptagons in such tilings.

\begin{theorem}\label{thm1}
A heptagon can tile a closed hyperbolic surface of Euler characteristic $\chi<0$, such that the number of tiles is the maximal possible $-6\chi$, if and only if the heptagon belongs to one of the four families described by the following edge and angle equalities:
\begin{align}
7_{2+3} \colon &
|\bar{0}|=|\bar{6}|, \;
|\bar{2}|=|\bar{5}|, \;
[0]+[0]+[0]=[1]+[2]+[6]=[3]+[4]+[5]=2\pi; \label{h5ev} \\
7_{3+3} \colon &
|\bar{0}|=|\bar{6}|,\;
|\bar{1}|=|\bar{3}|,\;
|\bar{2}|=|\bar{5}|,\;
[0]+[0]+[0]=[1]+[3]+[6]=[2]+[4]+[5]=2\pi; \label{h6ev} \\
7^1_{5+6} \colon &
|\bar{0}|=|\bar{2}|=|\bar{3}|=|\bar{4}|=|\bar{6}|,\;
|\bar{1}|=|\bar{5}|, \nonumber \\
&  [0]+[1]+[2]=[0]+[6]+[6]=[1]+[3]+[5]=[2]+[3]+[6] \nonumber \\
& =[3]+[4]+[4]=[4]+[5]+[5]=2\pi; \label{h11ev1} \\
7^2_{5+6} \colon &
|\bar{0}|=|\bar{1}|=|\bar{2}|=|\bar{4}|=|\bar{5}|=|\bar{6}|, \nonumber \\
& [0]+[1]+[2]=[0]+[1]+[6]=[0]+[5]+[5]=[2]+[2]+[6]= \nonumber \\
& [3]+[4]+[5]=[4]+[4]+[6]=2\pi. \label{h11ev2}
\end{align}
\end{theorem}

We remark that the assumption on heptagon is that it tiles {\em some} closed hyperbolic surface with maximal number of tiles. In fact, Theorem \ref{thm5} says that a generic (see Definition \ref{generic}) $7_{2+3}$-heptagon cannot tile non-orientable surfaces, and a generic $7_{3+3}$-heptagon cannot tile $g{\bb P}^2$ for odd $g$. 

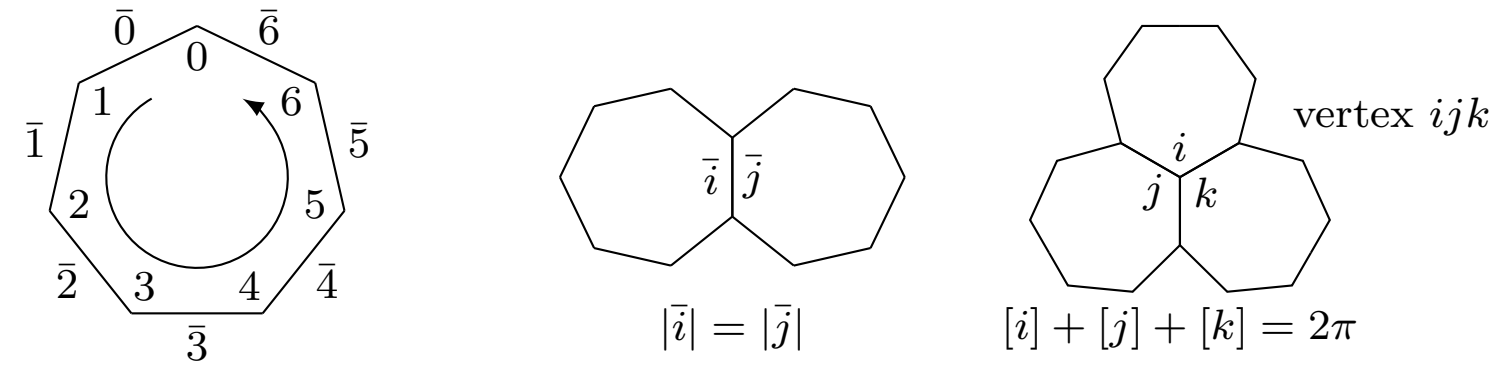
\begin{figure}[htp]
\centering
\begin{tikzpicture}[>=latex, scale=1, font=\footnotesize]

\foreach \a in {0,...,6}
{
\draw
	(38.57+51.43*\a:1) -- (90+51.43*\a:1);

\node at (90+51.43*\a:0.8) {\a};
\node at (115.7+51.43*\a:1.1) {$\bar{\a}$};

}

\draw[->]
	(120:0.6) arc (120:420:0.6);

\begin{scope}[xshift=3cm]

\foreach \a in {0,...,6}
\foreach \x in {0,1}
\draw[xshift=1.08*\x cm]
	(-25.714+25.714*\x+51.43*\a:0.6) -- (25.714+25.714*\x+51.43*\a:0.6);
	
\node at (0.4,0) { $\bar{i}$};
\node at (0.68,0) {$\bar{j}$};
\node at (0.54,-1) {$|\bar{i}|=|\bar{j}|$};

\end{scope}

\begin{scope}[xshift=6.5cm]

\foreach \a in {0,1,2}
\foreach \x in {1,-1}
\draw[rotate=120*\a, xscale=\x]
	(0,0) -- (30:0.45) -- (0.5,0.65) -- (0.25,1) -- (0,1);
	
\node at (90:0.2) {$i$};
\node at (210:0.2) {$j$};
\node at (-30:0.2) {$k$};
\node at (0,-1) {$[i]+[j]+[k]=2\pi$};
\node at (1.4,0.4) {vertex $ijk$};

\end{scope}
		
\end{tikzpicture}
\caption{Corners and edges of a prototile heptagon, and edge and angle equalities.} 
\label{prototile}
\end{figure}

The statement of Theorem \ref{thm1} includes some notations and terminologies that we need to explain. On the left of Figure \ref{prototile}, we circularly label the {\em corners} of the prototile heptagon $P$ by $i\in {\bb Z}_7$. We denote the {\em angle value} of the corner $i$ by $[i]$. We also denote the edge connecting $i$ to $i+1$ by $\bar{i}$, and denote the {\em length} of the edge by $|\bar{i}|$.  The isometry between a tile and $P$ translates the labels of the corners of $P$ to the labels of the corners of the tile. 

The labels are not intrinsic part of the tiling. Therefore labels can be changed by the {\em rotation} $i\mapsto c+i$ or the {\em flip} $i\mapsto c-i$. The discussions and results in this paper are stated up to such {\em relabellings}.

The equalities in Theorem \ref{thm1} are the geometrical conditions that $P$ should satisfy. The equalities consist of {\em edge equalities} of the form $|\bar{i}|=|\bar{j}|$, with $i\ne j$, and {\em angle equalities} of the form $[i]+[j]+[k]=2\pi$. The edge equalities come from the edges shared by adjacent tiles in edge-to-edge tilings, in the middle of Figure \ref{prototile}. Therefore $|\bar{i}|=|\bar{j}|$ means that the $i$-th and $j$-th edges can be glued together. All the edge equalities form an {\em edge system} $E$, defined in Section \ref{edge} as a collection of equalities of the form $|\bar{i}|=|\bar{j}|$. The angle equalities come from the sums of the angle values of the all corners at vertices, on the right of Figure \ref{prototile}. Therefore $[i]+[j]+[k]=2\pi$ means that $i$-th, $j$-th and $k$-th corners can form a {\em vertex} $ijk$ in the tiling. In general, an angle equality $[i_1]+[i_2]+\dots+[i_d]=2\pi$ is equivalent to a possible vertex $i_1i_2\cdots i_d$ in a tiling. Therefore we use the angle equality $[i_1]+[i_2]+\dots+[i_d]=2\pi$ and the vertex notation $i_1i_2\cdots i_d$ interchangeably in this paper. Of course, the assumption of the maximal number of tiles restricts to $d=3$. All the angle equalities form a {\em corner system} $V$, in Section \ref{corner}. 

Both edge system $E$ and corner system $V$ consist of linear equalities, and the rank of the whole edge and corner system $\Sigma=E\cup V$ (in Section \ref{edge-corner}) is the sum of the ranks of $E$ and $V$. In \eqref{h5ev}, we have $\text{rank}E=2$ and $\text{rank}V=3$, which is reflected in the subscript $2+3$ of the notation $7_{2+3}$ for the heptagon. Of course, $7$ in $7_{2+3}$ denotes the heptagon. We remark that, if $\text{rank}\Sigma\ge 13$, then it is easy to see that only the regular heptagon $P_7(\frac{2\pi}{3})$ of angle $\frac{2\pi}{3}$ and edge length $2\text{arccosh}(\frac{2}{\sqrt{3}}\cos\frac{\pi}{7})$ satisfies $\Sigma$. Therefore we only need to consider edge and corner systems of ranks up to 12.

A tiling ${\mc T}$ by a heptagonal prototile $P$ gives an edge and corner system $\Sigma({\mc T})=E({\mc T})\cup V({\mc T})$ satisfied by $P$. In fact, the systems in Theorem \ref{thm1} are $\Sigma({\mc T})$ for the unique tilings ${\mc T}$ in Theorem \ref{thm3}. In particular, the corner systems in Theorem \ref{thm1} are exactly the vertices appearing in ${\mc T}$.

Conversely, if a polygon $P$ satisfies an edge and corner system $\Sigma=E\cup V$, what are the tilings by congruent copies of $P$? Two questions immediate arise: (1) Whether $E$ provides the only ways the edges of $P$ can be glued together; (2) Whether $V$ provides the only ways the corners of $P$ can form vertices. For edge-to-edge tilings (which is the case in the context of this paper), equal edge length is the only criterion for edges to be glued together. This means that the answer to (1) is affirmative if and only if $P$ does not satisfy any edge equality not in $E$. On the other hand, angle sum being $2\pi$ is the only criterion for corners to form a vertex. Therefore affirmative answer to both (1) and (2) means the following. 

\begin{definition}\label{generic}
A heptagon is {\em generic} with respect to an edge and corner system $\Sigma$, if it satisfies $\Sigma$, and does not satisfy any edge equality (of the form $|\bar{i}|=|\bar{j}|$) or angle equality (of the form $[i_1]+[i_2]+\dots+[i_d]=2\pi$) that is not in $\Sigma$. 
\end{definition}

For any polygon $P$, let $E$ and $V$ be all the edge and angle equalities satisfied by $P$. Then $P$ is generic  with respect to the edge and corner system $E\cup V$, and is not generic with respect to any other system. For example, the regular heptagon $P_7(\frac{2\pi}{3})$ satisfies all four systems in Theorem \ref{thm1}, and is not generic with respect to any one of the four systems. It is generic only with respect to the edge system $|\bar{i}|=|\bar{j}|$ for all $i,j\in {\bb Z}_7$ together with the corner system $[i]+[j]+[k]=2\pi$ for all $i,j,k\in {\bb Z}_7$.

If $P$ satisfies a corner system $V$, then $P$ also satisfies any angle equality that is a linear combination of angle equalities in $V$. Therefore $V$ in the definition above is necessarily complete in the following sense. 

\begin{definition}\label{completion}
Suppose $V$ is a corner system. The {\em completion} $\bar{V}$ of $V$ consists of all the angle equalities $[i_1]+[i_2]+\dots+[i_d]=2\pi$ that are linear combinations of equalities in $V$. We say $V$ is {\em complete} if $\bar{V}=V$.
\end{definition}

The corner systems in \eqref{h5ev}, \eqref{h6ev}, \eqref{h11ev1} are complete. We also find examples \eqref{5generic} and \eqref{6generic} of heptagons that are generic with respect to \eqref{h5ev} and \eqref{h6ev}. Therefore, in a tiling by a generic $7_{2+3}$-heptagon, we know only $\bar{0}$ and $\bar{6}$ can be glued together,  only $\bar{2}$ and $\bar{5}$ can be glued together, and $0^3=000$, $126$, $345$ are the only possible vertices. This provides the guidance for constructing tilings by $P$. The comment also applies to generic $7_{3+3}$-heptagon.

Generic heptagons with respect to \eqref{h11ev1} and \eqref{h11ev2} do not exist. In fact, both systems have high enough rank, that they completely determine {\em unique non-regular} $7^1_{5+6}$- and $7^2_{5+6}$-heptagon, with edge lengths and angle values given by \eqref{111heptagon} and \eqref{112heptagon}. The angle values in \eqref{111heptagon} imply that the non-regular $7^1_{5+6}$-heptagon also satisfies angle equalities of degree $>3$, such as $[0]+[0]+[2]+[5]=2\pi$ in the complete corner system \eqref{h111v}. Therefore the non-regular $7^1_{5+6}$-heptagon is generic only with respect to the edge system in \eqref{h11ev1} together with the complete corner system \eqref{h111v}. For the non-regular $7^2_{5+6}$-heptagon, the problem is that the corner system in \eqref{h11ev2} is not complete. The reason is that the angle equalities imply $[2]=[4]=[6]=\tfrac{2\pi}{3}$, which means that $[i]+[j]+[k]=2\pi$ is in the completion for any $i,j,k\in \{2,4,6\}$. The completion of the corner system \eqref{h11ev2} is \eqref{h112v}, in which any one $\alpha\in \{2,4,6\}$. The values \eqref{112heptagon} for the unique non-regular $7^2_{5+6}$-heptagon shows that it is generic with respect to the edge system in \eqref{h11ev2} together with the complete corner system \eqref{h112v}.

In Section \ref{rank9}, we need to consider an edge and corner system \eqref{h9ev}. The corner part consists of $01^2,06^2,016$, and $ijk$ for any $i,j,k\in \{2,3,4,5\}$, and is complete. However, Proposition \ref{h9hidden} shows that there is a ``hidden equality'' $|\bar{1}|=|\bar{5}|$ implied by \eqref{h9ev}, for geometrical instead of linear algebra reason. Since $|\bar{1}|=|\bar{5}|$ is not in \eqref{h9ev}, generic heptagons with respect to \eqref{h9ev} do not exist. In fact, Proposition \ref{h9hidden} says that any heptagon satisfying the system \eqref{h9ev} also satisfies the system \eqref{h5ev}. 

\begin{definition}
We say an edge and corner system $\Sigma'$ {\em implies} another edge and corner system $\Sigma$, if a polygon that satisfies $\Sigma'$ also satisfies $\Sigma$.
\end{definition}

The simplest implication is by linear algebra, which is equivalent to $\text{rank}(\Sigma'\cup\Sigma)=\text{rank}\Sigma'$. However, the implication can also be by geometrical argument, as illustrated in Proposition \ref{h9hidden}. The most extreme argument is to actually calculate the heptagon by hyperbolic geometry. For many edge and corner systems of ranks 11 or 12, the calculation shows that the only heptagon satisfying the system is the regular heptagon $P_7(\frac{2\pi}{3})$.

Next we discuss how to prove Theorem \ref{thm1}. The universal cover of a hyperbolic surface tiling is a tiling of the Poincar\'e disk ${\bb H}^2$, with the same edge and corner system. Theorem \ref{thm1} is a consequence of the similar result on tilings of the Poincar\'e disk.

\begin{theorem}\label{thm2}
Suppose $\Sigma=E\cup V$ is an edge and corner system, with balanced $V$. Suppose any heptagon satisfying $\Sigma$ can tile the Poincar\'e disk ${\bb H}^2$, such that the tiling is edge-to-edge and all vertices have degree 3. Then up to relabelling, $\Sigma$ implies one of the four edge and corner systems in Theorem \ref{thm1}.
\end{theorem}

The concept of {\em balanced} corner system is defined in Definition \ref{balance} in Section \ref{corner}. The condition is a consequence of the fact that, in a tiling of closed surface,  each corner label appears the same finite number of times (in fact equal to the number of tiles). Therefore balanced corner system is a necessary condition for tilings of closed surfaces and their universal covers. 

Theorem \ref{thm2} is similar to the following well known classification of convex polygons that can tile the plane. 

\begin{theorem*}
There are exactly fifteen families of convex pentagons that can tile the plane, and exactly three families of convex hexagons that can tile the plane.
\end{theorem*}

The classification has a fascinating and long history \cite{zong}. The problem is trivial for triangle and quadrilateral, is non-trivial for hexagon, and is extremely difficult for pentagon. The classification for convex pentagon was finally settled by Mann, McLoud-Mann and von Derau \cite{mmd}, and Rao \cite{rao}. The classification for convex hexagon was universally considered as already proved by Reinhardt \cite{reinhardt} in 1918. In fact, all the existing proofs, including Reinhardt \cite{reinhardt} and Bollob\'as \cite{bollobas}, are significantly incomplete. The proof is only recently completed by Zhu, Wang and Yan \cite{zwy}.

\begin{figure}[htp]
\centering
\begin{tikzpicture}[>=latex,scale=1]

%% 6

\begin{scope}[xshift=-5cm, scale=0.8]

\foreach \a in {0,...,5}
\foreach \b in {0,...,5}
\foreach \c/\d in {0/1.5, 30/0.866, 30/1.732}
\draw[rotate=60*\b, shift={(\c:\d)},rotate=60*\a]
	(60:0.5) -- (0:0.5);

\node at (0,0) {1};

\foreach \a in {2,...,7}
\node at (-90+60*\a:0.866) {\a};

\foreach \a in {8,10,12,14,16,18}
\node at (120+30*\a:1.5) {\a};

\foreach \a in {9,11,13,15,17,19}
\node at (120+30*\a:1.732) {\a};

\end{scope}

%% 7

\foreach \a in {0,...,6}
\draw[rotate=51.43*\a]
	(78:1.2) -- (90:0.9) -- (90:0.5) -- (38.57:0.5) -- (38.57:0.9) -- (50.57:1.2) -- (64.286:1.3) -- (78:1.2)
	(78:1.2) -- (82:1.5) -- (87:1.7) -- (93:1.7) -- (98:1.5) -- (102:1.2)
	(82:1.5) -- (78:1.75) -- (74:1.9) -- (68.286:1.85) -- (64.286:1.6) -- (64.286:1.3)
	(98:1.5) -- (102:1.75) -- (106:1.9) -- (111.714:1.85) -- (115.714:1.6);

\node at (0,0) {1};

\foreach \a in {2,...,8}
\node at (51.43*\a-90:0.85) {\a};

\foreach \a in {9,12,15,18,21,24,27}
\node at (17.143*\a-115.714:1.3) {\footnotesize \a};

\foreach \a in {10,13,16,19,22,25,28}
\node at (17.143*\a-115.714:1.5) {\scriptsize \a};
\foreach \a in {11,14,17,20,23,26,29}
\node at (17.143*\a-115.714:1.5) {\scriptsize \a};

\end{tikzpicture}
\caption{${\mc D}_2$-patches for hexagon and heptagon.}
\label{2layers}
\end{figure}
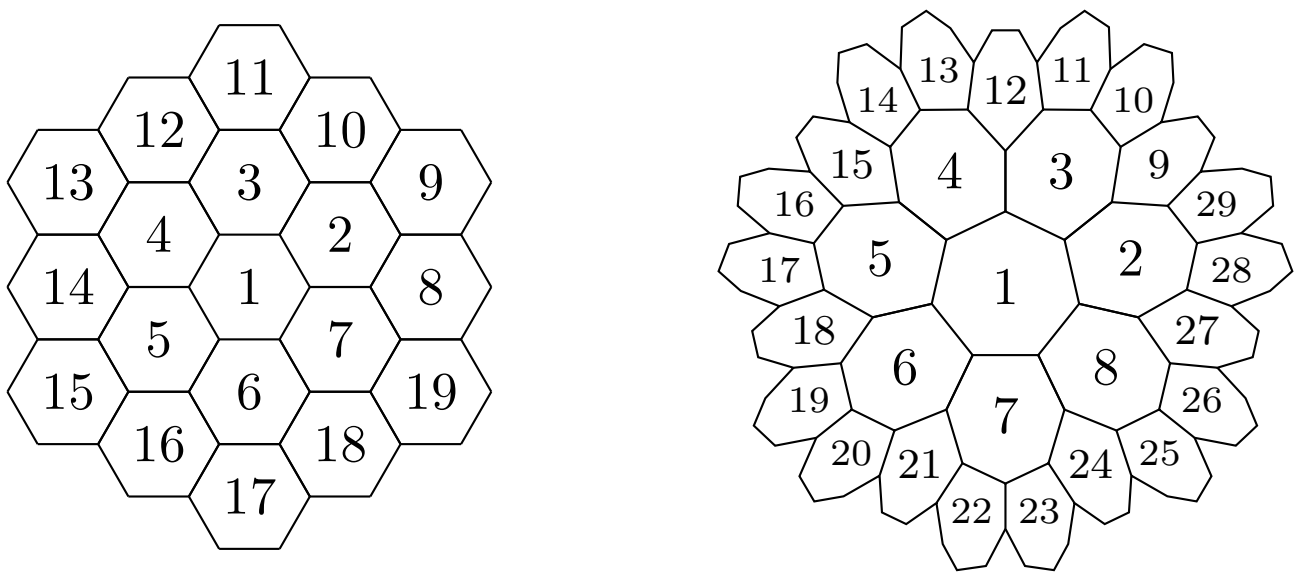

The hexagonal case is simpler than the pentagonal case because there is a lemma (Theorem 1 of \cite{bollobas}, and Lemma 1 of \cite{zwy}) similar to the maximal case of the finiteness Lemma \ref{bound} in Section \ref{finite}: A tiling of the plane by congruent hexagons has arbitrarily large patch, such that the patch is edge-to-edge and all vertices in the patch have degree 3. Combinatorially, the patch is part of the uniform tiling $\{6,3\}$ of the plane, on the left of figure \ref{2layers}. The picture has a central tile labeled 1, and first layer of tiles labeled 2 through 7, and second layer of tiles labeled 8 through 19. There are similar third, fourth layers, and so on. We call the central tile and the first layer together the ${\mc D}_1$-patch (tiles 1 through 7), and call the central tile, the first layer, and second layer together the ${\mc D}_2$-patch (tiles 1 through 19), and so on. 

The key idea in the proof \cite{zwy} for hexagonal tilings of the plane is quite simple: Any such tiling contains a ${\mc D}_2$-patch. We may read off a system of edge equalities $E$ and corner equalities $V$ from any possible ${\mc D}_2$-patch. Then we verify that the resulting edge and corner system $\Sigma=E\cup V$ always implies one of three families in Reinhardt \cite{reinhardt}.

Actually the number of possible ${\mc D}_2$-patches is very large, and the brute force approach following the idea above is computationally intensive. We actually employ a more computationally efficient approach, by first considering the ${\mc D}_1$-patches and then examining the additional equalities that need to be satisfied along the boundary of the ${\mc D}_1$-patch. 

Our proof of Theorem \ref{thm2} follows the same approach. The ${\mc D}_k$-patches in the tilings in Theorem \ref{thm2} are combinatorially part of the uniform tiling $\{7,3\}$ of the Poincar\'e disk, described by the right of Figure \ref{2layers}. In fact, we need to consider ${\mc D}_3$-patches for some of the more difficult cases. See Table \ref{table3}.

We remark that classification for (convex) hexagonal tilings of the plane is similar to Theorem \ref{thm2}, {\em without} the balanced assumption. At the end, it turns out the corner part in all three Reinhardt families are balanced and come from hexagonal tilings of the torus \cite{ywy}. If we similarly remove the balanced condition in Theorem \ref{thm2}, however, we find a number of tilings of the Poincar\'e disk by congruent heptagons with unbalanced corner systems. These heptagons do not satisfy any one of the four systems in Theorem \ref{thm1}, and therefore cannot tile closed surfaces. The study of these tilings is the subject of another paper. 

Now we try to find all tilings of the Poincar\'e disk and hyperbolic surfaces by a given heptagon. A bigger edge and corner system allows more flexibility in fitting tiles together. If the system is too big, then the flexibility in the fitting would make it impossible to describe all the tilings. Therefore we are content with describing tilings by congruent heptagons in Theorem \ref{thm1}. 

\begin{theorem}\label{thm3}
Edge-to-edge tilings of the Poincar\'e disk by generic $7_{2+3}$- and $7_{3+3}$-heptagons are unique. Moreover, tilings of the Poincar\'e disk by the unique non-regular $7^1_{5+6}$- and $7^2_{5+6}$-heptagons are edge-to-edge and unique, with the vertices given by exactly \eqref{h11ev1} and \eqref{h11ev2}.
\end{theorem}

Figures \ref{H5B}, \ref{H6B}, \ref{H111B}, \ref{H112B} give respective schematic and real drawings of the unique tilings ${\mc T}(\Sigma)$ of ${\bb H}^2$ for the heptagons in the theorem, for $\Sigma=7_{2+3},7_{3+3},7^1_{5+6},7^2_{5+6}$. The unique tilings have the corresponding Fuchsian symmetry groups $G(\Sigma)$. The groups are transitive on all tiles for $\Sigma=7_{2+3},7_{3+3}$, and have three orbits for $\Sigma=7^1_{5+6},7^2_{5+6}$. The three orbits are represented by tiles in three colors in Figures \ref{H111B} and \ref{H112B}. 

Theorem \ref{thm3} leads to the classification of tilings of closed hyperbolic surfaces by the heptagons in the theorem. These exactly correspond to the conjugate classes of torsionless subgroups of finite indices in $G(\Sigma)$. For example, tilings of $2{\bb T}^2$ by generic $7_{2+3}$-heptagons are the quotients of ${\mc T}(7_{2+3})$ by the conjugate classes of torsionless subgroups $G_k(7_{2+3}, 2{\bb T}^2)$ of index 12 in $G(7_{2+3})$. We use the standard computer package GAP to calculate these conjugate classes. At the end of the paper, we present all tilings of the surfaces $2{\bb T}^2$ and $4{\bb P}^2$ of Euler characteristic $-2$ in Figures \ref{H5C}, \ref{H6C}, \ref{H111C}, \ref{H112C}. 

\begin{theorem}\label{thm5}
The numbers of tilings of hyperbolic surfaces of small genus by generic $7_{2+3}$- and $7_{3+3}$-heptagons, and by the non-regular $7^1_{5+6}$- and $7^2_{5+6}$-heptagons, are given by Table \ref{table5}. Moreover,
\begin{enumerate}
\item Generic $7_{2+3}$-heptagons and the non-regular $7^1_{5+6}$-heptagon can only tile orientable surfaces.
\item Generic $7_{3+3}$-heptagons and the non-regular $7^2_{5+6}$-heptagon cannot tile non-orientable surfaces of odd genus.
\end{enumerate}
\end{theorem}

\begin{table}[htp]
\renewcommand{\arraystretch}{1.2}
\centering
\begin{tabular}{|c||c|c|c|c|c|c|}
\hline
& minimal & covering &
minimal & covering &
minimal & covering \\
\hline
heptagon & 
\multicolumn{2}{c|}{$2{\bb T}^2$} &
\multicolumn{2}{c|}{$3{\bb T}^2$} &
\multicolumn{2}{c|}{$4{\bb T}^2$} \\
\hline \hline
$7_{2+3}$ & 
39 & 0 & 
2301 & 330 & 
490044 & 2601 \\
\hline
$7_{3+3}$ & 
5 & 0 & 
78 & 52 & 
7126 & 226  \\
\hline
$7^1_{5+6}$ & 
10 & 0 & 
130 & 70 & 
5780 & 600  \\
\hline
$7^2_{5+6}$ & 
1 & 0 & 
4 & 16 & 
60 & 53  \\
\hline \hline
heptagon & 
\multicolumn{2}{c|}{$4{\bb P}^2$} &
\multicolumn{2}{c|}{$6{\bb P}^2$} & 
\multicolumn{2}{c|}{$8{\bb P}^2$} \\
\hline\hline
$7_{3+3}$ & 
31 & 0 &
3278 & 298 &
1892964 & 1865 \\
\hline
$7^2_{5+6}$ & 
10 & 0 &
128 & 56 &
8440 & 436 \\
\hline
\end{tabular}
\caption{Numbers of tilings of surfaces by generic heptagons.}
\label{table5}
\end{table}

We remark that the number of tilings in Table \ref{table5} is counted as equivalence classes of combinatorial tilings. A {\em combinatorial tiling} subject to an edge and corner system $\Sigma=E\cup V$ has all corners of the tiles labeled, such that edge labels are matched according to $E$, and corner labels form vertices according to $V$. Suppose ${\mc T}$ and ${\mc T}'$ are combinatorial tilings of closed surfaces $S$ and $S'$, subject to edge and corner systems $\Sigma$ and $\Sigma'$. The two tilings are {\em equivalent}, if there is a homeomorphism $\Lambda$ between $S$ and $S'$, and a rotation or flip $\lambda$ of corner labels that gives equivalence between $\Sigma$ and $\Sigma'$, such that $\Lambda$ takes tiles of ${\mc T}$ to tiles of ${\mc T}'$, and takes labels $i$ of the corners in ${\mc T}$ to the labels $\lambda(i)$ of the corresponding corners in ${\mc T}'$. The number of tilings in the table is the number of equivalence classes of the combinatorial tilings.

There may not be a satisfying geometrical interpretation of the number of tilings. The first problem is that the underlying surfaces are homeomorphic but may not be isometric. Therefore the comparison of tilings need to go through certain deformations of the metric and the heptagon. Second, the deformation between heptagons that are generic with respect to $\Sigma$ may have to go through some heptagons that are not generic with respect to $\Sigma$.

We remark that the Fuchsian group $G(\Sigma)$ is actually determined by the combinatorial information from $\Sigma$. The conjugate classes of torsionless subgroups of finite indices in $G(\Sigma)$ is still combinatorial information. This is the true reason for the counting of tilings to be combinatorial.

In Table \ref{table5}, we also distinguish minimal and covering tilings. This arises from the consideration that any covering of a tiling of a hyperbolic surface is a tiling of another hyperbolic surface by the same prototile and with the same edge and corner system. For example, various two fold covers of 39 tilings of $2{\bb T}^2$ by generic $7_{2+3}$-heptagon give 330 tilings of $3{\bb T}^2$. This leads to the following concept.

\begin{definition}\label{minimal}
A tiling of a surface is {\em minimal} if it is not the covering of a tiling of a surface with strictly smaller number of tiles.
\end{definition}

Table \ref{table5} suggests most tilings are minimal. This is confirmed in a different setting by Sun and Yan \cite{hy1}. Therefore listing all minimal tilings is not much easier than listing all tilings.

Table \ref{table5} also suggests the numbers of tilings for orientable and non-orientable tilings of the same Euler characteristic have similar scale. 

We conclude the introduction by some general remarks. 

First, in Theorem \ref{thm5}, we have no tilings for $g{\bb P}^2$ with odd $g$. Therefore there should be another notion of generic heptagons for these surfaces, and Theorem \ref{thm1} should be amended by another theorem that classifies all heptagons that can tile $g{\bb P}^2$ for all odd $g$. These families are obtained by adding suitable equalities to the families in Theorem \ref{thm1}.

Second, the symmetries on tilings by the unique non-regular $7^1_{5+6}$- and $7^2_{5+6}$-heptagons have three orbits. The fundamental domain of the symmetry can be regarded as a nonagon, and we obtain tilings by congruent nonagons. In fact, Figures \ref{H111B} and \ref{H112B} present only one choice of the three color fundamental domain. There are many other possible choices. All of these lead to families of nonagons that can tile hyperbolic surfaces. The investigation of this phenomenon appears in a separate paper \cite{lwwy3}.  

Third, seven is only the minimal number of sides for the problem of tilings by maximal number of tiles. This paper gives a systematic method that can be applied for octagonal tilings, nonagonal tilings, etc. If we restrict the degree of vertices to be $\ge 4$, then we can prove similar finiteness theorem for $n\ge 5$, and maximal number of tilings will be equivalent to all vertices having degree 4 in this context. The systematic method can also be applied to this problem. 

Finally, the programs used in the classification work can be found in https://github.com/wwuj/Heptagon-Tiling-Solver. We also use the computer package GAP for calculating tilings of specific surfaces.

\section{Edge and Corner Systems for Maximal Number of Tiles}

In Section \ref{finite}, we prove that tilings of closed hyperbolic surfaces by maximal number of tiles is edge-to-edge, and all vertices have degree $3$. 

The way to prove Theorem \ref{thm2} and then Theorem \ref{thm1} is to first prescribe all the possible edge and corner systems, and then construct tilings according to these systems. In Section \ref{edge}, we list all the possible edge systems. In Section \ref{corner}, we list all the possible corner systems. In Section \ref{edge-corner}, we combine the two and lists all the possible edge and corner systems.  

\subsection{Finite Problem}
\label{finite}

\begin{lemma}\label{bound}
In a tiling of a closed surface of Euler characteristic $\chi<0$ by $n$-gons with $n\ge 7$, if the degrees of full vertices are $\ge 3$, and the degrees of half vertices are $\ge 2$,  then the number $f$ of tiles is bounded
\[
\frac{-2\chi}{n-2}<f\le \frac{-6\chi}{n-6}.
\]
Moreover, $f$ equals the maximal number $\frac{-6\chi}{n-6}$ if and only if the tiling is edge-to-edge and all vertices have degree $3$.
\end{lemma}

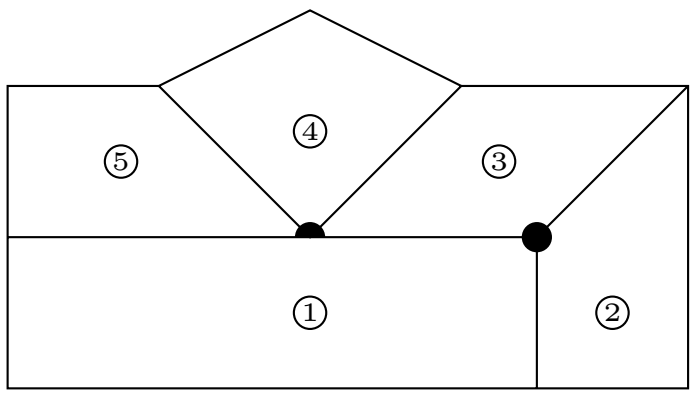
\begin{figure}[htp]
\centering
\begin{tikzpicture}[>=latex,scale=1]

\draw
	(0,0) -- (3.5,0) -- (4.5,1)
	(3,1) -- (2,1.5) -- (1,1) -- (0,1) -- (0,-1) -- (4.5,-1) -- (4.5,1) -- (3,1) -- (2,0) -- (1,1)
	(3.5,0) -- (3.5,-1);

\fill
	(3.5,0) circle (0.1)
	(2.1,0) arc (0:180:0.1);

\foreach \x/\y/\a in {2/-0.5/1, 4/-0.5/2, 3.25/0.5/3, 2/0.7/4, 0.75/0.5/5}
\node[inner sep=0.5, draw, shape=circle] at (\x, \y) {\tiny \a};

\end{tikzpicture}
\caption{Full vertex and half vertex.}
\label{terms}
\end{figure}

The lemma includes non-edge-to-edge tilings. Figure \ref{terms} shows a non-edge-to-edge tiling by five quadrilaterals, with a {\em full vertex} $\bullet$ and a {\em half vertex} $\hbullet$. The degree of a vertex is the number of corners at the vertex. The full vertex $\bullet$ has corners from \circled{1}, \circled{2}, \circled{3}, and the degree is 3. The half vertex $\hbullet$ has corners from \circled{3}, \circled{4}, \circled{5}, and the degree is also 3. The tile \circled{1} is also at the half vertex, but its corner are not at the half vertex. Therefore \circled{1} is not counted towards the degree of the half vertex. A tiling is edge-to-edge if and only if all vertices are full.

The condition on the degrees of full and half vertices in the lemma means that all vertices are {\em locally} the meeting place of at least three tiles. Here we emphasize the condition is local because different corners at a vertex may actually be joined by a long path within a large tile. This is exactly what happens in single tile tilings.

\begin{proof}
Let ${\mc T}$ be a tiling in the lemma. Let $v$, $h$, $e$, $f$ be the numbers of full vertices, half vertices, edges and tiles in ${\mc T}$. Let $v_d$ and $h_d$ be the numbers of full and half vertices of degree $d$. We have $d\ge 3$ in $v_d$, and $d\ge 2$ in $h_d$, and
\[
v=v_3+v_4+v_5+\cdots,\quad
h=h_2+h_3+h_4+\cdots.
\]
Moreover, we have the Euler formula
\[
v-e+f=\chi.
\]

Let ${\mc T}_E$ be the same tiling ${\mc T}$, but with all half vertices regarded as full vertices. Then ${\mc T}_E$ is an edge-to-edge tiling. A full vertex of degree $d$ in ${\mc T}$ becomes a vertex of degree $d$ in ${\mc T}_E$. A half vertex of degree $d$ in ${\mc T}$ becomes a vertex of degree $d+1$ in ${\mc T}_E$. Since the number of edges of ${\mc T}_E$ is $e+h$, we get
\begin{equation}\label{eqa}
2(e+h)
=3v_3+4v_4+5v_5+\cdots+3h_2+4h_3+5h_4+\cdots
=\sum_{d\ge 3}dv_d+\sum_{d\ge 2}(d+1)h_d, 
\end{equation} 
On the other hand, if $h(t)$ is the number of half vertices on the boundary of a tile $t$ of ${\mc T}$, then $t$ is an $(n+h(t))$-gon in ${\mc T}_E$. Therefore we get
\begin{equation}\label{eqb}
2(e+h)
=\sum_{t\in {\mc T}}(n+h(t))
=nf+\sum_{t\in {\mc T}}h(t)=nf+h. 
\end{equation}
The last equality in is due to the fact that each half vertex is on the side of exactly one tile. 

By \eqref{eqb} and the Euler equality, we get
\begin{equation}\label{eq1}
-2\chi
=-2v+2e-2f
=-(2v+h)+(n-2)f
<(n-2)f.
\end{equation}
For $n\ge 3$, this implies the inequality $f>\frac{-2\chi}{n-2}$ in the lemma.

By \eqref{eqa} and \eqref{eqb}, we get
\begin{equation}\label{eq2}
nf
=\sum_{d\ge 3}dv_d+\sum_{d\ge 2}(d+1)h_d-h
=\sum_{d\ge 3}dv_d+\sum_{d\ge 2}dh_d.
\end{equation}
Then by \eqref{eq1} and \eqref{eq2}, we get
\begin{align*}
6\chi+(n-6)f
&=3((2v+h)-(n-2)f)+(n-6)f
=3(2v+h)-2nf \\
&=3\left(2\sum_{d\ge 3}v_d+\sum_{d\ge 2}h_d\right)
-2\left(\sum_{d\ge 3}dv_d+\sum_{d\ge 2}dh_d\right) \\
&=-2\sum_{d\ge 3}(d-3)v_d-\sum_{d\ge 2}(2d-3)h_d
\le 0.
\end{align*}
For $n\ge 7$, the inequality is the same as $f\le \frac{-6\chi}{n-6}$ in the lemma. Moreover, $f=\frac{-6\chi}{n-6}$ if and only if the inequality becomes equality, which means $v_d=0$ for $d\ge 4$ and $h_d=0$ for all $d\ge 2$. 
\end{proof}

A consequence of Lemma \ref{bound} is $\frac{-6\chi}{n-6}\ge f\ge 1$, which is
\[
n\le 6(1-\chi).
\]
Combined with Lemma \ref{bound}, we conclude that, the tilings of a fixed surface by $n$-gons, with $n\ge 7$, happens only for finitely many $n$ and $f$. 

From now on, we concentrate on the case $n=7$ and $f=\frac{-6\chi}{7-6}=-6\chi$. Then the tiling is edge-to-edge, and an edge equality $|\bar{i}|=|\bar{j}|$ means the $i$-th and $j$-th edges can be glued together in the tiling. Moreover, all vertices have degree 3, and an angle equality is of the form $[i]+[j]+[k]=2\pi$, and a vertex is of the form $ijk$. 

\subsection{Edge System}
\label{edge}

\begin{lemma}\label{lcondition}
In a tiling by congruent polygons, suppose all vertices have degree $3$, and $u,a,b,v$ are consecutive edges. If the length of each of $a,b$ is different from all the other edges, then $u$ and $v$ have the same length.
\end{lemma}

\begin{figure}[htp]
\centering
\begin{tikzpicture}[>=latex, scale=0.8, font=\footnotesize]

\foreach \x in {1,-1}
{
\begin{scope}[xscale=\x]

\draw
	(90:1) -- (0,0) -- (-30:1) -- ++(-60:1);

\draw[dashed]
	(1.36,-1.36) to[out=-90, in=0] (0,-2)
	(90:1) -- ++(30:1) to[out=30, in=30] (1.4,-0.5) -- (-30:1);

\end{scope}
}

\node at (0.95,-1) {$u$};
\node at (-0.95,-1) {$v$};
\node at (0.37,-0.4) {$a$};
\node at (-0.37,-0.4) {$b$};
\node at (0.2,0.4) {$x$};

\foreach \x/\y/\a in {0/-1.3/1, 1/0.5/2, -1/0.5/3}
\node[inner sep=0.5, draw, shape=circle] at (\x,\y) {\a};
				
\end{tikzpicture}
\caption{Constraint on edge lengths.} 
\label{length}
\end{figure}
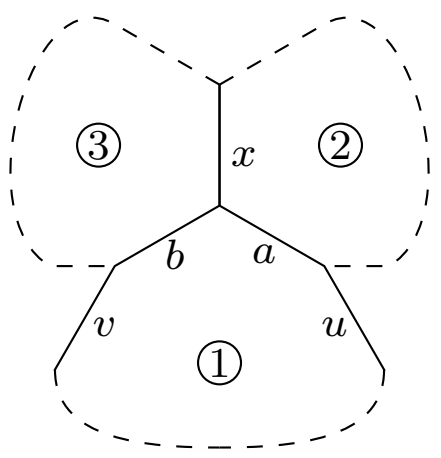

\begin{proof}
Consider the consecutive edges $u,a,b,v$ of tile \circled{1} in Figure \ref{length}. We have two other tiles \circled{2} and \circled{3} at the degree 3 vertex connecting $a$ and $b$. The two tiles share an edge $x$. Since $x$ and $a$ are adjacent in the tile \circled{2}, and $a$ is different from all the other edge lengths in the polygon, we know $x=b$ or $u$. By applying the same argument to $x$ and $b$ in \circled{3}, we also get $x=a$ or $v$. Since each of $a,b$ are different from all the other edges, we conclude $x=u$ and $x=v$. Therefore $u=v$. 
\end{proof}

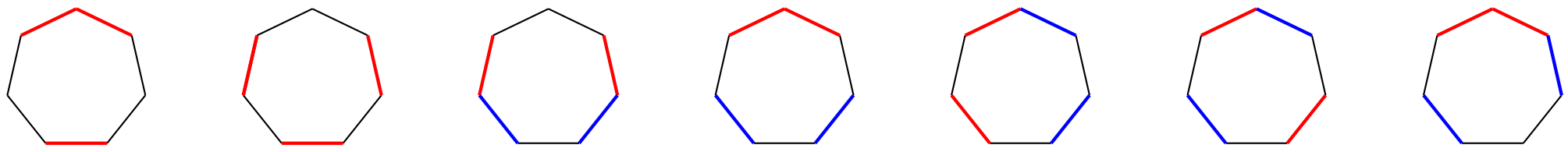
\begin{figure}[htp]
\centering
\begin{tikzpicture}[>=latex, scale=1, font=\tiny]

\foreach \a/\x in 
	{1/0, 2/0, 4/0, 5/0,
	0/1, 2/1, 4/1, 6/1,
	0/2, 3/2, 6/2,
	1/3, 3/3, 5/3,
	1/4, 3/4, 5/4,
	1/5, 3/5, 5/5,
    1/6, 3/6, 4/6}
\draw[xshift=2*\x cm]
	(90+51.4286*\a:0.6) -- (90+51.4286+51.4286*\a:0.6);

\foreach \a/\x in 
	{0/0, 3/0, 6/0,
	1/1, 1/1, 3/1, 5/1,
	1/2, 5/2,
	0/3, 6/3,
	0/4, 2/4,
	0/5, 4/5,
    0/6, 6/6}
\draw[red, xshift=2*\x cm, thick]
	(90+51.4286*\a:0.6) -- (90+51.4286+51.4286*\a:0.6);

\foreach \a/\x in 
	{2/2, 4/2,
	2/3, 4/3,
	4/4, 6/4,
	2/5, 6/5,
    2/6, 5/6}
\draw[blue, xshift=2*\x cm, thick]
	(90+51.4286*\a:0.6) -- (90+51.4286+51.4286*\a:0.6);
					
\end{tikzpicture}
\caption{All the rank 2 edge systems that satisfy Lemma \ref{lcondition}. Red edges have equal length, and blue edges have equal length, and these are the only equal lengths.} 
\label{length2}
\end{figure}

A consequence of the lemma is that the maximal number of distinct edge lengths for heptagonal tiling is five. Therefore $\text{rank}E\ge 2$. Moreover, the maximal rank of $E$ is 6, because this means all edge lengths are equal, and the tiling is equilateral. Figure \ref{length2} shows all the possible edge systems $E$ of rank 2. 

Table \ref{table1} lists the number of heptagonal edge systems satisfying the lemma.

\begin{table}[htp]
	\centering
	\begin{tabular}{|c||c|c|c|c|c|}
	\hline
	rank & 2 & 3 & 4 & 5 & 6  \\
	\hline
	number of $E$ & 7 & 29 & 31 & 8 & 1 \\
	\hline
	\end{tabular}
\caption{Numbers of possible heptagonal edge systems $E$ satisfying Lemma \ref{lcondition}. }
\label{table1}
\end{table}

\subsection{Corner System}
\label{corner}

By Lemma \ref{bound}, in tilings by maximal number of $f=-6\chi$ of heptagons, all vertices have degree 3. As explained in the overview, we only consider complete corner systems, for which we may use the following notations for (degree 3) vertices interchangeably:
\begin{itemize}
\item A vertex $ijk$ consists of corners $i,j,k$. For example, $01^2$ means the vertex consists of one corner 0 and two corners 1.
\item A vertex is a vector in ${\bb R}^7$ with non-negative integers as the coordinates, with the $i$-th coordinate being the number of corner $i$ at the vertex. For example, $(1,2,0,0,0,0,0)$ means the vertex $01^2$. 
\item A vertex is an angle equality $[i]+[j]+[k]=2\pi$. For example, the equation $[0]+[1]+[1]=2\pi$ means the vertex $01^2$.
\end{itemize}
We denote by $V({\mc T})$ all the vertices in a tiling ${\mc T}$, and call it the corner system of ${\mc T}$. 

For heptagonal tiling and in vector form, the degree 3 requirement means that all $\vec{v}\in V$ satisfy
\[
\vec{v}\cdot \vec{1}=3,\text{ where }
\vec{1}=(1,1,1,1,1,1,1)\in {\bb R}^7.
\]
We also have $\text{rank}V\le 7$. Since all angle values being $\frac{2\pi}{3}$ always satisfy the equation version of $V$, we know $\text{rank}V=7$ means this is the unique solution. On the other hand, since all seven corners should appear in the tiling, we have $\text{rank}V\ge 3$.

The possible corner systems for tilings in this paper should consist of only degree 3 vertices and should be complete. We remark that, if all vertices in $V$ have degree 3, then any angle equality $[i_1]+[i_2]+\dots+[i_d]=\theta$ derived as a linear combination of the angle equalities in $V$ will have $\theta=\frac{2d\pi}{3}$ on the right side. Therefore all the vertices in the completion $\bar{V}$ also have degree $3$.

Sometimes a complete corner system is a very long list. Then we may present a corner system as a collection of angle equalities, and the complete corner system consists of all the equalities $[i]+[j]+[k]=2\pi$ that are linear combinations of the collection.

We list all the possible corner systems by their ranks. We clearly have $\text{rank}\bar{V}=\text{rank}V$. A corner system of rank $r$ is {\em minimal} if it consists of exactly $r$ vertices. Therefore vertices in a minimal corner system are linearly independent. We only need to list minimal corner systems, with the understanding that any corner system of rank $r$ is between a minimal corner system $V$ of rank $r$ and its completion $\bar{V}$. Minimal corner systems of rank 3 are obtained by applying all the permutations of ${\bb Z}_7$ to  the following
\begin{equation}\label{vs1}
\begin{split}
\{012,345,6^3\},\;
\{012,345,06^2\},\;
\{012,345,0^26\},\;
\{012,234,56^2\},\;\\
\{012,34^2,56^2\},\;
\{012,345,016\},\;
\{012,345,036\},\;
\{012,036,045\}.
\end{split}
\end{equation}

Two corner systems $V_1$ and $V_2$ are {\em equivalent} if they imply each other through linear algebra. This is equivalent to $\bar{V}_1=\bar{V}_2$, and can be verified by $\text{rank}V_1=\text{rank}(V_1\cup V_2)=\text{rank}V_2$. Our prescribed corner systems will only be presented as equivalent classes of minimal corner systems.

Besides requiring that all corners appear, the possible corner systems should also satisfy a balance condition that comes from counting the total number of each corner. Suppose $\vec{v}\in V({\mc T})$ appears $x_{\vec{v}}$ times in a tiling ${\mc T}$ of a closed surface, and $f$ is the number of tiles. Then each corner $i$ appears $f$ times in the tiling:
\[
f=\sum_{\vec{v}\in V({\mc T})}x_{\vec{v}}v_i,\quad
i\in {\bb Z}_7,
\]
where $v_i$ is the $i$-th coordinate of $\vec{v}$, i.e., the number of times the corner $i$ appears at $\vec{v}$. The seven equalities together is the {\em vertex counting equation}
\begin{equation}\label{vceq}
f\vec{1}=\sum_{\vec{v}\in V({\mc T})}x_{\vec{v}}\vec{v}.
\end{equation}
Since $x_{\vec{v}}$ are non-negative integers, after dividing $f$, we get the following condition, which already appeared in \cite{jlwy,ly1,rao}.

\begin{definition}\label{balance}
A corner system $V$ is {\em balanced} if the unit vector $\vec{1}$ is a linear combination of vectors in $V$, with non-negative coefficients.
\end{definition}

The corner system $\{012,345,6^3\}$ in \eqref{vs1} is balanced because
\[
(1,1,1,1,1,1,1)=1(1,1,1,0,0,0,0)+1(0,0,0,1,1,1,0)+\tfrac{1}{3}(0,0,0,0,0,0,3).
\]
The balance condition for the corner system $\{012,345,06^2\}$ means
\[
(1,1,1,1,1,1,1)
=x_1(1,1,1,0,0,0,0)
+x_2(0,0,0,1,1,1,0)
+x_3(1,0,0,0,0,0,2).
\]
We find that the 0th, 1st and 6th coordinate equations are contradictory. The equation does not even have real number solution. In fact, $\{012,345,6^3\}$ is the only balanced corner system in \eqref{vs1}. 

If $V$ is balanced, and $V\sub V'$, then $V'$ is also balanced. Therefore in prescribing corner systems, we only consider complete and balanced corner systems. Moreover, we have $\sigma(\bar{V})=\overline{\sigma(V)}$ for any permutation $\sigma$ of ${\bb Z}_7$. Table \ref{table2} gives the numbers of equivalent classes of complete and balanced heptagonal corner systems, up to the permutations of ${\bb Z}_7$. 

\begin{table}[htp]
	\centering
	\begin{tabular}{|c||c|c|c|c|c|}
	\hline
	rank & 3 & 4 & 5 & 6 & 7  \\
	\hline
	number of $V$ & 1 & 5 & 23 & 42 & 1 \\
	\hline
	\end{tabular}
\caption{Numbers of equivalence classes of complete and balanced heptagonal corner systems, up to permutations of corner labels. }
\label{table2}
\end{table}

\subsection{Edge and Corner System}
\label{edge-corner}

Now we combine an edge system in Table \ref{table1} and a corner system in Table \ref{table2}. Since the tables list the numbers of systems up to equivalences, we need to pick one representative $E$ and one representative $V$ of the equivalence classes, and then consider all $\lambda(E)\cup\sigma(V)$. Here $\lambda$ is the symmetry in assigning circularly ordered labels to the corners of heptagon, which is either the rotation $\lambda(i)=c+i$ or the flip $\lambda(i)=c-i$. The corresponding map on the edges is $\lambda(\bar{i})=\overline{c+i}$ or $\lambda(\bar{i})=\overline{c-i-1}$. Moreover, $\sigma$ is a permutation of ${\bb Z}_7$. 

For example, if we assign 0 to the corner bounded by two red edges in the first edge system in Figure \ref{length2}, and assign the other circularly ordered labels in counterclockwise orientation, then the edge system is
\[
E\colon |\bar{0}|=|\bar{3}|=|\bar{6}|.
\]
Applying the rotation $\lambda(i)=2+i\colon {\bb Z}_7\to {\bb Z}_7$ to the edge system gives
\[
\lambda(E)\colon |\bar{2}|=|\bar{1}|=|\bar{5}|.
\]
On the other hand, applying the permutation $\sigma\colon 0,1,2,3,4,5,6\to 3,1,5,2,0,6,4$ to $V=\{012,345,6^3\}$, we get the corner system
\[
\sigma(V)\colon [1]+[3]+[5]=[0]+[2]+[6]=[4]+[4]+[4]=2\pi.
\]
Then $\lambda(E)\cup \sigma(V)$ is the following system
\begin{equation}\label{eq3}
\lambda(E)\cup \sigma(V)\colon 
|\bar{2}|=|\bar{1}|=|\bar{5}|,\;
[1]+[3]+[5]=[0]+[2]+[6]=[4]+[4]+[4]=2\pi.
\end{equation}

After forming the system $\lambda(E)\cup \sigma(V)$, we need to construct tilings ${\mc T}$, such that $E({\mc T})\sub \lambda(E)$ and $V({\mc T})\sub \sigma(\bar{V})=\overline{\sigma(V)}$. As far as tilings are concerned, two systems $E_1\cup V_1$ and $E_2\cup V_2$ are {\em equivalent} if they are related by the symmetry of the circular label assignment. In other words, there is a rotation or flip $\lambda$, such that $E_2=\lambda(E_1)$ and $V_2=\lambda(V_1)$. Since $\lambda(E)\cup\sigma(V)$ is equivalent to $E\cup \lambda^{-1}\sigma(V)$, we may fix the representative $E$ and only consider $E\cup \sigma(V)$ for all the permutations $\sigma$. Then $E\cup \sigma(V)$ and $E\cup \sigma'(V)$ are equivalent if $\sigma'(V)=\lambda\sigma(V)$ for some $\lambda$ satisfying $\lambda(E)=E$. 

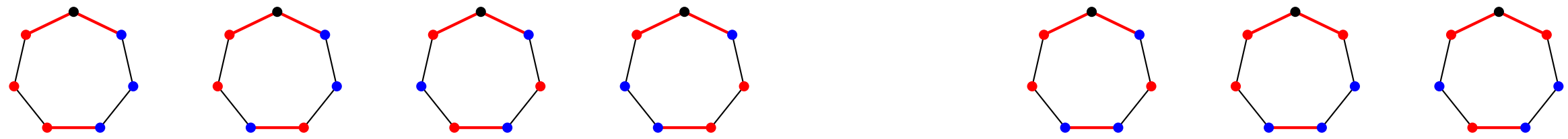
\begin{figure}[htp]
\centering
\begin{tikzpicture}[>=latex, scale=1, font=\tiny]

\foreach \x in {0,1,2,3,5,6,7}
{
\begin{scope}[xshift=2*\x cm]

\foreach \a in {1,2,4,5}
\draw
	(90+51.4286*\a:0.6) -- (90+51.4286+51.4286*\a:0.6);

\foreach \a in {0,3,6}
\draw[red, thick]
	(90+51.4286*\a:0.6) -- (90+51.4286+51.4286*\a:0.6);

\end{scope}
}

\foreach \x in {0,1,2,3,5,6,7}
\fill
	(2*\x,0.6) circle (0.05);

\foreach \a/\x in 
	{1/0, 2/0, 3/0,
	1/1, 2/1, 4/1,
	1/2, 3/2, 5/2,
	1/3, 4/3, 5/3,
	1/5, 2/5, 5/5,
	1/6, 2/6, 6/6,
	1/7, 3/7, 6/7
	}
\fill[red, xshift=2*\x cm]
	(90+51.4286*\a:0.6) circle (0.05);

\foreach \a/\x in 
	{4/0, 5/0, 6/0,
	3/1, 5/1, 6/1,
	2/2, 4/2, 6/2,
	2/3, 3/3, 6/3,
	3/5, 4/5, 6/5,
	3/6, 4/6, 5/6,
	2/7, 4/7, 5/7}
\fill[blue, xshift=2*\x cm]
	(90+51.4286*\a:0.6) circle (0.05);
						
\end{tikzpicture}
\caption{Some edge and corner systems of rank 5, where the edge system is the first of Figure \ref{length2}. The edge system is indicated by three red edges of equal length. The corner system is $\{{\color{red} {\bullet}{\bullet}{\bullet}}, {\color{blue} {\bullet}{\bullet}{\bullet}}, \bullet^3\}$.} 
\label{landa}
\end{figure}

Let us consider the case $\text{rank}(E\cup V)=5$. This means $\text{rank}E=2$ and $\text{rank}V=3$. One example is $\lambda(E)\cup \sigma(V)$ in \eqref{eq3}, which is equivalent to 
\[
E\cup \lambda^{-1}\sigma(V)\colon 
|\bar{0}|=|\bar{3}|=|\bar{6}|,\;
[1]+[3]+[6]=[0]+[4]+[5]=[2]+[2]+[2]=2\pi.
\]
We fix the $E$ part and consider the other possibilities $E\cup \sigma(V)$. The only nontrivial $\lambda$ preserving $E$ is the flip $\lambda(i)=0-i=-i$. For the vertex sets of rank 3, we have only one choice of the representative $V=\{012,345,6^3\}$. The following are all $\sigma(V)$ satisfying $\lambda\sigma(V)=\sigma(V)$:
\[
\{123,456,0^3\},\;
\{124,356,0^3\},\;
\{135,246,0^3\},\;
\{145,236,0^3\}.
\]
These are the four edge and corner systems on the left of Figure \ref{landa}.

Any $\sigma(V)$ is of the form $\{i_1i_2i_3,i_4i_5i_6,i_0^3\}$, where $i_0,i_1,\dots,i_6$ is a permutation of $0,1,\dots,6$. The following are all $\sigma(V)$ satisfying $i_0=0$ and $\lambda\sigma(V)\ne \sigma(V)$ 
\[
\{125,346,0^3\},\;
\{126,345,0^3\},\;
\{136,245,0^3\};\quad
\{652,431,0^3\},\;
\{651,432,0^3\},\;
\{641,532,0^3\}.
\]
Since the later three are obtained by applying $\lambda(i)=-i$ to the first three, we only need to consider $E\cup \sigma(V)$ for the first three $\sigma(V)$. These are the three edge and corner systems on the right of Figure \ref{landa}.

If $i_0\ne 0$, then any $\sigma(V)=\{i_1i_2i_3,i_4i_5i_6,i_0^3\}$ satisfies $\lambda\sigma(V)\ne \sigma(V)$. Considering the exchange of $i_1i_2i_3$ and $i_4i_5i_6$, for each $i_0\ne 0$, the total number of such $\sigma(V)$ is half of the number of choices of $i_1i_2i_3$. Therefore the number of these $\sigma(V)$ is $\frac{1}{2}\cdot 6\cdot \binom{6}{3}=60$. Since $\lambda\sigma(V)\ne \sigma(V)$, the number of edge and corner systems that we need to consider for the tiling purpose is $\frac{1}{2}\cdot 60=30$.

If $E$ is the first edge system in Figure \ref{length2}, then as far as finding tilings is concerned, the total number of $E\cup V$ of rank $5$ that we need to consider is $4+3+30=37$. 

The fourth column of Table \ref{table3} gives the numbers of equivalence classes of all possible edge and corner systems. The number 37 is part of the number 292 in the table.

\begin{table}[htp]
	\centering
	\begin{tabular}{|c|c|c||c|c||c|c|c|c|c|c|}
	\hline
	\multirow{2}{*}{rank} & \multirow{2}{*}{rank$E$} & \multirow{2}{*}{rank$V$} & \multicolumn{2}{|c||}{$E\cup V$} & \multicolumn{6}{|c|}{${\mc D}$ patch} \\
	\cline{4-11}
	 &&& initial & updated & ${\mc D}_1$ & $\pa{\mc D}_1$ & ${\mc D}_2$ & $\pa{\mc D}_2$ & ${\mc D}_3$ & $\pa{\mc D}_3$ \\
	\hline
	\hline
	5 & 2 & 3 & 292 & & 1 & \multicolumn{5}{|c|}{$7_{2+3}$} \\
	\hline
	\hline
	\multirow{2}{*}{6} & 3 & 3 & 1535 & 1525 & 1 & \multicolumn{5}{|c|}{$7_{3+3}$}   \\
	\cline{2-11}
	 & 2 & 4 & 11421 & 11373 & 0 &  \multicolumn{5}{|c|}{}  \\
	\hline
	\hline
	\multirow{3}{*}{7} & 4 & 3 & 1543 & 1512 & 0 &  \multicolumn{5}{|c|}{}  \\
	\cline{2-6}
	 & 3 & 4 & 61155 & 60628 & 0  & \multicolumn{5}{|c|}{}  \\
	\cline{2-6}
	 & 2 & 5 & 123923 & 123321 & 0 & \multicolumn{5}{|c|}{}  \\
	\cline{1-6}
    \multirow{4}{*}{8} & 5 & 3 & 329 & 307 & 0 &  \multicolumn{5}{|c|}{}  \\
	\cline{2-7}
	 & 4 & 4 & 61209 & 59733 & 3 & 0 & \multicolumn{4}{|c|}{}  \\
	\cline{2-7}
	 & 3 & 5 & 665215 & 658639 & 1 & 0 & \multicolumn{4}{|c|}{}  \\
     \cline{2-7}
	 & 2 & 6 & 381312 & 379590 & 0 &  \multicolumn{5}{|c|}{} \\
	\cline{1-6}
    \multirow{5}{*}{9} & 6 & 3 & 7 & 5 & 0 & \multicolumn{5}{|c|}{} \\
	\cline{2-7}
	 & 5 & 4 & 12852 & 11813 & 11 & 0 & \multicolumn{4}{|c|}{} \\
	\cline{2-8}
	 & 4 & 5 & 665417 & 647098 & 95 & 13 & 0 & \multicolumn{3}{|c|}{} \\
     \cline{2-11}
	 & 3 & 6 & 2048724 & 2029979 & 7 & 4 & \multicolumn{4}{|c|}{rank 10, special case of $7_{2+3}$}\\
    \cline{2-11}
	 & 2 & 7 & 7 & 6 & 0 & \multicolumn{5}{|c|}{} \\
	\hline
	\hline
    \multirow{4}{*}{10} & 6 & 4 & 216 & 123 & 8 & 0 & \multicolumn{4}{|c|}{} \\
	\cline{2-9}
	 & 5 & 5 & 139426 & 126703 & 375 & 35 & 0 & 0  & \multicolumn{2}{|c|}{} \\
	\cline{2-10}
	 & 4 & 6 & 2048916 & 1996977 & 674 & 223 & 1 & 1 & 0 & \\
     \cline{2-10}
	 & 3 & 7 & 29 & 18 & 1 & 1 & 0 & \multicolumn{3}{|c|}{} \\
	\cline{1-9}
    \multirow{3}{*}{11} & 6 & 5 & 2258 & 1174 & 249 & 19 & 1 & 0 & \multicolumn{2}{|c|}{} \\
	\cline{2-11}
	 & 5 & 6 & 428988 & 393333 & 3276 & 810 & 20 & 16 & 6 & 6 (4 regular)\\
	\cline{2-11}
	 & 4 & 7 & 31 & 11 & 1 & \multicolumn{5}{|c|}{regular} \\
	\hline
    \multirow{2}{*}{12} & 6 & 6 & 6852 & 3930 & 2214  & 414 & 86 & 44 & 32  &  31 (31 regular) \\
	\cline{2-11}
	 & 5 & 7 & 8 & 1 & 0 & \multicolumn{5}{|c|}{} \\
	\hline
	\end{tabular}
\caption{Numbers of possible edge and angle systems. No tiling if the rightmost number is $0$. }
\label{table3}
\end{table}

\section{Tilings of Rank 5}
\label{rank5}

\begin{figure}[htp]
\centering
\begin{tikzpicture}[>=latex,scale=1, font=\tiny]

\foreach \a in {0,...,6}
\fill[gray!50, rotate=51.43*\a]
	(0,0.5) -- (0,0.75) arc (90:208:0.25);

\foreach \a in {0,...,6}
\draw[rotate=51.43*\a]
	(78:1.2) -- (90:0.9) -- (90:0.5) -- (38.57:0.5) -- (38.57:0.9) -- (50.57:1.2) -- (64.286:1.3) -- (78:1.2)
	(78:1.2) -- (82:1.5) -- (87:1.7) -- (93:1.7) -- (98:1.5) -- (102:1.2)
	(82:1.5) -- (78:1.75) -- (74:1.9) -- (68.286:1.85) -- (64.286:1.6) -- (64.286:1.3)
	(98:1.5) -- (102:1.75) -- (106:1.9) -- (111.714:1.85) -- (115.714:1.6);

\foreach \a in {0,...,6}
\foreach \b in {0,...,6}
{
\node at (90+51.43*\a:0.35) {\a};

\coordinate (A\a) at (90+51.43*\a:0.5);
\coordinate (B\a) at (90+51.43*\a:0.9);
\coordinate (C\a) at (102+51.43*\a:1.2);
\coordinate (D\a) at (129.43+51.43*\a:1.2);
\coordinate (E\a) at (115.714+51.43*\a:1.3);
\coordinate (F\a) at (115.714+51.43*\a:1.6);

\coordinate (I\a) at (82+51.43*\a:1.5);
\coordinate (J\a) at (98+51.43*\a:1.5);
\coordinate (K\a) at (87+51.43*\a:1.7);
\coordinate (L\a) at (93+51.43*\a:1.7);

\coordinate (U\a) at (68.286+51.43*\a:1.85);
\coordinate (V\a) at (74+51.43*\a:1.9);
\coordinate (W\a) at (78+51.43*\a:1.75);

\coordinate (P\a) at (111.714+51.43*\a:1.85);
\coordinate (Q\a) at (106+51.43*\a:1.9);
\coordinate (R\a) at (102+51.43*\a:1.75);

\coordinate (X\a 0) at (100+51.43*\a:0.6);
\coordinate (X\a 1) at (99+51.43*\a:0.87);
\coordinate (X\a 2) at (105+51.43*\a:1.1);
\coordinate (X\a 3) at (115.714+51.43*\a:1.17);
\coordinate (X\a 4) at (126.43+51.43*\a:1.1);
\coordinate (X\a 5) at (132.43+51.43*\a:0.87);
\coordinate (X\a 6) at (131.43+51.43*\a:0.6);

\coordinate (Y\a 0) at (90+51.43*\a:1.05);
\coordinate (Y\a 1) at (83+51.43*\a:1.2);
\coordinate (Y\a 2) at (85.5+51.43*\a:1.4);
\coordinate (Y\a 5) at (94.5+51.43*\a:1.4);
\coordinate (Y\a 6) at (97+51.43*\a:1.2);

\coordinate (U\a 0) at (75.5+51.43*\a:1.32);
\coordinate (U\a 1) at (68.4+51.43*\a:1.37);
\coordinate (U\a 2) at (68.4+51.43*\a:1.55);
\coordinate (U\a 6) at (77+51.43*\a:1.5);

\coordinate (V\a 0) at (104.5+51.43*\a:1.32);
\coordinate (V\a 1) at (103+51.43*\a:1.5);
\coordinate (V\a 5) at (111.6+51.43*\a:1.55);
\coordinate (V\a 6) at at (111.6+51.43*\a:1.37);
}

\draw[red, thick]
	(A6) -- (A0) -- (A1)
	(A0) -- (B0)
	(A2) -- (B2) -- (D1)
	(B2) -- (C2)
	(C3) -- (E3) -- (D3)
	(B4) -- (C4) -- (E4)
	(B6) -- (D5) -- (E5)
	(C0) -- (J0) -- (L0)
	(I1) -- (K1) -- (L1)
	(J3) -- (L3) -- (K3)
	(C4) -- (J4)
	(D4) -- (I5) -- (K5)
	(D5) -- (I6)
	(F6) -- (U0) -- (V0)
	(J0) -- (R0)
	(U1) -- (V1) -- (W1)
	(J1) -- (R1) -- (Q1)
	(I2) -- (W2) -- (V2)
	(E2) -- (F2) -- (P2)
	(F2) -- (U3)
	(E3) -- (F3)
	(I5) -- (W5)
	(F5) -- (P5) -- (Q5) 
	(P6) -- (Q6) -- (R6)
	;

\draw[blue, thick]
	(C0) -- (E0)
	(A5) -- (A6)
	(A2) -- (A3)
	(D1) -- (E1)
	(D2) -- (E2)
	(A1) -- (B1)
	(B3) -- (C3)
	(A4) -- (B4)
	(B5) -- (D4)
	(C5) -- (E5)
	(B0) -- (D6)
	(B6) -- (C6)
	(L0) -- (K0)
	(D0) -- (I1)
	(C1) -- (J1)
	(I2) -- (K2)
	(C2) -- (J2)
	(I3) -- (K3)
	(D3) -- (I4)
	(J4) -- (L4)
	(J5) -- (L5)
	(L6) -- (K6)
	(E6) -- (F6)
	(I0) -- (W0)
	(P0) -- (Q0)
	(F0) -- (U1)
	(F1) -- (P1)
	(U2) -- (V2)
	(P2) -- (Q2)
	(W3) -- (V3)
	(J3) -- (R3)
	(F3) -- (P3)
	(U4) -- (V4)
	(E4) -- (F4)
	(R4) -- (Q4)
	(W5) -- (V5)
	(R5) -- (Q5)
	(U6) -- (F5)
	(I6) -- (W6)
	(J6) -- (R6);

\foreach \a in {0,...,6}
{
%\fill 
%	(B\a) circle (0.05);

%\filldraw[fill=white]
%	(F\a) circle (0.05)
%	(I\a) circle (0.05)
%	(J\a) circle (0.05);
}
	
\foreach \a/\b/\x in 
	{0/0/0, 0/1/1, 0/2/2, 0/3/3, 0/4/4, 0/5/5, 0/6/6, 
	1/0/2, 1/1/3, 1/2/4, 1/3/5, 1/4/6, 1/5/0, 1/6/1, 
	2/0/6, 2/1/0, 2/2/1, 2/3/2, 2/4/3, 2/5/4, 2/6/5, 
	3/0/4, 3/1/5, 3/2/6, 3/3/0, 3/4/1, 3/5/2, 3/6/3, 
	4/0/5, 4/1/6, 4/2/0, 4/3/1, 4/4/2, 4/5/3, 4/6/4, 
	5/0/3, 5/1/4, 5/2/5, 5/3/6, 5/4/0, 5/5/1, 5/6/2, 
	6/0/1, 6/1/2, 6/2/3, 6/3/4, 6/4/5, 6/5/6, 6/6/0}
\node at (X\a\b) {\x};

\foreach \a/\b/\y in 
	{0/0/2, 0/1/3, 0/2/4, 0/5/0, 0/6/1,
	1/0/4, 1/1/5, 1/2/6, 1/5/2, 1/6/3, 
	2/0/0, 2/1/1, 2/2/2, 2/5/5, 2/6/6, 
	3/0/3, 3/1/4, 3/2/5, 3/5/1, 3/6/2, 
	4/0/1, 4/1/2, 4/2/3, 4/5/6, 4/6/0, 
	5/0/5, 5/1/6, 5/2/0, 5/5/3, 5/6/4, 
	6/0/6, 6/1/0, 6/2/1, 6/5/4, 6/6/5}
\node at (Y\a\b) {\y};

\foreach \a/\b/\u in 
	{0/0/4, 0/1/5, 0/2/6, 0/6/3,
	1/0/3, 1/1/4, 1/2/5, 1/6/2,
	2/0/2, 2/1/3, 2/2/4, 2/6/1,
	3/0/5, 3/1/6, 3/2/0, 3/6/4,
	4/0/6, 4/1/0, 4/2/1, 4/6/5,
	5/0/1, 5/1/2, 5/2/3, 5/6/0,
	6/0/0, 6/1/1, 6/2/2, 6/6/6}
\node at (U\a\b) {\u};

\foreach \a/\b/\v in 
	{0/0/6, 0/1/0, 0/5/4, 0/6/5, 
	1/0/5, 1/1/6, 1/5/3, 1/6/4, 
	2/0/2, 2/1/3, 2/5/0, 2/6/1, 
	3/0/1, 3/1/2, 3/5/6, 3/6/0, 
	4/0/0, 4/1/1, 4/5/5, 4/6/6, 
	5/0/3, 5/1/4, 5/5/1, 5/6/2, 
	6/0/4, 6/1/5, 6/5/2, 6/6/3}
\node at (V\a\b) {\v};

\begin{scope}[xshift=5cm, scale=1.2]

\foreach \a in {0,...,6}
{
\fill[gray!50, rotate=51.43*\a]
	(0,0.5) -- (0,0.75) arc (90:205:0.25);

\draw[rotate=51.43*\a]
	(78:1.2) -- (90:0.9) -- (90:0.5) -- (38.57:0.5) -- (38.57:0.9) -- (50.57:1.2) -- (64.286:1.3) -- (78:1.2);

\node at (90+51.43*\a:0.35) {\a};
\node at (100+51.43*\a:0.6) {$j_{\a}$};
}

\node at (0,0) {\scriptsize $t$};

\foreach \a in {0,...,6}
\node at (115.7+51.43*\a:0.9) {\scriptsize $t_{\a}$};

\end{scope}

\end{tikzpicture}
\caption{Unique ${\mc D}_1$-patch $0264531$ for generic $7_{2+3}$-heptagon; the subsequent ${\mc D}_2$-patch; and the notation for ${\mc D}_1$-patch. Edges of the same color have the same length.}
\label{H5A}
\end{figure}
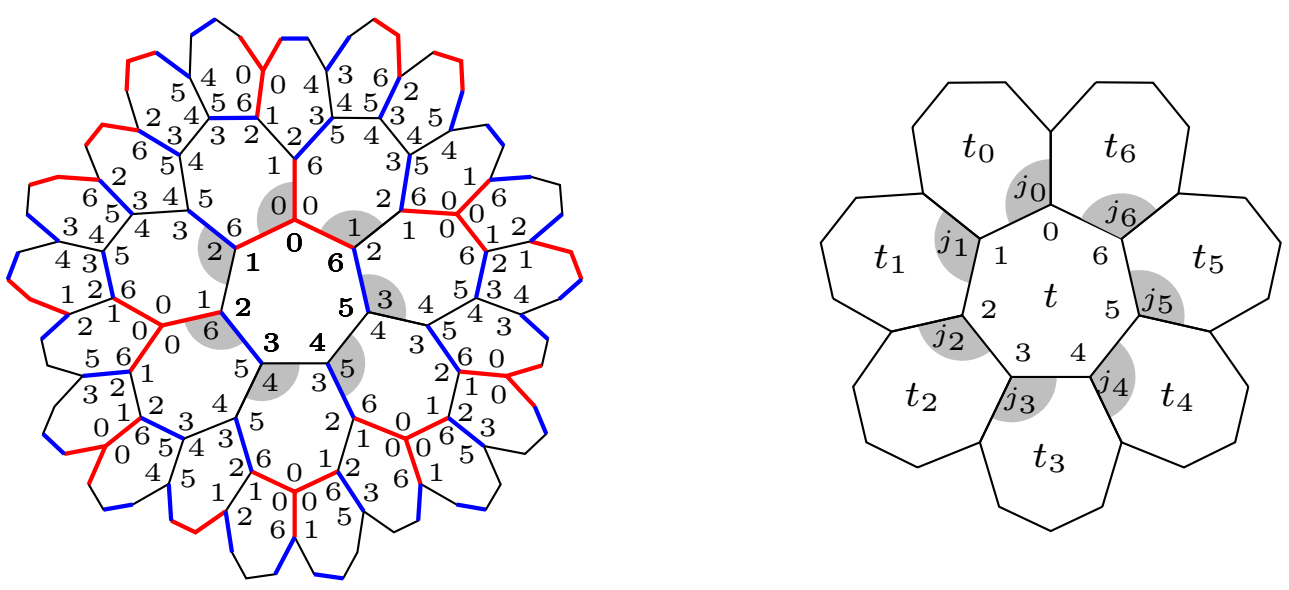

For each of the 292 edge and corner systems $E\cup V$ of rank 5 in Table \ref{table3}, we construct ${\mc D}_1$-patches satisfying $E$ and the completion $\bar{V}$. This means assigning corner labels to all eight tiles in a ${\mc D}_1$-patch, such that the following hold:
\begin{itemize}
\item The fourteen shared edges in the ${\mc D}_1$-patch give fourteen edge length equalities. These are either the trivial equality $|\bar{i}|=|\bar{i}|$, or belong to $E$. 
\item The seven vertices in the ${\mc D}_1$-patch belong to $\bar{V}$. 
\end{itemize}
Using computer program, we find that the edge and corner system \eqref{h5ev}, the first in Theorem \ref{thm1}, is the only system admitting ${\mc D}_1$-patches satisfying the above. Moreover, we find only one ${\mc D}_1$-patch for the corresponding generic $7_{2+3}$-heptagon, given by the ${\mc D}_1$ part of the left of Figure \ref{H5A}. 

We adopt the notation in \cite{zwy} to denote a ${\mc D}_1$-patch as $j_0^{\sigma_0}j_1^{\sigma_1}j_2^{\sigma_2}j_3^{\sigma_3}j_4^{\sigma_4}j_5^{\sigma_5}j_6^{\sigma_6}$ of the labels of the shaded corners. Specifically, on the right of Figure \ref{H5A}, we denote the central tile by $t$, and denote by $t_i$ the $i$-th {\em neighboring tile} of $t$. In other words, $t$ and $t_i$ share the $i$-th edge of $t$. The tile $t_i$ has two corners along the $i$-th edge of $t$. Let $j_i$ be the label of the corner of $t_i$ next to the $i$-th corner of $t$. Moreover, we introduce 
\[
\sigma_i=\begin{cases}
+, & t,t_i\text{ have same orientation} \\
-, & t,t_i\text{ have different orientation}
\end{cases}.
\] 
Then we get the sequence $j_0^{\sigma_0}j_1^{\sigma_1}j_2^{\sigma_2}j_3^{\sigma_3}j_4^{\sigma_4}j_5^{\sigma_5}j_6^{\sigma_6}$. For example,  the ${\mc D}_1$-patch on the left of Figure \ref{H5A} is $0^+2^+6^+4^+5^+3^+1^+$, and the ${\mc D}_1$-patch in Figure \ref{H6B} is $0^+3^-5^-1^-5^+2^-1^+$. In fact, we will abbreviate $j^+$ and $j^-$ as $j$ and $\tilde{j}$. Then the two ${\mc D}_1$-patches are abbreviated as $0264531$ and $0\tilde{3}\tilde{5}\tilde{1}5\tilde{2}1$.

\begin{figure}[htp]
\centering
\begin{tikzpicture}[>=latex,scale=1]

\pgftext{
	\includegraphics[scale=0.4]{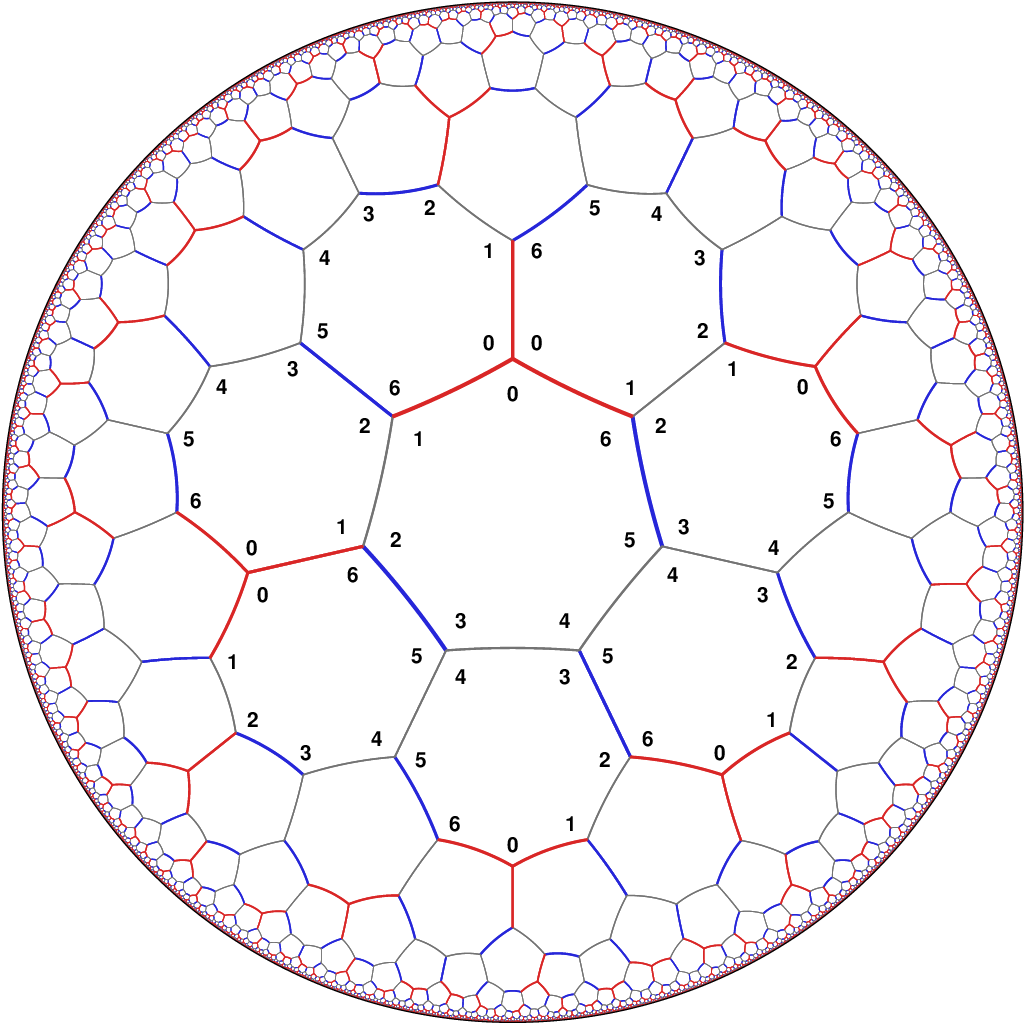}
	};
	
\foreach \a in {0,...,6}	
{
\draw[->]
	(115.715+51.43*\a:0.65) -- (115.715+51.43*\a:1.5);

\node at (115.715+51.43*\a:0.45) {\scriptsize $g_{\a}$};
}

\begin{scope}[xshift=8cm]

\pgftext{
	\includegraphics[scale=0.4]{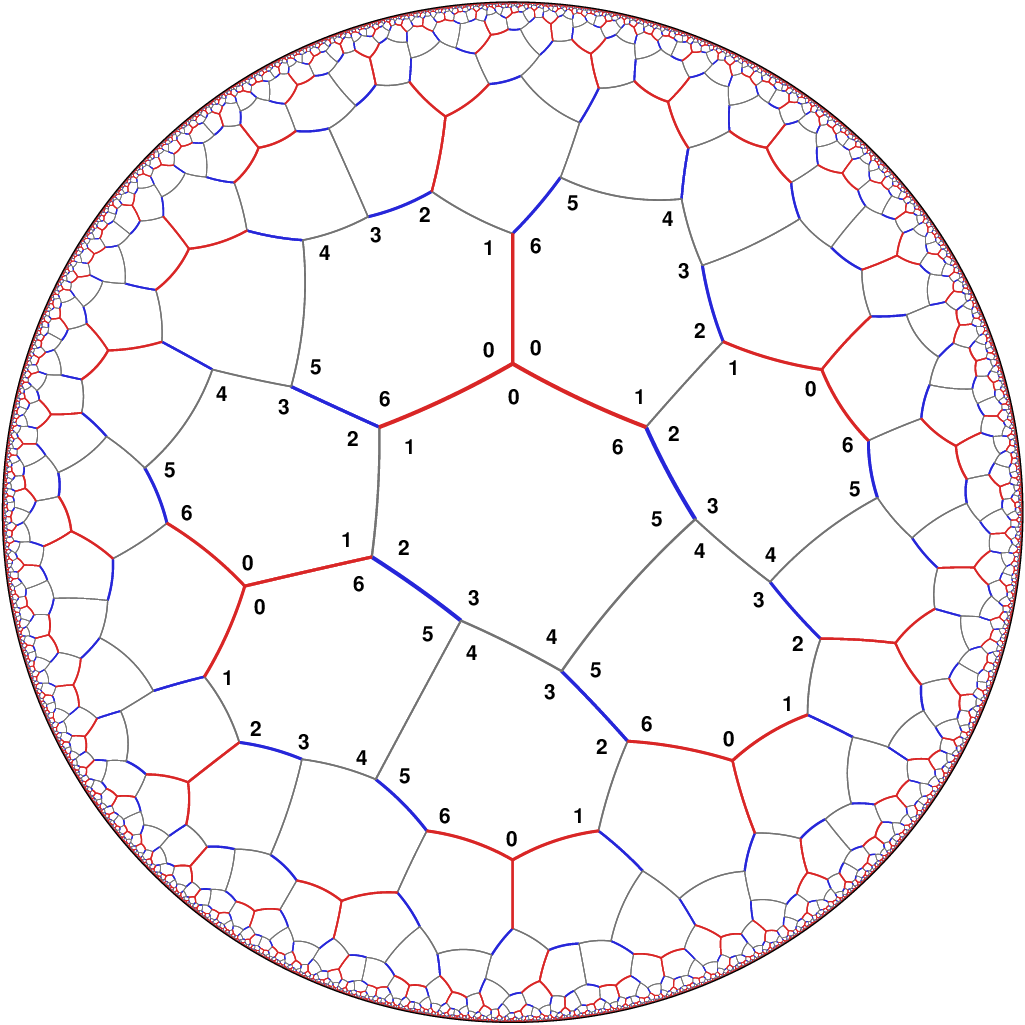}
	};
	
\end{scope}

\end{tikzpicture}
\caption{Schematic and real drawings of the unique tiling ${\mc T}(7_{2+3})$ of ${\bb H}^2$ by generic $7_{2+3}$-heptagon, and canonical generators $g_i$ of the symmetry group $G(7_{2+3})$. Red edges have the same lengths, blue edges have the same lengths. Each black edge has different lengths from the other six.}
\label{H5B}
\end{figure}

Since the ${\mc D}_1$-patch is unique, we apply $0264531$ to each of the seven tiles around the central tile and get seven ${\mc D}_1$-patches. We find these seven ${\mc D}_1$-patches to be compatible, and their union gives the unique ${\mc D}_2$-patch on the left of Figure \ref{H5A}. This implies that the application of the unique ${\mc D}_1$-patch can be repeated infinitely, and produces the unique tiling ${\mc T}(7_{2+3})$ of the Poincar\'e disk ${\bb H}^2$ in Figure \ref{H5B}.

The discussion so far is combinatorial only, and the left of Figure \ref{H5B} is only a schematic picture. To make the discussion meaningful, we need to show the existence of generic $7_{2+3}$-heptagon. By standard calculation method (see the first part of Section \ref{calculation}), we find one example satisfying the edge and corner system $7_{2+3}$, with the specific values:
\begin{align}
|\bar{0}| &=|\bar{6}|= 0.62600,\;
|\bar{1}|\approx 0.55362,\;
|\bar{2}|=|\bar{5}|\approx  0.45765,\;
|\bar{3}|\approx0.47608,\;
|\bar{4}|\approx 0.87086; \nonumber \\
[0] &=120^{\circ},\;
[1]=110^{\circ},\;
[2]=115^{\circ},\;
[3]=165^{\circ},\;
[4]=95^{\circ},\;
[5]=100^{\circ},\;
[6]=135^{\circ}. \label{5generic}
\end{align}
The right of Figure \ref{H5B} is the faithful drawing of the unique tiling ${\mc T}(7_{2+3})$ of ${\bb H}^2$ by this heptagon. It shows the heptagon is simple (the boundary does not cross itself). Moreover, the edge length values imply that the edge lengths do not satisfy any equality not in the edge part of $7_{2+3}$. It is also easy to verify that the angle values do not satisfy any degree 3 angle equality not in the corner part of $7_{2+3}$. Since all angle values are $>\frac{\pi}{2}$, there is no vertex of degree $>3$. Therefore the heptagon is a generic $7_{2+3}$-heptagon. All generic $7_{2+3}$-heptagons tile the Poincar\'e disk in the unique way in Figure \ref{H5B}. 

Let $G(7_{2+3})$ be the symmetry group of the unique tiling ${\mc T}(7_{2+3})$, and let $t$ be the central tile in Figures \ref{H5A} and \ref{H5B}. The uniqueness of the tiling implies that the following map is a one-to-one correspondence
\begin{equation}\label{oneone}
g\in G(7_{2+3})\mapsto gt\in {\mc T}(7_{2+3}).
\end{equation}
We let $g_i\in G(7_{2+3})$ be the symmetry that sends $t$ to the $i$-th neighboring tile $t_i$ (see the right of Figure \ref{H5A}). In other words, the image of \eqref{oneone} is $g_it=t_i$. 

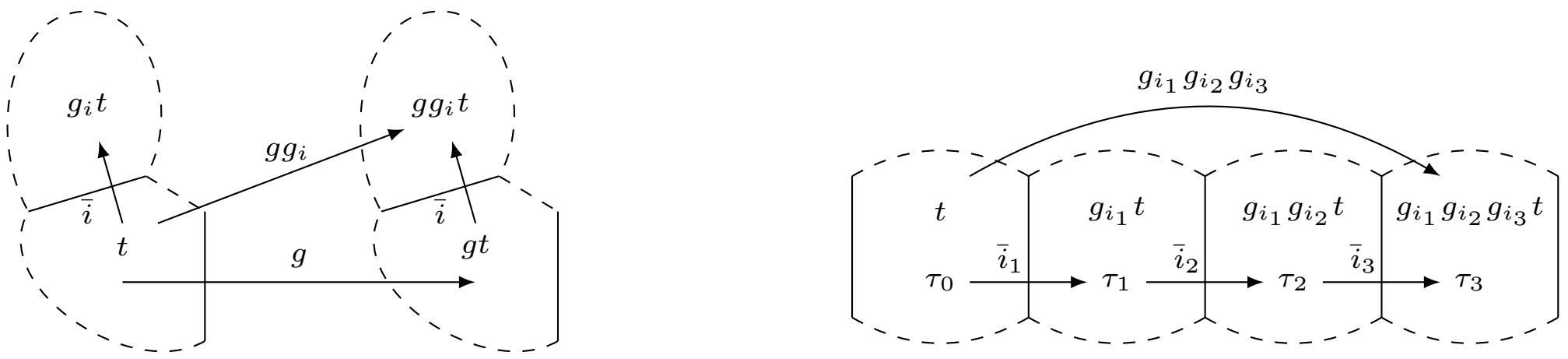
\begin{figure}[htp]
\centering
\begin{tikzpicture}[>=latex, scale=1, font=\scriptsize]

\begin{scope}[shift={(-7,0)}]

\foreach \a in {0,1}
{
\begin{scope}[xshift=3*\a cm]

\draw
	(0,-0.8) -- (0,0.3)
	(-1.5,0.3) -- (-0.5,0.6);

\draw[dashed]
	(0,-0.8) to[out=210, in=240] 
	(-1.5,0.3) to[out=120, in=180] 
	(-1,2) to[out=0, in=60] 
	(-0.5,0.6) -- (0,0.3);

\node at (-1,0.3) {$\bar{i}$};

\draw[->] (-0.7,0.2) -- (-0.9,0.9);

\end{scope}
}	

\foreach \x/\y/\z/\w in {-0.7/-0.3/2.3/-0.3, -0.4/0.2/1.7/1}
\draw[->]
	(\x,\y) -- (\z,\w);

\foreach \x/\y/\z in {0.8/-0.1/g, 0.7/0.8/gg_i,
	-0.7/0/t, 2.3/0/gt, -1/1.2/g_it, 2/1.2/gg_it}
\node at (\x,\y) {$\z$};

\end{scope}

\foreach \x in {-1,...,3}
\draw (1.5*\x,-0.6) -- ++(0,1.2);

\foreach \x in {0,...,3}
{
\draw[dashed]
	(1.5*\x,0.6) to[out=150, in=30] ++(-1.5,0)
	(1.5*\x,-0.6) to[out=210, in=-30] ++(-1.5,0);
\node at (-0.75+1.5*\x,-0.3) {$\tau_{\x}$};
}

\foreach \x in {1,2,3}	
{
\draw[->]
	(-2+1.5*\x,-0.3) -- ++(1,0);
\node at (-1.65+1.5*\x,-0.1) {$\bar{i}_{\x}$};
}

\foreach \x/\y in {0/t, 1/g_{i_1}t, 2/g_{i_1}g_{i_2}t, 3/g_{i_1}g_{i_2}g_{i_3}t}
\node at (-0.75+1.5*\x,0.3) {$\y$};

\draw[->]
	(-0.5,0.6) to[out=30, in=150] (3.5,0.6);

\node at (1.5,1.4) {$g_{i_1}g_{i_2}g_{i_3}$};

\end{tikzpicture}
\caption{Composition of symmetries.}
\label{composition}
\end{figure}

On the left of Figure \ref{composition}, a symmetry $g\in G(7_{2+3})$ takes a tile $t$ to a tile $gt$. Then $g$ also takes the $i$-th neighboring tile $g_it$ of $t$ to the $i$-th neighboring tile of $gt$. Therefore the symmetry that takes $t$ to the $i$-th neighboring tile of $gt$ is $gg_i$. 

On the right of Figure \ref{composition}, we consider a sequence of tiles $t=\tau_0,\tau_1,\tau_2,\dots,\tau_k$, such that $\tau_p$ is the $i_p$-th neighboring tile of $\tau_{p-1}$. We denote the sequence by
\begin{equation}\label{sequence}
t=\tau_0\xrightarrow{\bar{i}_1}
\tau_1\xrightarrow{\bar{i}_2}
\tau_2\xrightarrow{\bar{i}_3}
\dots\xrightarrow{\bar{i}_k}
\tau_k.
\end{equation}
By repeatedly applying the fact on the left of Figure \ref{composition}, we get 
\[
\tau_k=g_{i_1}g_{i_2}\cdots g_{i_k}t.
\]

Any tile in ${\mc T}(7_{2+3})$ is connected to $t$ by a sequence. In other words, the tile can be the right end $\tau_k$ of a sequence \eqref{sequence}. Therefore the fact above and the one-to-one correspondence imply that $g_0,g_1,\dots,g_6$ generate the symmetry group $G(7_{2+3})$.

Next, we apply the right of Figure \ref{composition} to the sequence $t\xrightarrow{\bar{0}} t_0\xrightarrow{\bar{6}} t$ on the left of Figure \ref{H5B} (and using the $t,t_i$ notation on the right of Figure \ref{H5A}), we get $g_0g_6=1$. This is an inverse relation obtained from the red edges. Similarly, we get another inverse relation $g_2g_5=1$ from the blue edges. Moreover, the edge $\bar{1}$ is only glued to $\bar{1}$. This ``self glueing'' gives the sequence $t\xrightarrow{\bar{1}} t_1\xrightarrow{\bar{1}} t$, and the corresponding inverse relation $g_1^2=g_1g_1=1$. We also get $g_3^2=g_4^2=1$ from the other two self glueing edges. 

Now consider the vertices between $t,t_i,t_{i+1}$. The vertex 162 between $t,t_0,t_1$ on the left of Figure \ref{H5B} gives the sequence $t\xrightarrow{\bar{1}} t_1\xrightarrow{\bar{2}} t_0\xrightarrow{\bar{6}} t$, and the equality $g_1g_2g_6=1$. Similarly, at the vertex 345 between $t,t_2,t_3$, we get
the sequence $t\xrightarrow{\bar{3}} t_3\xrightarrow{\bar{4}} t_2\xrightarrow{\bar{5}} t$, and the equality $g_3g_4g_5=1$. We also get $g_0^3=g_0g_0g_0=1$ from the vertex 000 between $t,t_0,t_6$. The relations obtained from the other vertices of $t$ are equivalent to the ones already obtained.

The method above for obtaining the canonical presentation of the symmetry group is similar to the one in Section 2.4 of \cite{lwwy1}. The symmetry group is 
\begin{align*}
G(7_{2+3})
&=\langle g_0,g_1,g_2,g_3,g_4,g_5,g_6 \; | \;
g_0g_6=g_2g_5=g_1^2=g_3^2=g_4^2=g_0^3=g_1g_2g_6=g_3g_4g_5=1
\rangle \\
&=\langle g_1,g_3,g_4 \; | \;
g_1^2=g_3^2=g_4^2=(g_1g_3g_4)^3=1
\rangle,
\end{align*}
with
\[
g_0=g_1g_3g_4,\;
g_2=g_3g_4,\;
g_5=g_4g_3,\;
g_6=g_4g_3g_1.
\]

A tiling of a closed surface of Euler characteristic $\chi$ by a generic $7_{2+3}$-heptagon is the quotient of ${\bb H}^2$ by a torsionless subgroup of $G(7_{2+3})$ of finite index $f=-6\chi$. The equivalence classes of tilings are the same as the conjugate classes of such subgroups. Since all tiles in ${\mc T}(7_{2+3})$ have the same orientation, the surface must be orientable. We use standard computer package GAP to find the numbers of tilings of $g{\bb T}^2$, $g=2,3,4$, in Table \ref{table5}.

\begin{figure}[htp]
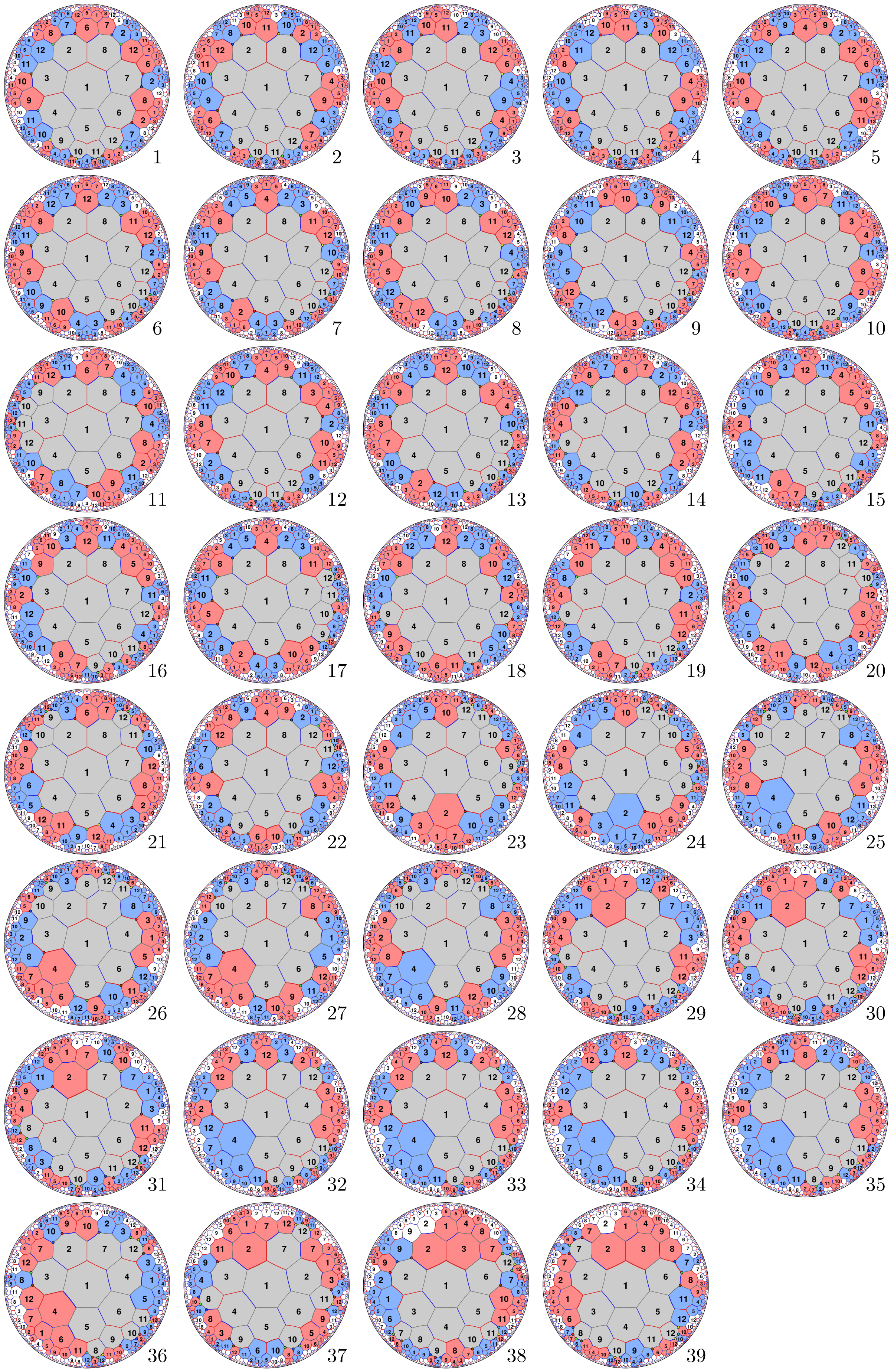

\centering
\begin{tikzpicture}[>=latex,scale=1]

%% 1-5

\foreach \a in {1,...,5}
{
\begin{scope}[shift={(-3.2*1+3.2*\a,0)}]

\pgftext{
	\includegraphics[scale=0.166]{H5domain/H5T\a}
	};

\node at (1.25,-1.25) {\a};

\end{scope}
}

%% 6-10

\foreach \a in {6,...,10}
{
\begin{scope}[shift={(-3.2*6+3.2*\a,-3.05)}]

\pgftext{
	\includegraphics[scale=0.166]{H5domain/H5T\a}
	};

\node at (1.25,-1.25) {\a};

\end{scope}
}

%% 11-15

\foreach \a in {11,...,15}
{
\begin{scope}[shift={(-3.2*11+3.2*\a,-3.05*2)}]

\pgftext{
	\includegraphics[scale=0.166]{H5domain/H5T\a}
	};

\node at (1.25,-1.25) {\a};

\end{scope}
}

%% 16-20

\foreach \a in {16,...,20}
{
\begin{scope}[shift={(-3.2*16+3.2*\a,-3.05*3)}]

\pgftext{
	\includegraphics[scale=0.166]{H5domain/H5T\a}
	};

\node at (1.25,-1.25) {\a};

\end{scope}
}

%% 21-25

\foreach \a in {21,...,25}
{
\begin{scope}[shift={(-3.2*21+3.2*\a,-3.05*4)}]

\pgftext{
	\includegraphics[scale=0.166]{H5domain/H5T\a}
	};

\node at (1.25,-1.25) {\a};

\end{scope}
}

%% 26-30

\foreach \a in {26,...,30}
{
\begin{scope}[shift={(-3.2*26+3.2*\a,-3.05*5)}]

\pgftext{
	\includegraphics[scale=0.166]{H5domain/H5T\a}
	};

\node at (1.25,-1.25) {\a};

\end{scope}
}

%% 31-35

\foreach \a in {31,...,35}
{
\begin{scope}[shift={(-3.2*31+3.2*\a,-3.05*6)}]

\pgftext{
	\includegraphics[scale=0.166]{H5domain/H5T\a}
	};

\node at (1.25,-1.25) {\a};

\end{scope}
}

%% 36-39

\foreach \a in {36,...,39}
{
\begin{scope}[shift={(-3.2*36+3.2*\a,-3.05*7)}]

\pgftext{
	\includegraphics[scale=0.166]{H5domain/H5T\a}
	};

\node at (1.25,-1.25) {\a};

\end{scope}
}

\end{tikzpicture}
\caption{Fundamental domains of ${\mc T}_k(7_{2+3},2{\bb T}^2)$, $k=1,2,\dots,39$, of tilings of $2{\bb T}^2$ by generic $7_{2+3}$-heptagon.}
\label{H5C}
\end{figure}

There are 39 tilings ${\mc T}_k(7_{2+3},2{\bb T}^2)$ of $2{\bb T}^2$ by generic $7_{2+3}$-heptagon. Figure \ref{H5C} shows the gray fundamental domain and the red/blue {\em neighboring domains}, which are the translations of the fundamental domain that are glued to the fundamental domain. The numerical labels of the tiles indicate the twelve orbits of the corresponding torsionless subgroup $G_k(7_{2+3},2{\bb T}^2)$ of index 12. We may choose any region consisting of twelve tiles labelled 1 through 12 as the fundamental domain. Figure \ref{H5C} shows the more compact choice. 

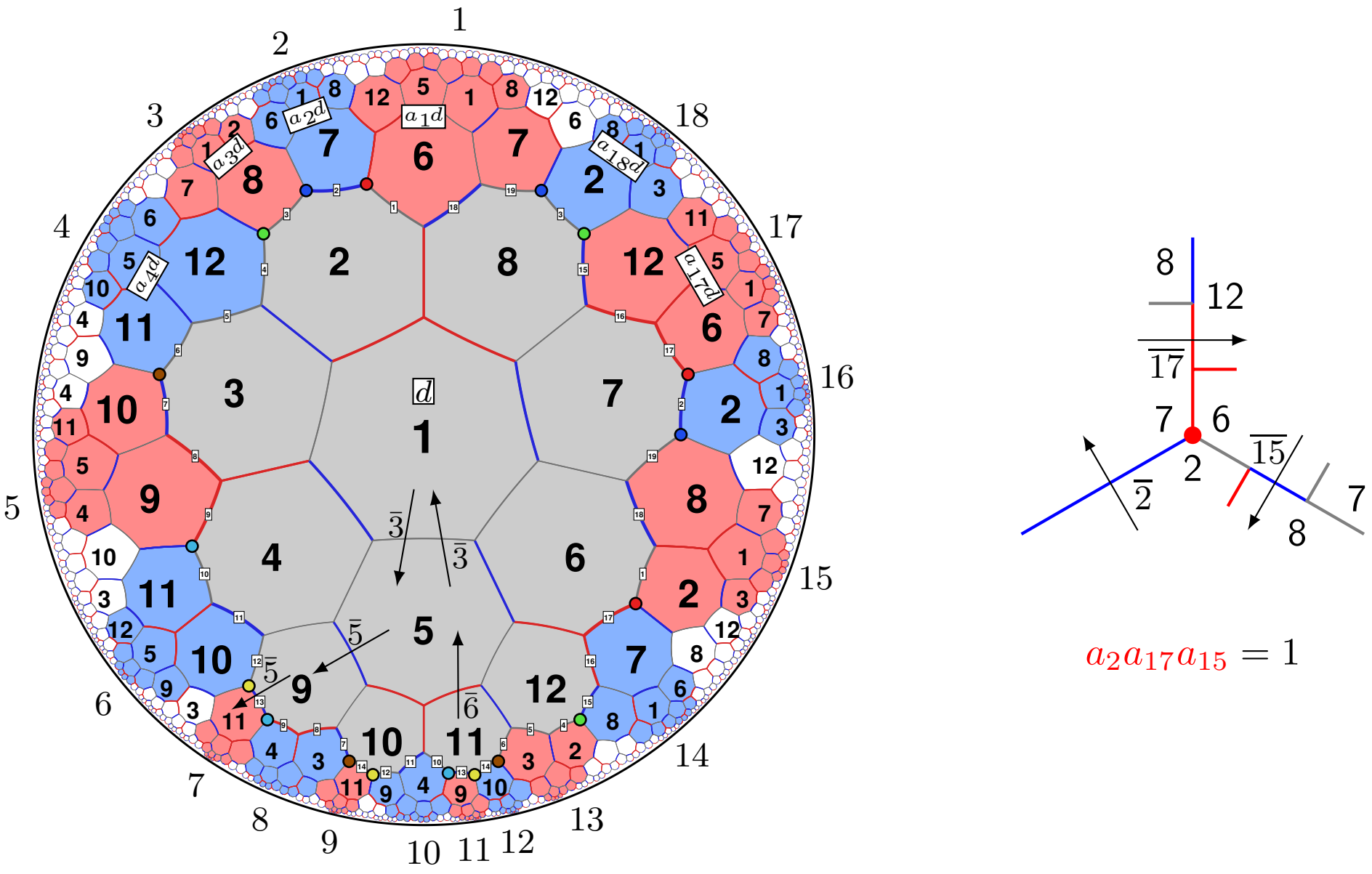
\begin{figure}[htp]
\centering
\begin{tikzpicture}[>=latex,scale=1]

\pgftext{
	\includegraphics[scale=0.4]{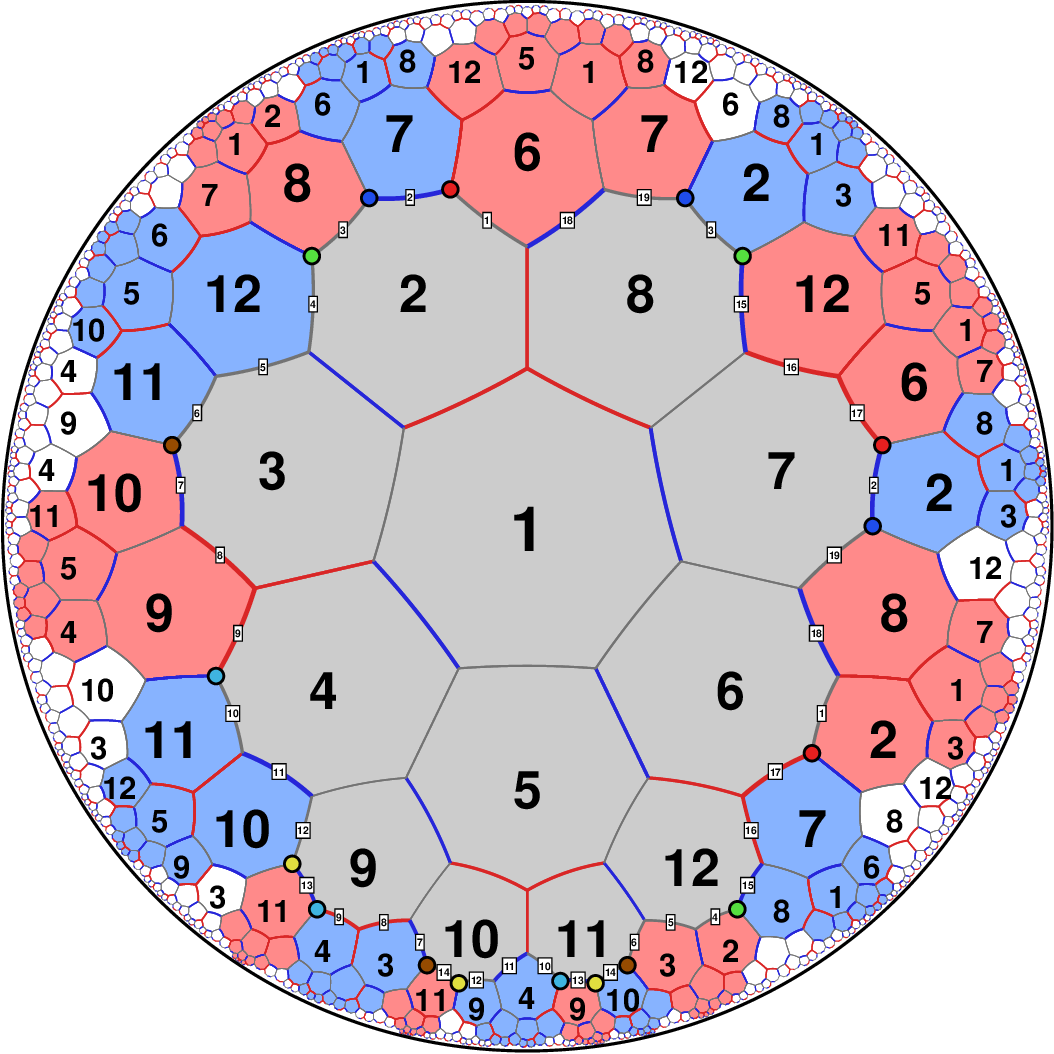}
	};

\foreach \x/\y/\a/\b in {
	90/2.9/1/0, 
	110/3.1/2/20, 
	125/3.1/3/40, 
	150/2.9/4/60, 
	30/2.9/17/-60, 
	55/3.1/18/-35}
\node[rotate=\b, inner sep=0.5, fill=white, draw, shape=rectangle] at (\x:\y) {\tiny $a_{\a}\!d$};

\node[inner sep=0.5, fill=white, draw, shape=rectangle] at (0,0.4) {\footnotesize $d$};

\foreach \a/\b in {
	85/1, 110/2, 130/3, 150/4, 190/5, 220/6, 237/7, 247/8, 257/9, 270/10, 277/11, 283/12, 293/13, 310/14, -20/15, 8/16, 30/17, 50/18}
\node at (\a:3.8) {\small \b};

\foreach \a/\x/\b/\y in {-100/0.5/-100/1.4, -100/1.8/-115/2.4, -119/2.5/-125/3.05, -83/2.6/-80/1.8, -80/1.4/-80/0.5}
\draw[->]
	(\a:\x) -- (\b:\y);

\foreach \a/\x/\b/\y in {-107/0.85/3, -109/1.9/5, -123/2.55/5, -80/2.5/6, -73/1.15/3}
\node at (\a:\x) {\footnotesize $\bar{\b}$};

%% red vertex

\begin{scope}[xshift=7cm]

\draw[red, thick]
	(0,0) -- (90:1.2)
	(90:0.6) -- ++(0:0.4)
	(-30:0.6) -- ++(-120:0.4);
	
\draw[blue, thick]
	(90:1.2) -- (90:1.8)
	(0,0) -- (210:1.8)
	(-30:0.6) -- (-30:1.2);

\draw[gray, thick]
	(90:1.2) -- ++(180:0.4)
	(0,0) -- (-30:0.6) 
	(-30:1.8) -- (-30:1.2) -- ++(60:0.4);

\foreach \a in {0,1,2}
\draw[->]
	(120+120*\a:1) -- (60+120*\a:1);

\fill[red]
	(0,0) circle (0.08);

\foreach \a/\b in {0/17, 1/2, 2/15}
\node at (110+120*\a:0.7) {\small $\overline{\b}$};

\foreach \x/\y/\a in {30/0.3/6, 150/0.3/7, -90/0.3/2, 77/1.3/12, 99/1.6/8, -43/1.3/8, -21/1.6/7}
\node at (\x:\y) {\sf \a};

\node at (0,-2) {${\color{red} a_{2}a_{17}a_{15}}=1$};
	
\end{scope}

\end{tikzpicture}
\caption{Fundamental domain of the first tiling (in Figure \ref{H5C}) ${\mc T}_1(7_{2+3},2{\bb T}^2)$ of $2{\bb T}^2$ by generic $7_{2+3}$-heptagon, and neighboring translations of the fundamental domain (with eighteen labels outside the Poincar\'e disk). On the right is the relation at the vertex {\color{red} $\bullet$}, for the generators $a_k$ of the subgroup $G_1(7_{2+3},2{\bb T}^2)$.}
\label{H5E}
\end{figure}

The left of Figure \ref{H5E} is the detail of the first tiling ${\mc T}_1(7_{2+3},2{\bb T}^2)$ of $2{\bb T}^2$ in Figure \ref{H5C}. The fundamental domain consists of twelve tiles that are obtained by applying the following elements of $G_1(7_{2+3},2{\bb T}^2)$ to the central tile that is labeled by 1.
\begin{center}
\begin{tabular}{|c||c|c|c|c|c|c|c|c|c|c|c|c|}
\hline
tile label & 1 & 2 & 3 & 4 & 5 & 6 & 7 & 8 & 9 & 10 & 11 & 12   \\
\hline
group element & 1 & $g_0$ & $g_1$ & $g_2$ & $g_3$ & $g_4$ & $g_5$ & $g_6$ & $g_3g_5$ & $g_3g_6$ & $g_3g_0$ & $g_3g_1$   \\
\hline
\end{tabular}
\end{center} 

The gray fundamental domain $d$ in the left of Figure \ref{H5E} is a 38-gon, and the surface $2{\bb T}^2$ is obtained by glueing 19 pairs of edges of the 38-gon in opposing way. The edge pairs are determined as follows: We have an edge labeled \fbox{1} between gray tile {\color{gray} 2} and red/blue tile {\color{red} 6}. This is paired with the edge (again labeled \fbox{1}) between gray tile {\color{gray} 6} and red/blue tile {\color{red} 2}. In general, the edge between gray tile {\color{gray} $i$} and red/blue tile {\color{red} $j$} is paired with the edge between gray tile {\color{gray} $j$} and red/blue tile {\color{red} $i$}. 

We start with the edge pair labeled \fbox{1}, and label the other pairs as we go counterclockwise along the boundary of the fundamental domain $d$. We find some consecutively paired edges, which we call {\em compound edge pairs}. For example, between gray {\color{gray} 2, 8} and red {\color{red} 6, 7} we find a compound edge \fbox{19,18,1} consisting of \fbox{19}, \fbox{18}, \fbox{1} in counterclockwise direction. This is paired with the compound edge $\fbox{1,18,19}=\fbox{19,18,1}^{-1}$ between gray {\color{gray} 6, 7} and red {\color{red} 2, 8}. The compound edges are actually the intersections of the fundamental domain $d$ and the neighbring domains $a_id$. 

Using (non-straight) compound edges, the fundamental domain $d$ becomes an 18-gon. The vertices of this 18-gon are indicated by colored vertices. These vertices are exactly the meeting places of gray, red and blue domains. The surface $2{\bb T}^2$ is given by a planar diagram of the 18-gon
\footnotesize
\[
{\setlength{\fboxsep}{0.5pt}% padding space
 \setlength{\fboxrule}{0.1pt}% border thickness
\boxed{19,\!18,\!1}\;
\boxed{2}\;
\boxed{3}\;
\boxed{4,\!5,\!6}\;
\boxed{7,\!8,\!9}\;
\boxed{10,\!11,\!12}\;
\boxed{13}\;
\boxed{7,\!8,\!9}^{-1}
\boxed{14}\;
\boxed{10,\!11,\!12}^{-1}
\boxed{13}^{-1}
\boxed{14}^{-1}
\boxed{4,\!5,\!6}^{-1}
\boxed{15,\!16,\!17}\;
\boxed{19,\!18,\!1}^{-1}
\boxed{2}^{-1}
\boxed{15,\!16,\!17}^{-1}
\boxed{3}^{-1}.
}
\]\normalsize

We may regard the tiling ${\mc T}_1(7_{2+3},2{\bb T}^2)$ as a tiling of $2{\bb T}^2$ by the single 18-gon $d$ of the fundamental domain. Single tile tilings of closed surfaces are the same as planar diagrams, and the single tile tiling ${\mc T}_1(7_{2+3},2{\bb T}^2)$ corresponds exactly to the planar diagram above. The 18 vertices are colored to indicate the vertices derived from the planar diagram, and are the vertices of the single tile tiling of $2{\bb T}^2$. Since the surface $2{\bb T}^2$ is orientable, we may use Proposition 5 of \cite{lwwy1} to get the standard presentation of the symmetry group $G_1(7_{2+3},2{\bb T}^2)$ of the universal cover tiling. In fact, the method is same as the one described above, based on Figure \ref{composition}, except the central tile $t$ is now the 18-gon $d$. 

The neighboring red/blue domains are the translations $a_id$ of the central fundamental domain $d$, by isometries $a_i\in G_1(7_{2+3},2{\bb T}^2)$. We indicate $a_id$ within some neighboring domains, and indicate all 18 subscripts $i$ just outside the Poincar\'e disk. The following are the formulae for the isometries $a_i$
\begin{align*}
a_1 &=g_6g_4g_2, & 
a_4 &=g_0g_4g_0g_4, &
a_7 &=g_3g_5^2g_6g_3, & 
a_{10} &=g_3g_0g_1g_5, & 
a_{13} &=g_3g_0g_4g_1, & 
a_{16} &=g_5g_5g_6, & \\
a_2 &=g_0g_2^2, & 
a_5 &=g_1g_5g_0g_3, &
a_8 &=g_2g_3g_0g_1, &
a_{11} &=g_3g_0g_2^2g_3, & 
a_{14} &=g_4g_0g_2, & 
a_{17} &=g_5g_6g_4,  \\
a_3 &=g_0g_3g_0, & 
a_6 &=g_2g_1g_6g_3, &
a_9 &=g_3g_6g_3g_6g_3, &
a_{12} &=g_3g_0g_3g_0g_3, & 
a_{15} &=g_4g_1g_6, & 
a_{18} &=g_6g_3g_6. 
\end{align*}
For example, we obtain $a_7$ from the following path that connects the tile {\color{gray} 1} in the gray domain $d$ and the tile {\color{red} 1} in the red domain $a_7d$ 
\[
{\color{gray} 1} \xrightarrow{\bar{3}} 
{\color{gray} 5} \xrightarrow{\bar{5}} 
{\color{gray} 9} \xrightarrow{\bar{5}} 
{\color{red} 11} \xrightarrow{\bar{6}} 
{\color{red} 5} \xrightarrow{\bar{3}}
{\color{red} 1}.
\]
Here {\color{gray} $i$} denotes the tile labeled $i$ in the gray fundamental domain $d$, {\color{red} $i$} denotes the tile labeled $i$ in the red neighboring domain $a_7d$. On the left of Figure \ref{H5E}, we draw the first three arrows $\xrightarrow{\bar{3}} ,\xrightarrow{\bar{5}} ,\xrightarrow{\bar{5}}$. However, the domain $a_7d$ is too small for us to draw the arrows $\xrightarrow{\bar{6}}, \xrightarrow{\bar{3}}$ between the red tiles. We notice that, as far as the edge labels $\bar{6}$ and $\bar{3}$ are concerned, the arrows are the same as the arrows between the gray tiles of the same labels. It is these two arrows between the gray tiles that we draw on the left of Figure \ref{H5E}. Then the sequence of edges $\bar{3},\bar{5},\bar{5},\bar{6},\bar{3}$ gives the formula $a_7=g_3g_5^2g_6g_3$.

These $a_i$ form canonical generators of the subgroup $G_1(7_{2+3},2{\bb T}^2)$. The relation among the generators consists of two parts. The first is the inverse relations given by the edge pairings of the 18-gon
\[
a_1a_{15}=
a_2a_{16}=
a_3a_{18}=
a_4a_{13}=
a_5a_8 =
a_6a_{10}=
a_7a_{11}=
a_9a_{12}=
a_{14}a_{17}=1.
\]
The second is the relations given by the colored vertices
\[
{\color{red} a_{2}a_{17}a_{15}}
={\color{blue} a_3a_1a_{16}}
={\color{green} a_4a_{14}a_{18}}
={\color{Dbrown} a_5a_9a_{13}}
={\color{cyan} a_6a_{11}a_8}
={\color{Gold} a_7a_{12}a_{10}}
=1.
\]
The right of Figure \ref{H5E} shows the relation ${\color{red} a_{2}a_{17}a_{15}}=1$ given by the vertex {\color{red} $\bullet$}. We indicate the labels of the tiles (in sans serif font) around the vertex, and use the labels to determine the isometries in clockwise direction. For example, the isometry going from tiles {\color{gray} 7} and {\color{gray} 8} in the gray domain to the tiles {\color{red} 6} and {\color{red} 12} outside the gray domain is $a_{17}$.

By following the rules in Proposition 5 of \cite{lwwy1}, we can practically read off the presentations of the subgroups $G_k(7_{2+3},2{\bb T}^2)$ for the tilings of $2{\bb T}^2$ in Figure \ref{H5C}.

\section{Tilings of Rank 6}
\label{rank6}

After obtaining the edge and corner system \eqref{h5ev} and the tiling by the corresponding generic $7_{2+3}$-heptagon in Section \ref{rank5}, among edge and corner systems of rank 6 in Table \ref{table3}, we remove those that can imply \eqref{h5ev} through linear algebra, up to rotation and flip of labels. There are $10+48=58$ such systems, and the number of systems of rank 6 is updated to 1525 and 11373 in the fifth column of Table \ref{table3}. Then we use computer program to construct ${\mc D}_1$-patches for the systems in the updated list. 

\begin{figure}[htp]
\centering
\begin{tikzpicture}[>=latex,scale=1]

\pgftext{
	\includegraphics[scale=0.4]{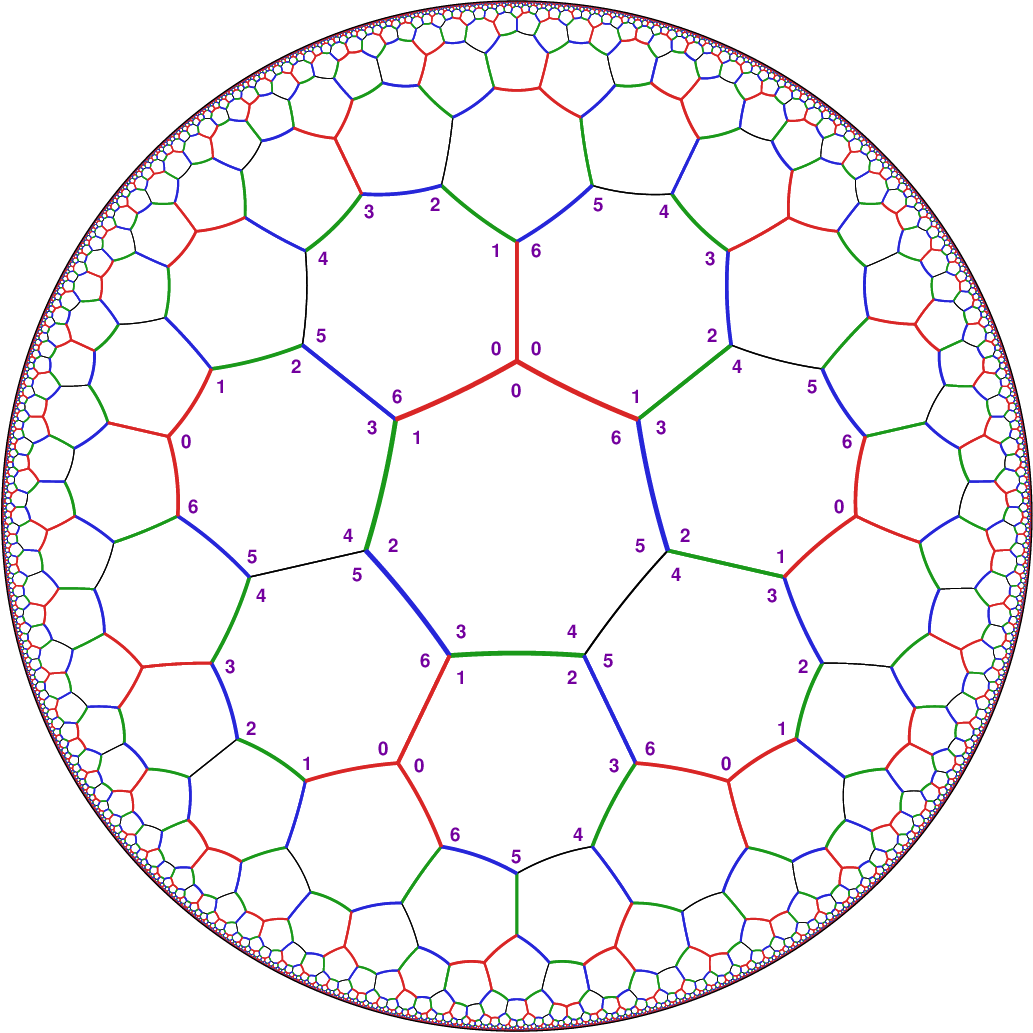}
	};
	
\foreach \a in {0,...,6}	
{
\draw[->]
	(115.715+51.43*\a:0.65) -- (115.715+51.43*\a:1.5);

\node at (115.715+51.43*\a:0.45) {\scriptsize $g_{\a}$};
}

\begin{scope}[xshift=8cm]

\pgftext{
	\includegraphics[scale=0.4]{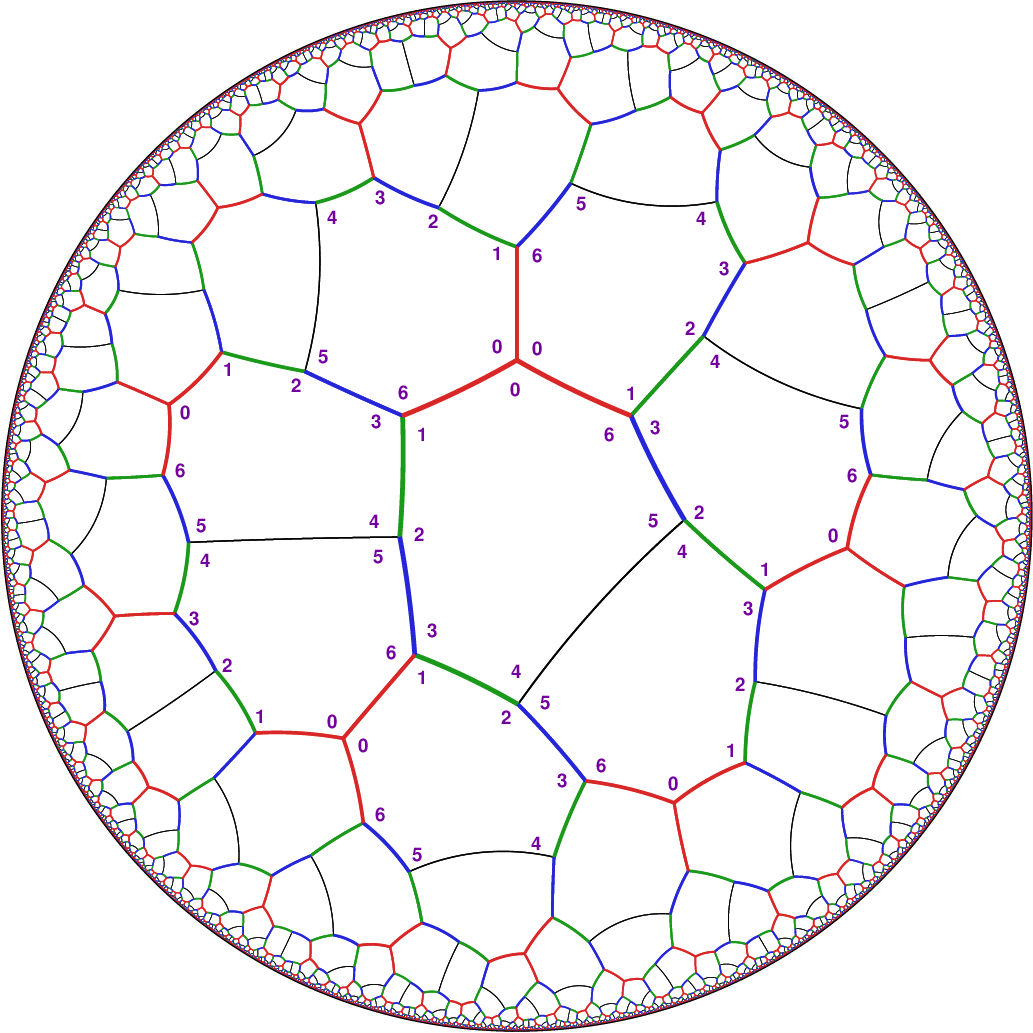}
	};
	
\end{scope}

\end{tikzpicture}
\caption{Schematic and real drawings of the unique tiling ${\mc T}(7_{3+3})$ of ${\bb H}^2$ by generic $7_{3+3}$-heptagon, and canonical generators $g_i$ of the symmetry group $G(7_{3+3})$. Edges of the same color have the same lengths. }
\label{H6B}
\end{figure}

We find that the edge and corner system \eqref{h6ev}, the second in Theorem \ref{thm1}, is the only system in the updated list that admits ${\mc D}_1$-patches. Moreover, we find only one ${\mc D}_1$-patch $0\tilde{3}\tilde{5}\tilde{1}5\tilde{2}1$ for the corresponding generic $7_{3+3}$-heptagon, given by the ${\mc D}_1$-patch at the center of the left of Figure \ref{H6B}. Similar to the tiling for generic $7_{2+3}$-heptagons, this uniquely determines a tiling ${\mc T}(7_{3+3})$ of the Poincar\'e disk ${\bb H}^2$ in Figure \ref{H6B}. We also find a generic $7_{3+3}$-heptagon with the specific values:
\begin{align}
|\bar{0}| &=|\bar{6}|=0.53500,\;
|\bar{1}|=|\bar{3}|\approx 0.50166,\;
|\bar{2}|=|\bar{5}|\approx 0.49762,\;
|\bar{4}|\approx 1.04756; \nonumber \\
[0] &=120^{\circ},\;
[1]=110^{\circ},\;
[2]=165^{\circ},\;
[3]=115^{\circ},\;
[4]=95^{\circ},\;
[5]=100^{\circ},\;
[6]=135^{\circ}. \label{6generic}
\end{align}
The right of Figure \ref{H6B} is the faithful picture of the unique tiling ${\mc T}(7_{3+3})$ of ${\bb H}^2$ by this specific heptagon. 

The symmetry group  $G(7_{3+3})$ of ${\mc T}(7_{3+3})$ is transitive on all the tiles, and is generated by seven symmetries $g_i$ satisfying $g_it=t_i$. Unlike $G(7_{2+3})$, where all seven generators preserve the orientation, for $G(7_{3+3})$, $g_0,g_4,g_6$ preserve the orientation and $g_1,g_2,g_3,g_5$ reverse the orientation. By considering the edges glued together, we get
\[
g_0g_6=g_1g_3=g_2g_5=g_4^2=1.
\]
By going around the vertices 000, 136, 245, we get 
\[
g_0^3=g_1g_2g_6=g_2g_4g_3=1.
\] 
The symmetry group is 
\begin{align*}
G(7_{3+3})
&=\langle g_0,g_1,g_2,g_3,g_4,g_5,g_6 \; | \;
g_0g_6=g_1g_3=g_2g_5=g_4^2=g_0^3=g_1g_2g_6=g_2g_4g_3=1
\rangle \\
&=\langle g_0,g_1 \; | \;
g_0^3=(g_0^{-1}g_1^2)^2=1
\rangle,
\end{align*}
with
\[
g_2=g_1^{-1}g_0,\;
g_3=g_1^{-1},\;
g_4=g_0^{-1}g_1^2,\;
g_5=g_0^{-1}g_1,\;
g_6=g_0^{-1}.
\]

A tiling of a closed surface of Euler characteristic $\chi$ by a generic $7_{3+3}$-heptagon is the quotient of ${\bb H}^2$ by a torsionless subgroup of $G(7_{3+3})$ of finite index $f=-6\chi$. Since tiles in ${\mc T}(7_{3+3})$ may not have the same orientation, the surface may or may not be orientable. The surface is orientable if and only if the numbers of $g_1,g_2,g_3,g_5$ (or only $g_1$ in the simplified presentation) in the generators of the corresponding torsionless subgroup are all even.

\begin{proposition}\label{even6}
If a generic $7_{3+3}$-heptagon tiles $g{\bb P}^2$, then $g$ is even.
\end{proposition}

\begin{proof}
Let $G$ be the torsionless subgroup of $G(7_{3+3})$ of index $f=-6\chi$, $\chi= 2-g$, that gives a tiling ${\mc T}$ of $g{\bb P}^2$ by a
generic $7_{3+3}$-heptagon. A group element $x\in G$ gives an element $\bar{x}\in S_f$ of the permutation group of all
tiles in ${\mc T}$. The map $x\mapsto \bar{x}$ is a homomorphism. 

Since $g_4^2 = 1$, we know $\bar{g}_4$ is a product of disjoint transpositions. Since the conjugates of $g_4$ are torsion 2 elements, and $G$ is torsionless, we know the conjugates are not in $G$. Therefore $\bar{g}_4$ acts freely on the $f=-6\chi$ tiles in ${\mc T}$. Then $\bar{g}_4$ is a product of $\frac{f}{2}=-3\chi$ disjoint transpositions, and the sign of $\bar{g}_4$ is $(-1)^{-3\chi}$. Since $g_4 = g_0^{-1}g_1^2$, we know the sign of $\bar{g}_0$ is also $(-1)^{-3\chi}$. On the other hand, since $g_0^3=1$, we know the sign of $\bar{g}_0$ is 1. Therefore $(-1)^{-3\chi}=1$, and $\chi= 2-g$ is even.
\end{proof}
 
\begin{figure}[htp]
\centering
\begin{tikzpicture}[>=latex,scale=1]

\pgftext{
	\includegraphics[scale=0.4]{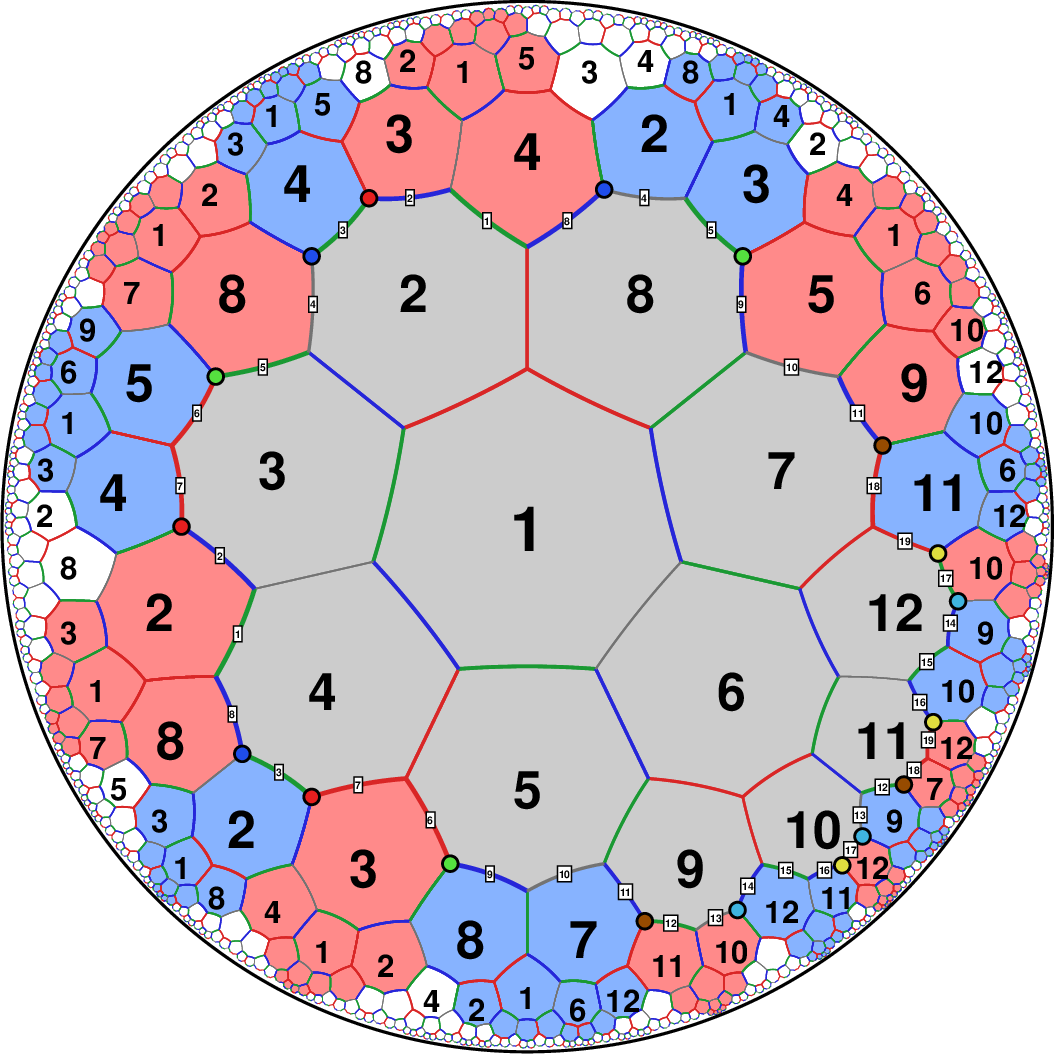}
	};

\begin{scope}[xshift=8cm]

\pgftext{
	\includegraphics[scale=0.4]{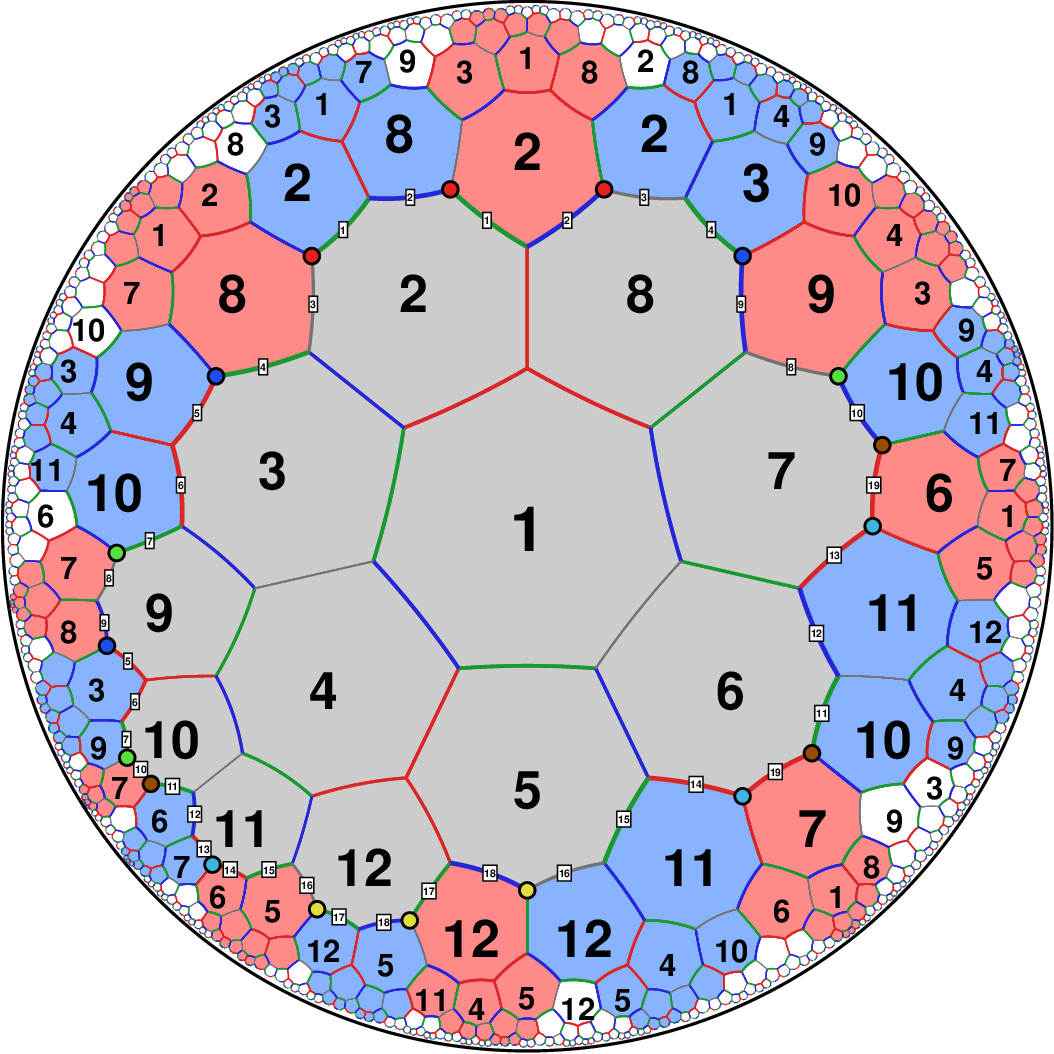}
	};

\foreach \a/\b in {
	90/\tilde{1}, 115/\tilde{2}, 140/3, 168/\tilde{4}, 187/\tilde{5}, 203/\tilde{6}, 212/7, 220/\tilde{8}, 235/9, 250/\tilde{10}, 263/\tilde{11}, 290/12, 310/\tilde{13}, -20/\tilde{14}, 0/\tilde{15}, 17/16, 40/\tilde{17}, 65/18}
\node at (\a:3.8) {\small $\b$};

\end{scope}

\end{tikzpicture}
\caption{First tiling ${\mc T}_1(7_{3+3},2{\bb T}^2)$ of $2{\bb T}^2$, and first tiling ${\mc T}_1(7_{3+3},4{\bb P}^2)$ of $4{\bb P}^2$, in Figure \ref{H6C}, by generic $7_{3+3}$-heptagon. On the right, the label $\tilde{k}$ outside the Poincar\'e disk means $a_k$ reverses orientation, i.e., the fundamental domain and the $k$-th neighboring domain have opposite orientation.}
\label{H6E}
\end{figure}

\begin{figure}[htp]
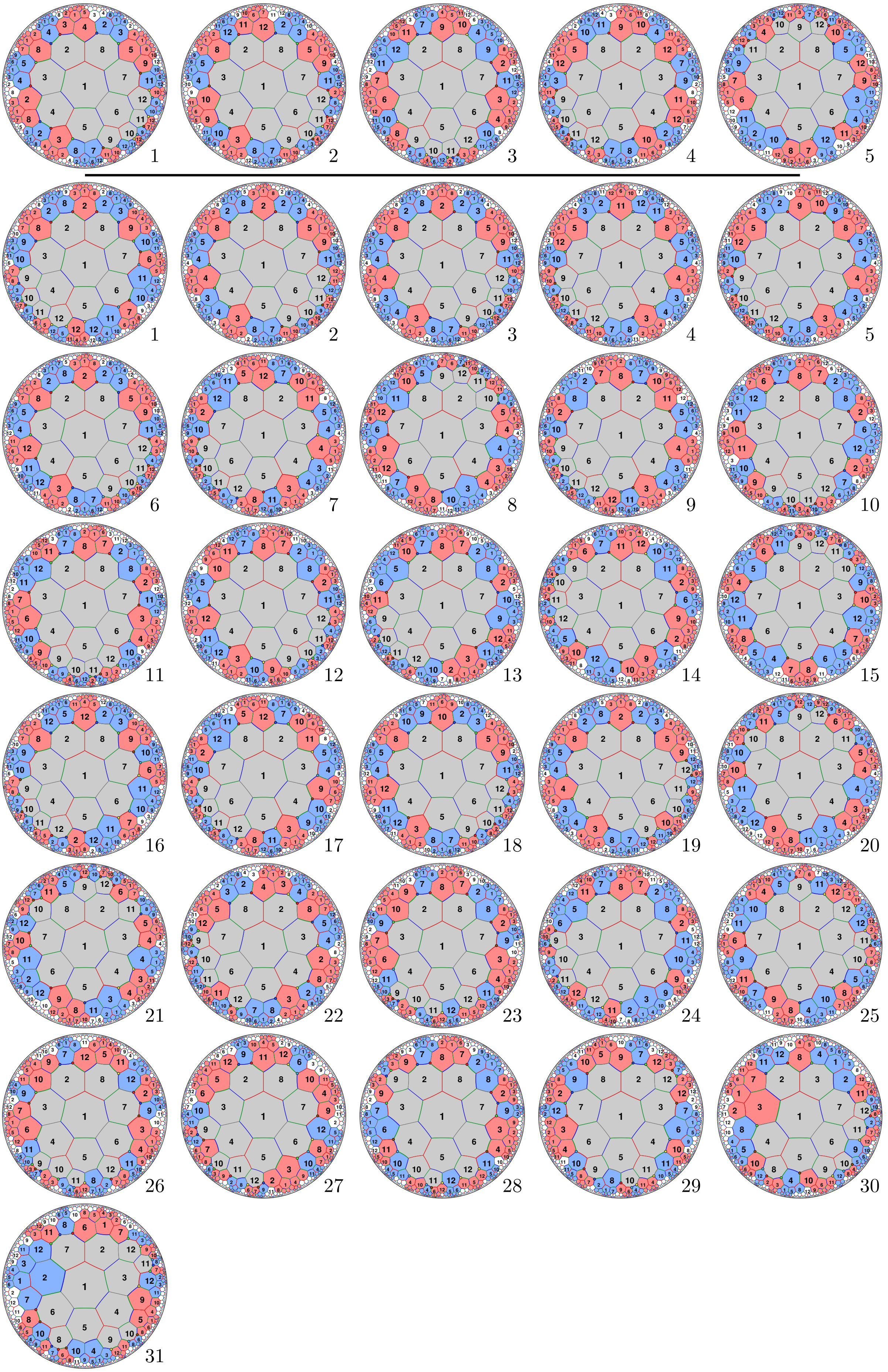

\centering
\begin{tikzpicture}[>=latex,scale=1]

%% 1-5 2T^2

\foreach \a in {1,...,5}
{
\begin{scope}[shift={(-3.2*1+3.2*\a,3.2)}]

\pgftext{
	\includegraphics[scale=0.166]{H6domainT/H6T\a}
	};

\node at (1.25,-1.25) {\a};

\end{scope}
}

%% 1-5 4T^2

\foreach \a in {1,...,5}
{
\begin{scope}[shift={(-3.2*1+3.2*\a,0)}]

\pgftext{
	\includegraphics[scale=0.166]{H6domainP/H6P\a}
	};

\node at (1.25,-1.25) {\a};

\end{scope}
}

%% 6-10

\foreach \a in {6,...,10}
{
\begin{scope}[shift={(-3.2*6+3.2*\a,-3.05*1)}]

\pgftext{
	\includegraphics[scale=0.166]{H6domainP/H6P\a}
	};

\node at (1.25,-1.25) {\a};

\end{scope}
}

%% 11-15

\foreach \a in {11,...,15}
{
\begin{scope}[shift={(-3.2*11+3.2*\a,-3.05*2)}]

\pgftext{
	\includegraphics[scale=0.166]{H6domainP/H6P\a}
	};

\node at (1.25,-1.25) {\a};

\end{scope}
}

%% 16-20

\foreach \a in {16,...,20}
{
\begin{scope}[shift={(-3.2*16+3.2*\a,-3.05*3)}]

\pgftext{
	\includegraphics[scale=0.166]{H6domainP/H6P\a}
	};

\node at (1.25,-1.25) {\a};

\end{scope}
}

%% 21-25

\foreach \a in {21,...,25}
{
\begin{scope}[shift={(-3.2*21+3.2*\a,-3.05*4)}]

\pgftext{
	\includegraphics[scale=0.166]{H6domainP/H6P\a}
	};

\node at (1.25,-1.25) {\a};

\end{scope}
}

%% 26-30

\foreach \a in {26,...,30}
{
\begin{scope}[shift={(-3.2*26+3.2*\a,-3.05*5)}]

\pgftext{
	\includegraphics[scale=0.166]{H6domainP/H6P\a}
	};

\node at (1.25,-1.25) {\a};

\end{scope}
}

%% 31

\foreach \a in {31}
{
\begin{scope}[shift={(-3.2*31+3.2*\a,-3.05*6)}]

\pgftext{
	\includegraphics[scale=0.166]{H6domainP/H6P\a}
	};

\node at (1.25,-1.25) {\a};

\end{scope}
}

\draw[very thick] 
	(0,1.6) -- ++(12.8,0);

\end{tikzpicture}
\caption{Fundamental domains of ${\mc T}_k(7_{3+3},2{\bb T}^2)$, $k=1,2,\dots,5$, and ${\mc T}_k(7_{3+3},4{\bb P}^2)$, $k=1,2,\dots,31$, by generic $7_{3+3}$-heptagon.}
\label{H6C}
\end{figure}

There are 5 tilings of $2{\bb T}^2$, and 31 tilings of $4{\bb P}^2$, by generic $7_{3+3}$-heptagon. Figure \ref{H6C} presents the gray fundamental domains and the red/blue neighboring domains. Figure \ref{H6E} is the details of the first tilings ${\mc T}_1(7_{3+3},2{\bb T}^2)$ and ${\mc T}_1(7_{3+3},4{\bb P}^2)$ in Figure \ref{H6C}. The fundamental domains consist of twelve tiles that are obtained by applying the following group elements to the central tile.
\begin{center}
\begin{tabular}{|c||c|c|c|c|c|c|c|c|c|c|c|c|}
\hline
tile label & 1 & 2 & 3 & 4 & 5 & 6 & 7 & 8 & 9 & 10 & 11 & 12   \\
\hline
group element for ${\mc T}_1(7_{3+3},2{\bb T}^2)$ & 1 & $g_0$ & $g_1$ & $g_2$ & $g_3$ & $g_4$ & $g_5$ & $g_6$ & $g_6g_4$ & $g_0g_4$ & $g_1g_4$ & $g_2g_4$   \\
\hline
group element for ${\mc T}_1(7_{3+3},4{\bb P}^2)$ & 1 & $g_0$ & $g_1$ & $g_2$ & $g_3$ & $g_4$ & $g_5$ & $g_6$ & $g_0g_2$ & $g_1g_2$ & $g_2g_2$ & $g_3g_2$   \\
\hline
\end{tabular}
\end{center}
The surface $2{\bb T}^2$ is given by the planar diagram\footnotesize
\[
{\setlength{\fboxsep}{0.5pt}% padding space
 \setlength{\fboxrule}{0.1pt}% border thickness
\boxed{8,\!1,\!2}\;
\boxed{3}\;
\boxed{4,\!5}\;
\boxed{6,\!7}\;
\boxed{8,\!1,\!2}^{-1}
\boxed{3}^{-1}
\boxed{6,\!7}^{-1}
\boxed{9,\!10,\!11}\;
\boxed{12,\!13}\;
\boxed{14,\!15,\!16}\;
\boxed{17}\;
\boxed{12,\!13}^{-1}
\boxed{18,\!19}\;
\boxed{14,\!15,\!16}^{-1}
\boxed{17}^{-1}
\boxed{18,\!19}^{-1}
\boxed{9,\!10,\!11}^{-1}
\boxed{4,\!5}^{-1}.
}
\]\normalsize
The surface $4{\bb P}^2$ is given by the planar diagram\footnotesize
\[
{\setlength{\fboxsep}{0.5pt}% padding space
 \setlength{\fboxrule}{0.1pt}% border thickness
\boxed{2,\!1}\;
\boxed{2,\!1}\;
\boxed{3,\!4}\;
\boxed{5,\!6,\!7}\;
\boxed{8,\!9}\;
\boxed{5,\!6,\!7}\;
\boxed{10}\;
\boxed{11,\!12,\!13}\;
\boxed{14,\!15,\!16}\;
\boxed{17,\!18}\;
\boxed{17,\!18}\;
\boxed{14,\!15,\!16}^{-1}
\boxed{19}\;
\boxed{11,\!12,\!13}\;
\boxed{19}\;
\boxed{10}^{-1}
\boxed{8,\!9}\;
\boxed{3,\!4}^{-1}.
}
\]\normalsize

By the same method as in Section \ref{rank5}, we can easily get canonical presentations of the corresponding torsionless subgroups $G_k(7_{3+3},2{\bb T}^2)$ and $G_k(7_{3+3},4{\bb P}^2)$ from the planar diagrams. We may apply Propositions 5 and 6 of \cite{lwwy1} respectively to orientable and non-orientable surfaces. The method is the same as the one based on Figure \ref{composition}, except the central tile $t$ is now the 18-gon of fundamental domain. 

Now we give the canonical presentation of $G_1(7_{3+3},4{\bb P}^2)$, based on the right of Figure \ref{H6E}. The generators are
\begin{align*}
a_1 &=g_0g_1g_6, & 
a_4 &=g_1g_0g_1g_5, &
a_7 &=g_2^4, & 
a_{10} &=g_3g_6g_2g_1, & 
a_{13} &=g_4g_0g_2, & 
a_{16} &=g_5^4, & \\
a_2 &=g_0g_2g_0, & 
a_5 &=g_1g_5g_4g_2, &
a_8 &=g_2g_1g_5g_4, &
a_{11} &=g_3g_5g_0g_1, & 
a_{14} &=g_4g_1g_5^2, & 
a_{17} &=g_5g_4g_1g_5,  \\
a_3 &=g_0g_4g_0, & 
a_6 &=g_1g_5g_6g_3, &
a_9 &=g_2g_0g_4g_1, &
a_{12} &=g_3g_4g_0g_1, & 
a_{15} &=g_5g_6g_4, & 
a_{18} &=g_6g_3^2. 
\end{align*}
The relations are
\[
a_1a_2=a_3a_{18}=a_4a_6=a_5a_{17}=a_7a_{16}=a_8a_{14}=a_9a_{12}=a_{10}a_{11}=a_{13}a_{15}=1,
\]
and
\[
{\color{red} a_2^2a_{18}}
={\color{blue} a_4a_5a_{18}}
={\color{green} a_5a_{16}a_6}
={\color{Dbrown} a_8a_{13}a_{16}}
={\color{cyan} a_8a_{15}a_{12}}
={\color{Gold} a_9a_{11}^2}
=1.
\]

\section{Tilings of Ranks 7 through 9}
\label{rank9}

After obtaining the edge and corner systems  \eqref{h5ev} and  \eqref{h6ev}, and the tilings by the corresponding generic $7_{2+3}$- and $7_{3+3}$-heptagons, among edge and corner systems of ranks 7, 8 and 9 in the fourth column of Table \ref{table3}, we remove those that can imply \eqref{h5ev} and  \eqref{h6ev} through linear algebra, up to rotation and flip of labels. Then the numbers of systems of ranks 7, 8 and 9 are updated in the fifth column of Table \ref{table3}. 

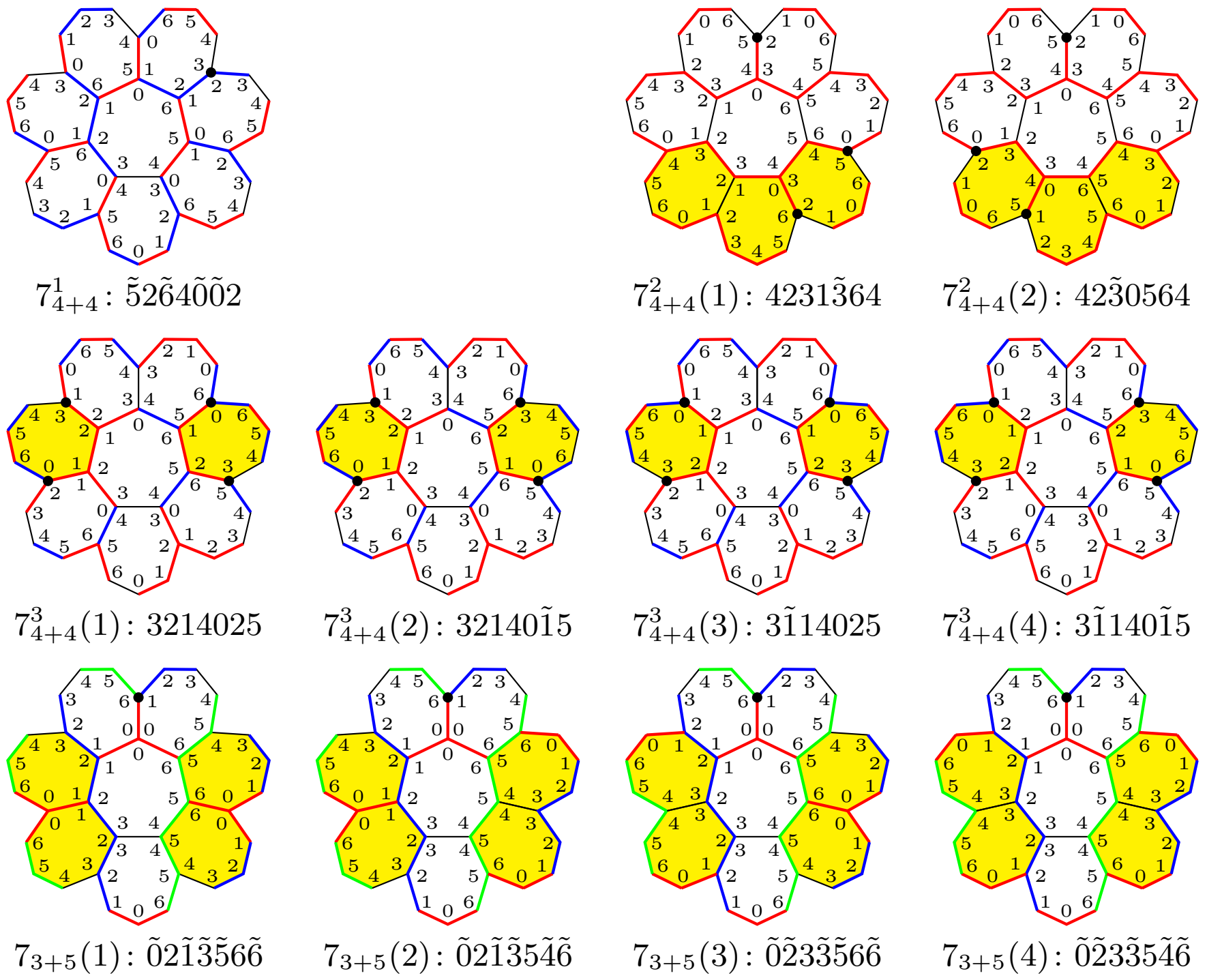
\begin{figure}[htp]
\centering
\begin{tikzpicture}[>=latex,scale=1, font=\tiny]

\foreach \a in {0,...,6}
\foreach \b in {0,...,6}
\foreach \x in {0,1,2,3}
\foreach \y in {0,1,2}
{
\begin{scope}[shift={(3*\x, -3.2*\y)}]

\coordinate (A\x\y\a) at (90+51.43*\a:0.5);
\coordinate (B\x\y\a 0) at (90+51.43*\a:0.5);
\coordinate (B\x\y\a 1) at (90+51.43*\a:0.9);
\coordinate (B\x\y\a 2) at (102+51.43*\a:1.2);
\coordinate (B\x\y\a 3) at (115.714+51.43*\a:1.3);
\coordinate (B\x\y\a 4) at (129.43+51.43*\a:1.2);
\coordinate (B\x\y\a 5) at (141.43+51.43*\a:0.9);
\coordinate (B\x\y\a 6) at (141.43+51.43*\a:0.5);

\coordinate (X\x\y\a 0) at (100+51.43*\a:0.6);
\coordinate (X\x\y\a 1) at (99+51.43*\a:0.87);
\coordinate (X\x\y\a 2) at (105+51.43*\a:1.1);
\coordinate (X\x\y\a 3) at (115.714+51.43*\a:1.17);
\coordinate (X\x\y\a 4) at (126.43+51.43*\a:1.1);
\coordinate (X\x\y\a 5) at (132.43+51.43*\a:0.87);
\coordinate (X\x\y\a 6) at (131.43+51.43*\a:0.6);

\end{scope}
}

\foreach \x in {2,3}
\fill[yellow]
	(A\x 02) -- (A\x 03) -- (A\x 04) -- (A\x 05) -- (B\x 046) -- (B\x 045) -- (B\x 044) -- (B\x 043) -- (B\x 042) -- (B\x 041) -- (B\x 034) -- (B\x 033) -- (B\x 032) -- (B\x 031) -- (B\x 024) -- (B\x 023) -- (B\x 022) -- (B\x 021);
	
\foreach \x in {0,1,2,3}
\fill[yellow]
	(B\x 110) -- (B\x 111) -- (B\x 112) -- (B\x 113) -- (B\x 114) -- (B\x 115) -- (B\x 116)
	(B\x 150) -- (B\x 151) -- (B\x 152) -- (B\x 153) -- (B\x 154) -- (B\x 155) -- (B\x 156)
	(A\x 21) -- (A\x 22) -- (A\x 23) -- (B\x 225) -- (B\x 224) -- (B\x 223) -- (B\x 222) -- (B\x 221) -- (B\x 214) -- (B\x 213) -- (B\x 212) -- (B\x 211)
	(A\x 24) -- (A\x 25) -- (A\x 26) -- (B\x 255) -- (B\x 254) -- (B\x 253) -- (B\x 252) -- (B\x 251) -- (B\x 244) -- (B\x 243) -- (B\x 242) -- (B\x 241);

\foreach \a in {0,...,6}
\foreach \x in {0,1,2,3}
\foreach \y in {0,1,2}
{
\draw[shift={(3*\x, -3.2*\y)}, rotate=51.43*\a]
	(78:1.2) -- (90:0.9) -- (90:0.5) -- (38.57:0.5) -- (38.57:0.9) -- (50.57:1.2) -- (64.286:1.3) -- (78:1.2);
	
\node[shift={(3*\x, -3.2*\y)}] at (90+51.43*\a:0.35) {\a};
}

\fill[white]
	(1.5,1.5) rectangle (4.5,-1.5);

\draw[red, thick]
	(A000) -- (A001)
	(A004) -- (A005) -- (A006);
	
\draw[blue, thick]
	(A000) -- (A006)
	(A001) -- (A002) -- (A003);

\foreach \x in {2,3}
\draw[red, thick]
	(B\x 013) -- (B\x 012) -- (B\x 011) -- (A\x 01) -- (A\x 00) -- (A\x 06)  --  (B\x 055) -- (B\x 054) -- (B\x 053)
	(A\x 00) -- (B\x 001)
	(B\x 002) -- (B\x 003) -- (B\x 004)
	(B\x 062) -- (B\x 063) -- (B\x 064)
	(B\x 014) -- (B\x 015) -- (A\x 02) -- (A\x 03) -- (A\x 04) -- (A\x 05) -- (B\x 051) -- (B\x 052)
	(B\x 023) -- (B\x 024)
	(B\x 032) -- (B\x 033) -- (B\x 034)
	(B\x 042) -- (B\x 043);

\foreach \x in {0,1,2,3}
{
\draw[red, thick]
	(A\x 10) -- (A\x 11) -- (A\x 12) -- (A\x 13)
	(B\x 151) -- (A\x 15) -- (A\x 16) -- (B\x 155)
	(B\x 102) -- (B\x 103)
	(B\x 104) -- (B\x 105) -- (B\x 106)
	(B\x 120) -- (B\x 121) -- (B\x 122)
	(B\x 124) -- (B\x 125) -- (B\x 132)
	(B\x 133) -- (B\x 134) -- (B\x 135) -- (B\x 136)
	(B\x 141) -- (B\x 142) -- (B\x 143)
	(B\x 162) -- (B\x 163) -- (B\x 164) -- (B\x 165)
	(A\x 21) -- (A\x 20) -- (A\x 26)
	(A\x 20) -- (B\x 201)
	(B\x 232) -- (B\x 233) -- (B\x 234);

\draw[blue, thick]
	(A\x 10) -- (A\x 16)
	(A\x 14) -- (A\x 15)
	(B\x 101) -- (B\x 102)
	(B\x 103) -- (B\x 104)
	(B\x 123) -- (B\x 124)
	(B\x 125) -- (B\x 126)
	(B\x 144) -- (B\x 145)
	(B\x 161) -- (B\x 162)
	(B\x 204) -- (B\x 205) -- (A\x 21) -- (A\x 22) -- (A\x 23) -- (B\x 231) -- (B\x 232)
	(B\x 263) -- (B\x 264) -- (B\x 265);

\draw[green, thick]
	(B\x 234) -- (B\x 235) -- (A\x 24) -- (A\x 25) -- (A\x 26) -- (B\x 261) -- (B\x 262)
	(B\x 201) -- (B\x 202) -- (B\x 203);
}
		
\foreach \x/\y/\a/\u/\v in 
	{0/0/0/0/1, 0/0/0/4/5, 0/0/1/2/3, 0/0/1/3/4, 0/0/2/0/1, 0/0/2/1/2, 0/0/3/0/1, 0/0/3/1/2, 0/0/3/3/4, 0/0/4/1/2, 0/0/4/2/3, 0/0/5/1/2, 0/0/5/2/3, 0/0/6/2/3, 0/0/6/3/4,
	2/0/2/1/2, 2/0/2/4/5, 2/0/3/1/2, 2/0/4/0/1, 2/0/4/3/4,
	3/0/2/2/3, 3/0/2/5/6, 3/0/3/4/5, 3/0/4/1/2, 3/0/4/4/5,
	0/1/1/3/4, 0/1/5/3/4,
	1/1/1/3/4, 1/1/1/2/3,
	2/1/1/2/3, 2/1/5/3/4,
	3/1/1/2/3, 3/1/1/2/3,
	0/2/1/4/5, 0/2/1/5/6, 0/2/2/1/2, 0/2/4/4/5, 0/2/4/5/6, 0/2/5/1/2,
	1/2/1/4/5, 1/2/1/5/6, 1/2/2/1/2, 1/2/4/1/2, 1/2/4/2/3, 1/2/5/3/4, 1/2/5/4/5,
	2/2/1/1/2, 2/2/1/2/3, 2/2/2/3/4, 2/2/2/4/5, 2/2/4/4/5, 2/2/4/5/6, 2/2/5/1/2,
	3/2/1/1/2, 3/2/1/2/3, 3/2/2/3/4, 3/2/2/4/5, 3/2/4/1/2, 3/2/4/2/3, 3/2/5/3/4, 3/2/5/4/5}
\draw[red, thick]
	(B\x\y\a\u) -- (B\x\y\a\v);

\foreach \x/\y/\a/\u/\v in 
	{0/0/0/2/3, 0/0/0/3/4, 0/0/0/5/6, 0/0/1/4/5, 0/0/2/3/4, 0/0/2/4/5, 0/0/3/2/3, 0/0/3/4/5, 0/0/3/5/6, 0/0/4/4/5, 0/0/4/5/6, 0/0/5/4/5, 0/0/5/5/6, 0/0/6/4/5,
	0/1/1/2/3, 0/1/1/4/5, 0/1/5/2/3, 0/1/5/4/5,
	1/1/1/2/3, 1/1/1/4/5, 1/1/5/1/2, 1/1/5/3/4,
	2/1/1/1/2, 2/1/1/3/4, 2/1/5/2/3, 2/1/5/4/5,
	3/1/1/1/2, 3/1/1/3/4, 3/1/5/1/2, 3/1/5/3/4}
\draw[blue, thick]
	(B\x\y\a\u) -- (B\x\y\a\v);

\draw[green, thick]
	;

\foreach \x/\y/\a/\u/\v in {
	0/0/0/0/5, 0/0/0/1/4, 0/0/0/2/3, 0/0/0/3/2, 0/0/0/4/1, 0/0/0/5/0, 0/0/0/6/6,
	0/0/1/0/2, 0/0/1/1/3, 0/0/1/2/4, 0/0/1/3/5, 0/0/1/4/6, 0/0/1/5/0, 0/0/1/6/1,
	0/0/2/0/6, 0/0/2/1/5, 0/0/2/2/4, 0/0/2/3/3, 0/0/2/4/2, 0/0/2/5/1, 0/0/2/6/0,
	0/0/3/0/4, 0/0/3/1/5, 0/0/3/2/6, 0/0/3/3/0, 0/0/3/4/1, 0/0/3/5/2, 0/0/3/6/3,
	0/0/4/0/0, 0/0/4/1/6, 0/0/4/2/5, 0/0/4/3/4, 0/0/4/4/3, 0/0/4/5/2, 0/0/4/6/1,
	0/0/5/0/0, 0/0/5/1/6, 0/0/5/2/5, 0/0/5/3/4, 0/0/5/4/3, 0/0/5/5/2, 0/0/5/6/1,
	0/0/6/0/2, 0/0/6/1/3, 0/0/6/2/4, 0/0/6/3/5, 0/0/6/4/6, 0/0/6/5/0, 0/0/6/6/1,
	2/0/2/0/3, 2/0/2/1/4, 2/0/2/2/5, 2/0/2/3/6, 2/0/2/4/0, 2/0/2/5/1, 2/0/2/6/2,
	2/0/3/0/1, 2/0/3/1/2, 2/0/3/2/3, 2/0/3/3/4, 2/0/3/4/5, 2/0/3/5/6, 2/0/3/6/0,
	2/0/4/0/3, 2/0/4/1/2, 2/0/4/2/1, 2/0/4/3/0, 2/0/4/4/6, 2/0/4/5/5, 2/0/4/6/4,
	3/0/2/0/3, 3/0/2/1/2, 3/0/2/2/1, 3/0/2/3/0, 3/0/2/4/6, 3/0/2/5/5, 3/0/2/6/4,
	3/0/3/0/0, 3/0/3/1/1, 3/0/3/2/2, 3/0/3/3/3, 3/0/3/4/4, 3/0/3/5/5, 3/0/3/6/6,
	3/0/4/0/5, 3/0/4/1/6, 3/0/4/2/0, 3/0/4/3/1, 3/0/4/4/2, 3/0/4/5/3, 3/0/4/6/4,
	0/1/1/0/2, 0/1/1/1/3, 0/1/1/2/4, 0/1/1/3/5, 0/1/1/4/6, 0/1/1/5/0, 0/1/1/6/1,
	0/1/5/0/2, 0/1/5/1/3, 0/1/5/2/4, 0/1/5/3/5, 0/1/5/4/6, 0/1/5/5/0, 0/1/5/6/1,
	1/1/1/0/2, 1/1/1/1/3, 1/1/1/2/4, 1/1/1/3/5, 1/1/1/4/6, 1/1/1/5/0, 1/1/1/6/1,
	1/1/5/0/1, 1/1/5/1/0, 1/1/5/2/6, 1/1/5/3/5, 1/1/5/4/4, 1/1/5/5/3, 1/1/5/6/2,
	2/1/1/0/1, 2/1/1/1/0, 2/1/1/2/6, 2/1/1/3/5, 2/1/1/4/4, 2/1/1/5/3, 2/1/1/6/2,
	2/1/5/0/2, 2/1/5/1/3, 2/1/5/2/4, 2/1/5/3/5, 2/1/5/4/6, 2/1/5/5/0, 2/1/5/6/1,
	3/1/1/0/1, 3/1/1/1/0, 3/1/1/2/6, 3/1/1/3/5, 3/1/1/4/4, 3/1/1/5/3, 3/1/1/6/2,
	3/1/5/0/1, 3/1/5/1/0, 3/1/5/2/6, 3/1/5/3/5, 3/1/5/4/4, 3/1/5/5/3, 3/1/5/6/2}
\node at (X\x\y\a\u) {\v};

\foreach \x in {2,3}
\foreach \a/\u/\v in {
	0/0/4, 0/1/5, 0/2/6, 0/3/0, 0/4/1, 0/5/2, 0/6/3,
	1/0/2, 1/1/3, 1/2/4, 1/3/5, 1/4/6, 1/5/0, 1/6/1,
	5/0/6, 5/1/0, 5/2/1, 5/3/2, 5/4/3, 5/5/4, 5/6/5,
	6/0/4, 6/1/5, 6/2/6, 6/3/0, 6/4/1, 6/5/2, 6/6/3}
\node at (X\x 0\a\u) {\v};

\foreach \x in {0,1,2,3}
{

\foreach \a/\u/\v in {
	0/0/3, 0/1/4, 0/2/5, 0/3/6, 0/4/0, 0/5/1, 0/6/2,
	2/0/1, 2/1/2, 2/2/3, 2/3/4, 2/4/5, 2/5/6, 2/6/0,
	3/0/4, 3/1/5, 3/2/6, 3/3/0, 3/4/1, 3/5/2, 3/6/3,
	4/0/0, 4/1/1, 4/2/2, 4/3/3, 4/4/4, 4/5/5, 4/6/6,
	6/0/5, 6/1/6, 6/2/0, 6/3/1, 6/4/2, 6/5/3, 6/6/4}
\node at (X\x 1\a\u) {\v};

\foreach \a/\u/\v in {
	0/0/0, 0/1/6, 0/2/5, 0/3/4, 0/4/3, 0/5/2, 0/6/1,
	3/0/3, 3/1/2, 3/2/1, 3/3/0, 3/4/6, 3/5/5, 3/6/4,
	6/0/6, 6/1/5, 6/2/4, 6/3/3, 6/4/2, 6/5/1, 6/6/0}
\node at (X\x 2\a\u) {\v};

}

\foreach \x in {0,1}
{
\foreach \a/\u/\v in {
	1/0/2, 1/1/3, 1/2/4, 1/3/5, 1/4/6, 1/5/0, 1/6/1,
	2/0/1, 2/1/0, 2/2/6, 2/3/5, 2/4/4, 2/5/3, 2/6/2}
\node at (X\x 2\a\u) {\v};

\draw[green, thick]
	(B\x 212) -- (B\x 213) -- (B\x 214)
	(B\x 222) -- (B\x 223) -- (B\x 224);
	
}

\foreach \x in {2,3}
{
\foreach \a/\u/\v in {
	1/0/2, 1/1/1, 1/2/0, 1/3/6, 1/4/5, 1/5/4, 1/6/3,
	2/0/3, 2/1/4, 2/2/5, 2/3/6, 2/4/0, 2/5/1, 2/6/2}
\node at (X\x 2\a\u) {\v};

\draw[green, thick]
	(B\x 213) -- (B\x 214) -- (B\x 215)
	(B\x 221) -- (B\x 222) -- (B\x 223);;
}

\foreach \x in {0,2}
{
\foreach \a/\u/\v in {
	4/0/5, 4/1/4, 4/2/3, 4/3/2, 4/4/1, 4/5/0, 4/6/6,
	5/0/6, 5/1/0, 5/2/1, 5/3/2, 5/4/3, 5/5/4, 5/6/5}
\node at (X\x 2\a\u) {\v};

\draw[blue, thick]
	(B\x 242) -- (B\x 243) -- (B\x 244)
	(B\x 252) -- (B\x 253) -- (B\x 254);

}

\foreach \x in {1,3}
{
\foreach \a/\u/\v in {
	4/0/5, 4/1/6, 4/2/0, 4/3/1, 4/4/2, 4/5/3, 4/6/4,
	5/0/4, 5/1/3, 5/2/2, 5/3/1, 5/4/0, 5/5/6, 5/6/5}
\node at (X\x 2\a\u) {\v};

\draw[blue, thick]
	(B\x 243) -- (B\x 244) -- (B\x 245) -- (B\x 252) -- (B\x 253) ;
}

\foreach \c in {0061, 2001, 2041, 2051, 3001, 3021, 3031, 0111, 0121, 0151, 0161, 1111, 1121, 1151, 1161, 2111, 2121, 2111, 2151, 2161, 3111, 3121, 3151, 3161, 0201, 1201, 2201, 3201}
\fill
	(B\c) circle (0.05);

\begin{scope}[yshift=-1.6cm, font=\small]

\node at (0,0) {$7^1_{4+4}\colon \tilde{5}2\tilde{6}4\tilde{0}\tilde{0}2$};
\node at (6,0) {$7^2_{4+4}(1)\colon 4231\tilde{3}64$};
\node at (9,0) {$7^2_{4+4}(2)\colon 42\tilde{3}0564$};
\node at (0,-3.2) {$7^3_{4+4}(1)\colon 3214025$};
\node at (3,-3.2) {$7^3_{4+4}(2)\colon 32140\tilde{1}5$};
\node at (6,-3.2) {$7^3_{4+4}(3)\colon 3\tilde{1}14025$};
\node at (9,-3.2) {$7^3_{4+4}(4)\colon 3\tilde{1}140\tilde{1}5$};
\node at (0,-6.4) {$7_{3+5}(1)\colon \tilde{0}2\tilde{1}\tilde{3}\tilde{5}6\tilde{6}$};
\node at (3,-6.4) {$7_{3+5}(2)\colon \tilde{0}2\tilde{1}\tilde{3}5\tilde{4}\tilde{6}$};
\node at (6,-6.4) {$7_{3+5}(3)\colon \tilde{0}\tilde{2}3\tilde{3}\tilde{5}6\tilde{6}$};
\node at (9,-6.4) {$7_{3+5}(4)\colon \tilde{0}\tilde{2}3\tilde{3}5\tilde{4}\tilde{6}$};

\end{scope}

\end{tikzpicture}
\caption{All ${\mc D}_1$-patches of rank 8.}
\label{H8A}
\end{figure}

We use computer program to construct ${\mc D}_1$-patches for the edge and corner systems in the updated list. We find no ${\mc D}_1$-patches for any of the systems of rank 7 (indicated by three 0 in the sixth column of Table \ref{table3}), and find exactly four edge and corner systems of rank 8 admitting ${\mc D}_1$-patches:
\begin{align*}
7^1_{4+4} &\colon
|\bar{0}|=|\bar{4}|=|\bar{5}|,\;
|\bar{1}|=|\bar{2}|=|\bar{6}|,\;
[5]=\tfrac{2\pi}{3},\;
[0]+[1]+[5]=[0]+[3]+[4]=[1]+[2]+[6]=2\pi; \\
7^2_{4+4} &\colon
|\bar{0}|=|\bar{2}|=|\bar{3}|=|\bar{4}|=|\bar{6}|,\;
[0]=\tfrac{2\pi}{3},\;
[0]+[3]+[4]=[1]+[2]+[3]=[4]+[5]+[6]=2\pi; \\
7^3_{4+4} &\colon
|\bar{0}|=|\bar{1}|=|\bar{2}|=|\bar{5}|,\;
|\bar{4}|=|\bar{6}|,\;
[1]=[2]=\tfrac{2}{3}\pi,\;
[0]+[3]+[4]=[1]+[5]+[6]=2\pi; \\
7_{3+5} &\colon
|\bar{0}|=|\bar{6}|,\;
|\bar{1}|=|\bar{2}|,\;
|\bar{4}|=|\bar{5}|,\;
[0]=\tfrac{2\pi}{3},\;
[1]=[3],\;
[4]=[6],\;
[1]+[2]+[3]=[4]+[5]+[6]=2\pi.
\end{align*}
We find generic heptagons, and completions of the corner parts are the following
\begin{align*}
\bar{V}(7^1_{4+4}) &=\{ 015, 034, 126, 555\}; \\
\bar{V}(7^2_{4+4}) &=\{ 034, 123, 456, 000\}; \\
\bar{V}(7^3_{4+4}) &=\{ 034, 156, 256, 111, 112, 122, 222\}; \\
\bar{V}(7_{3+5}) &=\{ 112, 123, 233, 445, 456, 566, 000\}.
\end{align*}
Based on these data, we construct all the corresponding ${\mc D}_1$-patches in Figure \ref{H8A}. For $7^2_{4+4}$, we may horizontally flip the yellow part. For $7^3_{4+4}$ and $7_{3+5}$, each of the two yellow parts can have two independent choices.

On the boundary of the first ${\mc D}_1$-patch $\tilde{5}2\tilde{6}4\tilde{0}\tilde{0}2$ for $7^1_{4+4}$, the vertex $\bullet=23\cdots$ does not belong to $\bar{V}(7^1_{4+4})$, i.e., there is no $k\in {\bb Z}_7$, such that $23k\in \bar{V}(7^1_{4+4})$. Therefore the ${\mc D}_1$-patch cannot be extended to a ${\mc D}_2$-patch. We call such argument {\em $\partial {\mc D}_1$-verification}: There are seven vertices on the boundary of a ${\mc D}_1$-patch, in which we already know two corners $i$ and $j$. The extension of the ${\mc D}_1$-patch to a ${\mc D}_2$-patch requires us to find a corner $k$, such that $ijk\in\bar{V}$. All the ${\mc D}_1$-patches in Figure \ref{H8A} fail the $\partial {\mc D}_1$-verification at $\bullet$. This is indicated by two 0 in the seventh column of Table \ref{table3}.

We remark that there is similar {\em $\partial {\mc D}_k$-verification} along the boundary of ${\mc D}_k$-patch that is necessary for extension to ${\mc D}_{k+1}$-patch. 

Among the edge and corner systems of rank 9 listed in the fifth column of Table \ref{table3}, we find $11+95+7=113$ admitting ${\mc D}_1$-patches. Among these, $13+4=17$ pass $\partial {\mc D}_1$-verification. However, only one has further extension to ${\mc D}_2$-patches:
\begin{equation}\label{h9ev}
7_{3+6}\colon
|\bar{0}|=|\bar{6}|,\;
|\bar{2}|=|\bar{3}|=|\bar{4}|,\;
[2]=[3]=[4]=[5]=\tfrac{2}{3}\pi,\;
[1]=[6],\;
[0]+[1]+[6]=2\pi.
\end{equation}
The following geometrical argument shows the system implies the system \eqref{h5ev}. Therefore we do not include the system in Theorem \ref{thm1}.

\begin{proposition}\label{h9hidden}
If a hyperbolic heptagon satisfies \eqref{h9ev}, then $|\bar{1}|=|\bar{5}|$. In particular, the heptagon is symmetric and satisfies \eqref{h5ev}.
\end{proposition}

\begin{proof}
The left of Figure \ref{H9A} is a heptagon satisfying \eqref{h9ev}. The two red edges have the same lengths, and the three blue edges have the same lengths. We will show that the heptagon is symmetric with respect to the horizontal flip. This implies the two normal edges have the same length.  

By equal edge lengths, and $[3]=[4]$, we know the triangle $\triangle 016$ and the quadrilateral $\square 2345$ are symmetric. Then by $[1]=[6]$ and $[2]=[5]$, we get a quadrilateral $\square 1256$ with updated corner labels $1',2',5',6'$ satisfying $[1']=[6']$ and $[2']=[5']$.

\begin{figure}[h!]
\centering
\begin{tikzpicture}[>=latex, scale=1]

\draw[dashed]
	(-1,0.5) -- (1,0.5)
	(-1.2,-0.5) -- (1.2,-0.5);

\foreach \x in {1,-1}
{
\begin{scope}[xscale=\x]

\draw[red, thick]
	(0,1) -- (1,0.5);

\draw
	(1.2,-0.5) -- (1,0.5);

\draw[blue, thick]
	(1.2,-0.5) -- (0.5,-1.2) -- (0,-1.2);
	
\end{scope}
}

\foreach \x/\y/\a in 
	{0/0.8/0, -0.85/0.35/1, 0.85/0.35/6, -1/-0.4/2, 1/-0.4/5, -0.4/-1/3, 0.4/-1/4}
\node at (\x,\y) {\footnotesize \a};

\begin{scope}[xshift=5cm]

\foreach \x in {1,-1}
{
\begin{scope}[xscale=\x]

\fill[gray!50]
	(1.3,-1) -- (0.9,1) -- (1.4,1) -- (1.8,-1);
	
\draw
	(1.3,-1) -- (0.9,1);

\end{scope}
}	
	
\draw[dashed]
	(-1,0.5) -- (1,0.5)
	(-1.2,-0.5) -- (1.2,-0.5);

\draw
	(0,-0.1) -- (0,-0.9)
	(0,0.1) -- (0,0.9);

\foreach \x/\y/\a in 
	{-0.85/0.35/1, 0.85/0.35/6, -1/-0.35/2, 1/-0.35/5}
\node at (\x,\y) {\footnotesize \a'};

\node at (-1.3,-1.2) {\small $L_{12}$};
\node at (1.3,-1.2) {\small $L_{56}$};

\node at (0,-1.1) {\small $R_{25}$};
\node at (0,1.1) {\small $R_{16}$};

\end{scope}

\end{tikzpicture}
\caption{Proof of Proposition \ref{h9hidden}.}
\label{H9A}
\end{figure}
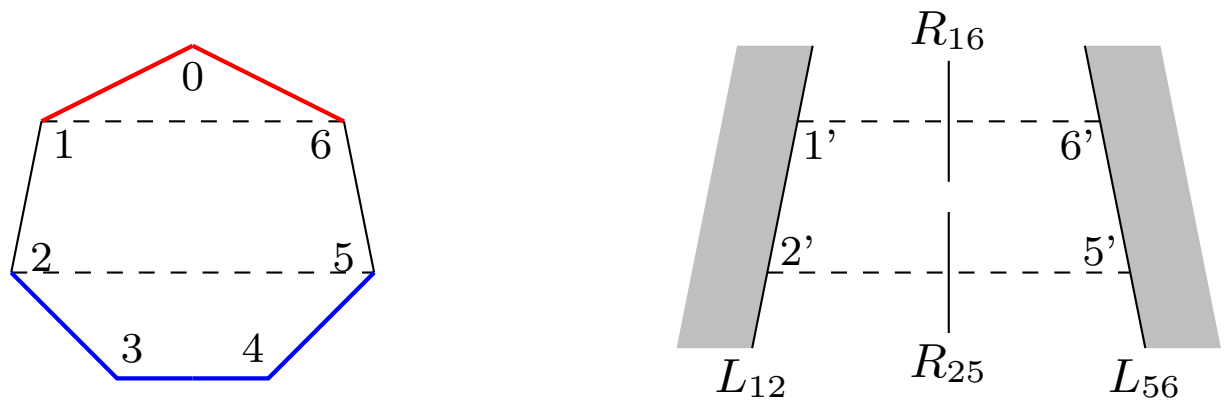

On the right of Figure \ref{H9A} is the quadrilateral $\square 1256=\square 1'2'5'6'$. We extend the edges $\bar{1}$ and $\bar{5}$ to geodesics $L_{12}$ and $L_{56}$, and also choose the sides of the two geodesics indicated by the gray regions. We consider the reflection $R_{16}$ that exchanges the two ends of the dashed line connecting $1'$ and $6'$, and the other reflection $R_{25}$ that exchanges the two ends of the dashed line connecting $2'$ and $5'$.

By $[1']=[6']$, the reflection $R_{16}$ exchanges $L_{12}$ and $L_{56}$, together with their chosen sides. By $[2']=[5']$, the reflection $R_{25}$ has the same property. It is an elementary fact of abstract geometry (including hyperbolic geometry) that the reflection with the property is unique. Therefore the two reflections are the same reflection. This reflection is the symmetry for the quadrilateral $\square 1256=\square 1'2'5'6'$ as well as the heptagon.
\end{proof}

\section{Ranks 10, 11, 12}
\label{rank11}

Among edge and corner systems of ranks 10 and 11 in the fourth column of Table \ref{table3}, we remove those that imply \eqref{h5ev}, \eqref{h6ev}, or \eqref{h9ev} (due to Proposition \ref{h9hidden}) through linear algebra, up to rotation and flip of labels. Then the number of systems of ranks 10 and 11 are updated in the fifth column of Table \ref{table3}. 

Among the edge and corner systems of rank 10 in the updated list, we find $8+375+674+1=1058$ admitting ${\mc D}_1$-patches. Applying $\partial {\mc D}_1$-verification reduces the number to $35+223+1=259$. Then we find only one system admitting ${\mc D}_2$-patches. The system survives $\partial {\mc D}_2$-verification, but does not admit ${\mc D}_3$-patch. We conclude that no systems of rank 10 in the updated list admit tilings of ${\bb H}^2$.

For ranks 11 and 12, a new tool can be used. A heptagon can be determined by 6 edge lengths and 5 angle values between these 6 edges. In principle, an edge and corner systems of rank $\ge 11$ should determine isolated cases of heptagons. Of course the regular heptagon $P_7(\frac{2\pi}{3})$ is always a solution. The problem is whether there are non-regular heptagons satisfying the system. 

The calculation of heptagons satisfying 11 conditions is very intensive. The construction of ${\mc D}_k$-patches based on the edge and corner systems is also computationally intensive as $k$ becomes larger. Therefore we first construct ${\mc D}_k$-patches for increasing $k$, until the number of systems admitting ${\mc D}_k$-patches becomes small. Then we symbolically calculate the heptagons satisfying these systems. Section \ref{calculation} explains the equation \eqref{hpoly} underlying the calculation, gives the detailed calculation for  the non-regular $7^1_{5+6}$-heptagon, and outlines the calculation for the non-regular $7^2_{5+6}$-heptagon.

Among the systems of rank 11 in the updated list, we find $249+3276+1=3526$ admitting ${\mc D}_1$-patches. The last one edge and corner system with $\text{rank}E=4$ and $\text{rank}V=7$ is the following
\[
|\bar{0}|=|\bar{1}|,\;
|\bar{2}|=|\bar{3}|=|\bar{4}|,\;
|\bar{5}|=|\bar{6}|,\;
[0]=[1]=[2]=[3]=[4]=[5]=[6]=\tfrac{2\pi}{3}.
\]
The equation \eqref{hpoly} is three polynomial equations of the hyperbolic sine and cosine of $k=|\bar{0}|$, $l=|\bar{2}|$, $m=|\bar{5}|$. Changing the variables to $x=e^k$, $y=e^l$, $z=e^m$, we get three polynomial equations in $x,y,z$. Then we symbolically calculate the solutions, and find that the only solution is the regular heptagon.

We apply $\partial {\mc D}_1$-verification to the other $3525$ systems and reduce the number to $19+810=829$. Then we further find $1+20=21$ systems admitting ${\mc D}_2$-patches. The only system with $\text{rank}E=6$ and $\text{rank}V=5$ cannot pass the $\partial {\mc D}_2$-verification. Moreover, 16 systems with $\text{rank}E=5$ and $\text{rank}V=6$ pass the $\partial {\mc D}_2$-verification. Among these remaining 16, only 6 edge and corner systems admit ${\mc D}_3$-patches. The 6 systems are \eqref{h11ev1}, \eqref{h11ev2}, and 4 systems with the common edge equalities 
\[
E\colon |\bar{1}| = |\bar{2}| = |\bar{3}| = |\bar{4}| = |\bar{5}| = |\bar{6}|,
\]
and the following vertex equalities
\begin{itemize}
\item $[0] + [0] + [5]  
= [0] + [1] + [2] 
= [1] + [1] + [6]
= [2] + [3] + [5] 
= [3] + [4] + [4]
= [4] + [6] + [6] 
= 2\pi$.
\item $[0] = [6]$, 
$[1] = [2]$, 
$[0] + [0] + [5] 
= [1] + [1] + [3] 
= [2] + [4] + [6]
= [3] + [4] + [5]  = 2\pi$.
\item $[1] = [3] = [6] = \tfrac{2}{3}\pi$, 
$[0] + [1] + [4] 
= [0] + [2] + [2]
= [4] + [5] + [5] 
=  2\pi$.
\item $[1] = \tfrac{2\pi}{3}$, 
$[3] = [5]$, 
$[4] = [6]$, 
$[0] + [1] + [2]
= [0] + [3] + [3]
= [1] + [3] + [6]
= [2] + [4] + [4]  = 2\pi$.
\end{itemize}
The 4 systems have two free edge lengths $k,l$ and one free angle $\theta$. The equation \eqref{hpoly} is three polynomial equations of the hyperbolic sine and cosine of $k,l$, and usual sine and cosine of $\theta$. Changing the variables to $x=e^k$, $y=e^l$, and $t=\tan\frac{\theta}{2}$, we get three polynomial equations in $x,y,t$. Then we symbolically calculate the solutions, and find that the only solution is the regular heptagon.

Finally, edge and corner systems of rank 12 are overdetermined, and we expect only regular heptagon solutions. Still, we follow the same process as above. We find $2214$ systems admitting ${\mc D}_1$-patches, all of which have $\text{rank}E=6$ and $\text{rank}V=6$. Then $414$ pass $\partial {\mc D}_1$-verification, 86 admit ${\mc D}_2$-patches, 44 pass $\partial {\mc D}_2$-verification, 32 admit ${\mc D}_3$-patches, and 31 pass $\partial {\mc D}_3$-verification. Then we symbolically calculate the heptagons for the 31 systems, and find regular heptagon is the only solution.

\subsection{Calculation of Heptagons}
\label{calculation}

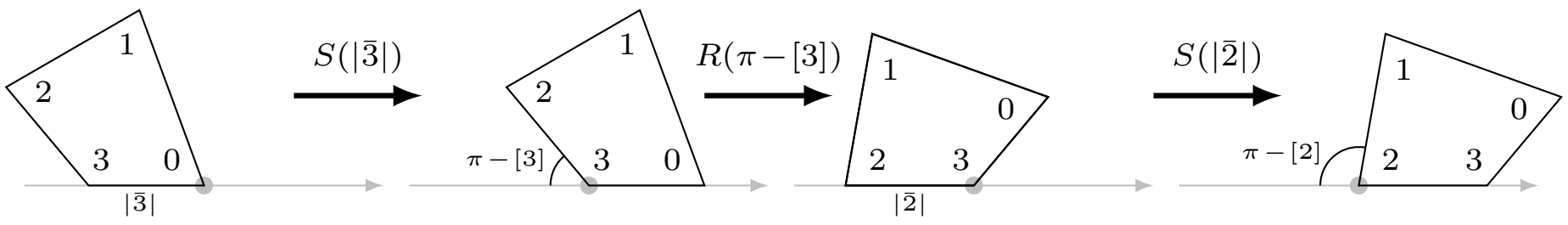
\begin{figure}[htp]
\centering
\begin{tikzpicture}[>=latex]

\foreach \x in {0,...,3}
{
\draw[->, gray!50]
	(-1.4+3*\x,0) -- ++(2.8,0);

\fill[gray!50] (3*\x,0) circle (0.07);
}

\draw[xshift=-0.9 cm]
	(0,0) -- ++(130:1) -- ++(30:1.2) -- (0:0.9) -- cycle;

\draw[xshift=3 cm]
	(0,0) -- ++(130:1) -- ++(30:1.2) -- (0:0.9) -- cycle;
		
\draw[xshift=6cm, rotate=50]
	(0,0) -- ++(130:1) -- ++(30:1.2) -- (0:0.9) -- cycle;

\draw[xshift=6cm, rotate=50]
	(0,0) -- ++(130:1) -- ++(30:1.2) -- (0:0.9) -- cycle;

\draw[xshift=10cm, rotate=50]
	(0,0) -- ++(130:1) -- ++(30:1.2) -- (0:0.9) -- cycle;

\draw[xshift=3 cm]
	(-0.3,0) arc (180:130:0.3);
\draw[xshift=9 cm]
	(-0.3,0) arc (180:80:0.3);
	
\foreach \x/\y/\n in {
	-0.25/0.2/0, -0.6/1.1/1, -1.25/0.73/2, -0.8/0.2/3, 
	3.65/0.2/0, 3.3/1.1/1, 2.65/0.73/2, 3.1/0.2/3, 
	6.25/0.6/0, 5.35/0.9/1, 5.25/0.2/2, 5.9/0.2/3, 
	10.25/0.6/0, 9.35/0.9/1, 9.25/0.2/2, 9.9/0.2/3
	}	
\node at (\x,\y) {\scriptsize \n};

\node at (-0.5,-0.15) {\tiny $|\bar{3}|$};
\node at (2.35,0.2) {\tiny $\pi\!-\![3]$};
\node at (5.5,-0.15) {\tiny $|\bar{2}|$};
\node at (8.4,0.25) {\tiny $\pi\!-\![2]$};

\draw[->, very thick]
	(0.7,0.7) -- ++(1,0);
\node at (1.2,1) {\scriptsize $S(|\bar{3}|)$};

\draw[->, very thick]
	(3.9,0.7) -- ++(1,0);
\node at (4.4,1) {\scriptsize $R(\pi\!-\![3])$};

\draw[->, very thick]
	(7.4,0.7) -- ++(1,0);
\node at (7.9,1) {\scriptsize $S(|\bar{2}|)$};
		
\end{tikzpicture}
\caption{A polygon is a sequence of shifts and rotations.}
\label{polygon}
\end{figure}

The governing equation for a hyperbolic $n$-gon is
\begin{equation}\label{hpoly}
\prod_{i=0}^{n-1}R(\pi-[i])S(|\bar{i}|)=I.
\end{equation}
The equation appeared in Section 3.3 of \cite{cly1} and Section 3.4 of \cite{cly2}. We illustrate the equation by the quadrilateral in Figure \ref{polygon}. The gray dot is the center of the Poincar\'e disk ${\bb H}^2$, and the arrow is horizontal and to the right. In the first picture, we put the quadrilateral in the position such that the corner 0 is the center of the Poincar\'e disk, and the edge $\bar{3}$ is horizontal and on the left of the center. Next, applying the horizontal shift $S(|\bar{3}|)$ that moves the quadrilateral to the right by the length $|\bar{3}|$, we get the second picture, where the corner 3 is at the center, and the edge $\bar{3}$ is still horizontal. Then we apply the rotation $R(\pi-[3])$ around the center to get the third picture, where the corner 3 is at the center, and the edge $\bar{2}$ is horizontal and on the left of the center. Next, we apply the horizontal shift $S(|\bar{2}|)$, and so on. All the shifts and rotations that correspond to the edges and corners of the quadrilateral combine to form the transformation $R(\pi-[0])S(|\bar{0}|)R(\pi-[1])S(|\bar{1}|)R(\pi-[2])S(|\bar{2}|)R(\pi-[3])S(|\bar{3}|)$. Since applying the transformation to the quadrilateral gives back the quadrilateral at the beginning, the composition is the identity.

The isometries of the Poincar\'e disk are M\"obius transforms. In terms of the complex numbers in the unit disk, the two transformations are
\[
S(u)\colon z\mapsto \frac{z\cosh\frac{u}{2}+\sinh\frac{u}{2}}{z\sinh\frac{u}{2}+\cosh\frac{u}{2}};\quad
R(v)\colon z\mapsto e^{iv}z.
\]
In practice, it is easier to use the corresponding matrices in $PSL_2({\bb R})$ (obtained by using the half plane model instead of the disk model for ${\bb H}^2$)
\begin{equation}\label{matrix}
S(u)=\begin{pmatrix}
e^{\frac{u}{2}} & 0 \\
0 & e^{-\frac{u}{2}}
\end{pmatrix},\quad
R(v)=\begin{pmatrix}
\cos \frac{v}{2} & -\sin \frac{v}{2} \\
\sin \frac{v}{2} & \cos \frac{v}{2}
\end{pmatrix}.
\end{equation}
Note that equalities in $PSL_2({\bb R})$ are $2\times 2$ matrix equalities up to the multiplication by $\pm I$.

Since the dimension of the isometry group of the Poincar\'e disk is 3, the $2\times 2$ matrix equation \eqref{hpoly} is actually a system of 3 equations. Moreover, the equation essentially means that line segments of certain length and the turning angles determine a closed loop. For the closed loop to represent a polygon suitable for tiling, we still need to check that the closed loop has no self intersection.

Next, we outline the calculation for $7^1_{5+6}$-heptagon. In the edge and corner system in \eqref{h11ev1}, we set $a=|\bar{0}|$, $b=|\bar{1}|$, and $\theta=[2]-\frac{2\pi}{3}$. Then all other edge lengths and angle values can be expressed in terms of $a,b,\theta$. The condition that all angles belong to the interval $(0,2\pi)$ mean $-\frac{2\pi}{15}<\theta<\frac{\pi}{9}$.

The $2\times 2$ matrix equation \eqref{hpoly} for the heptagon can be written as 
\[
S(b)MS(b)N=\sigma I,\quad
\sigma=\pm,
\]
where $M$ and $N$ do not involve $b$
\begin{align*}
M &=R(\pi-[2])S(a)R(\pi-[3])S(a)R(\pi-[4])S(a)R(\pi-[5]), \\
N &=R(\pi-[6])S(a)R(\pi-[0])S(a)R(\pi-[1]).
\end{align*}
We wish to eliminate $b$. Since $S(b)=\begin{pmatrix}\lambda & 0 \\ 0 & \lambda^{-1} \end{pmatrix}$, we know the off diagonal terms of $S(b)MS(b)$ and $M$ are the same. Moreover, since $\det N=1$, we know the off diagonal terms of $N^{-1}$ and $N$ differ by a sign. Comparing the off diagonal terms of $S(b)MS(b)$ and $\sigma N^{-1}$ gives equations involving only $a$ and $\theta$
\[
A=M_{12}+\sigma N_{12}=0,\quad
B=M_{21}+\sigma N_{21}=0.
\]
Calculation shows 
\[
A-B=2\sigma F(\theta,c),\quad
A+B=2\sinh\frac{a}{2}\;\sin\theta \;G(\theta,c),
\]
where $c=-\sigma\cosh\frac{a}{2}$, and
\begin{align*}
F(\theta,c)
&=c^3\left[\cos(\tfrac{\pi}{6}+\theta)+\sqrt{3}\cos 3\theta-\sin \theta\right]
-c^2\left[\sin(\tfrac{\pi}{6}-7\theta)+\cos \theta\right] \\
&\quad -c\left[\sqrt{3}\cos 3\theta-\sin \theta\right]
+\sin(\tfrac{\pi}{6}-7\theta), \\
G(\theta,c)
&=c^2\left[\cos \theta+\sqrt{3}\sin \theta-4\cos \theta\cos 2\theta\right]
-2c\cos(\tfrac{\pi}{6}+3\theta)
+4\cos \theta\cos 2\theta.
\end{align*}
Since $a>0$, the equations $A=0$ and $B=0$ become two problems:
\begin{enumerate}
\item $\theta=0$ and $F(0,c)=0$.
\item $F(\theta,c)=0$ and $G(\theta,c)=0$.
\end{enumerate}

If $\theta=0$, then $[i]=\frac{2\pi}{3}$ for all $i\in{\bb Z}_7$, and $F(0,c)=0$ has three roots $c=\frac{2}{\sqrt{3}}\cos\frac{k\pi}{7}$, where $k=1,3,5$. Since $c=-\sigma\cosh\frac{a}{2}$ has absolute value $>1$ and only the first root has this property, we get $c=\frac{2}{\sqrt{3}}\cos\frac{k\pi}{7}$, and $\sigma=-1$. This gives $a$, and we may further calculate $b$. The result shows the heptagon is regular.

For the equation $F(\theta,c)=G(\theta,c)=0$, we note that $F$ and $G$ are polynomials of $c$ with trigonometric functions of $\theta$ as coefficients. In principle, we may get an equation for $\theta$ by calculating the resultant of $F$ and $G$ in $c$. A more efficient way is to express the trigonometric function coefficients as polynomials of $t=\frac{\tan \theta}{\sqrt{3}}$, up to dividing a power of $\cos \theta$. Specifically, we find that $\tilde{F}(t,c)=\frac{F(\theta,c)}{\cos^7\theta}$ and $\tilde{G}(t,c)=\frac{G(\theta,c)}{\cos^3\theta}$ are polynomials of $c$, with polynomials of $t$ as coefficients. The range $-\frac{2\pi}{15}<\theta<\frac{\pi}{9}$ implies that $\cos \theta>0$ and $|t|<\frac{1}{3}$. Therefore $F(\theta,c)=G(\theta,c)=0$ is the same as $\tilde{F}(t,c)=\tilde{G}(t,c)=0$.

The resultant
\[
\text{Res}_c(\tilde{F},\tilde{G})
=(t+1)^3(3t-1)^2p(t)q(t),
\]
where $p(t)=9t^3-9t^2-9t+1$, and $q(t)$ is a polynomial of degree 15. We may show $q(t)<0$ for $|t|<\frac{1}{3}$. One way is to use the substitution $t=\frac{s-1}{3(s+1)}$, with $s>0$, to find that $(s+1)^{15}q(\frac{s-1}{3(s+1)})$ is a polynomial with all negative coefficients. Since $|t|<\frac{1}{3}$ also implies $t+1>0$ and $3t-1<0$, we conclude $p(t)=0$.

We substitute $t=\frac{\tan \theta}{\sqrt{3}}$ back, and get 
\[
p\left(\frac{\tan \theta}{\sqrt{3}}\right)
=\frac{2\cos(\frac{\pi}{3}+3\theta)}{\cos^3\theta}.
\]
Therefore we get one solution $\theta=\frac{\pi}{18}$ of $p(t)=0$. Then we get $c=\cosh\frac{a}{2}>1$ as the only positive root (the other solution is negative) of the quadratic polynomial
\[
G(\tfrac{\pi}{18},c)
=-4c^2\sin\tfrac{2\pi}{9}\cos\tfrac{\pi}{9}-c+4\cos\tfrac{\pi}{9}\cos\tfrac{\pi}{18}.
\]
The solution is $c=\frac{\cos\frac{\pi}{18}}{\cos\frac{\pi}{9}}
\approx 1.04801$. This is the value of $\cosh\frac{a}{2}$. Further calculation shows the value of $b$ is related to the value of $a$ by $\sinh\frac{a}{2}=\tanh\frac{b}{2}$.

In the edge and corner system \eqref{h11ev2} for $7^2_{5+6}$-heptagon, we set $a=|\bar{0}|$, $b=|\bar{3}|$, and $\theta=\frac{2\pi}{3}-[3]$. Since only one edge has length $b$, the equation \eqref{hpoly} is of the form
\[
S(b)M=\sigma I,\quad
\sigma=\pm,
\]
where $M$ does not involve $b$. Since $S(b)=\begin{pmatrix}\lambda & 0 \\ 0 & \lambda^{-1} \end{pmatrix}$, the equation implies $M_{12}=M_{21}=0$. The rest of the calculation is similar to the calculation for $7^1_{5+6}$-heptagon. The resultant we get is a constant multiple of
\[
(t+1)^4(3t-1)^3(3t+1)^2(1+3t^2)^5p(t),
\]
where $t=\frac{\tan\frac{\theta}{2}}{\sqrt{3}}$ and 
\begin{equation}\label{712solution}
p(t)=9t^5+81t^4+30t^3-18t^2-7t+1.
\end{equation}
The condition for all angles to be positive is $-\frac{\pi}{3}<\theta<\frac{\pi}{3}$. This means $|t|<\frac{1}{3}$. Therefore the resultant vanishes if and only if $p(t)=0$. Then we can show that $p(t)$ has unique root in the interval $(-\frac{1}{3},\frac{1}{3})$. The geometric data of the non-regular $7^2_{5+6}$-heptagon can be expressed in terms of this unique root.

\subsection{Tilings by $7^1_{5+6}$-Heptagon}
\label{H111}

By the calculation in Section \ref{calculation}, the non-regular $7^1_{5+6}$-heptagon has the following specific values
\begin{align}
7^1_{5+6}&\colon 
|\bar{0}|=|\bar{2}|=|\bar{3}|=|\bar{4}|=|\bar{6}|=2\text{arccosh}\!~c\approx 0.61729,\;
|\bar{1}|=|\bar{5}|=2\text{arctanh}\sqrt{c^2-1}\approx 0.64900; \nonumber \\
&\quad
[0]=\tfrac{\pi}{3},\;
[1]=\tfrac{17\pi}{18},\;
[2]=\tfrac{13\pi}{18},\;
[3]=\tfrac{4\pi}{9},\;
[4]=\tfrac{7\pi}{9},\;
[5]=\tfrac{11\pi}{18},\;
[6]=\tfrac{5\pi}{6}. \label{111heptagon}
\end{align}
Here $c=\frac{\cos\frac{\pi}{18}}{\cos\frac{\pi}{9}}
\approx 1.04801$ is from Section \ref{calculation}. Using the values, we find the corresponding complete corner system
\begin{equation}\label{h111v}
\bar{V}(7^1_{5+6})
=\{012, 066, 135, 236, 344, 455,
0025, 0334, 0355,
00333,
000000\}.
\end{equation}
For the purpose of this paper, we do not need to consider vertices of degree $\ge 3$ in \eqref{h111v}. Still, we prove the stronger result. 

\begin{proposition}\label{h111tiling}
A tiling by the non-regular $7^1_{5+6}$-heptagon is edge-to-edge, and every vertex has degree $3$.
\end{proposition}

\begin{proof}
We first prove the tiling is edge-to-edge. We mainly use the facts that a half vertex consists of three 0 corners (by the complete corner system \eqref{h111v}), and that corner 0 is bounded by two edges of length $a=2\text{arccosh}\!~c$. 

On the left of Figure \ref{H111X}, we consider a half vertex $H$ on one side $AB$ of a tile \circled{1}. The vertex consists of three 0 corners, two of which belong to the tiles \circled{2} and \circled{3}. The data \eqref{111heptagon} shows $b=2\text{arctanh}\sqrt{c^2-1}$ satisfies $a<b<2a$. Therefore the length $|AB|<2a$. This implies either $|AH|<a$ or $|BH|<a$. Without loss of generality, we assume $|AH|<a$. The inequality $|AH|<a$ implies $A$ is a half vertex on a side of the tile \circled{2}. Then the corner of \circled{1} at the half vertex $A$ is 0, which implies $|AB|=a$. Therefore $|BH|<|AB|=a$, and $B$ is also a half vertex on a side of the tile \circled{3}. Then the corner of \circled{1} at the half vertex $B$ is also 0. Now \circled{1} has two 0 corners at $A$ and $B$. This is a contradiction, and proves the tiling is edge-to-edge.

\begin{figure}[htp]
\centering
\begin{tikzpicture}[>=latex, scale=1, font=\footnotesize]

\draw
	(0,0) -- (60:1.4)
	(0,0) -- (120:1.4)
	(-2,0) -- (2,0)
	(-1.4,0) -- ++(-60:1.4)
	(1,0) -- ++(-120:1.4);

\draw[dashed]
	(60:1.4) to[out=0, in=90] (2,0)
	(120:1.4) to[out=180, in=90] (-2,0)
	(-0.7,-1.2) to[out=-30, in=210] (0.3,-1.2);

\fill
	(0.1,0) arc (0:180:0.1)
	(-1.3,0) arc (0:-180:0.1)
	(1.1,0) arc (0:-180:0.1);

\foreach \x/\y/\u in {0/0.3/0, -0.25/0.15/0, 0.25/0.15/0, -1.15/-0.15/0, 0.75/-0.15/0, 0/-0.17/H, -1.4/0.17/A, 1/0.17/B}
\node at (\x,\y) {$\u$};

\foreach \x/\y/\a in {-0.2/-0.7/1, -1/0.7/2, 1/0.7/3}
\node[inner sep=0.5, draw, shape=circle] at (\x,\y) {\a};

\begin{scope}[shift={(6,0)}]

\draw
	(0,0) -- (-30:1.4) -- ++(30:1.4)
	(0,0) -- (210:1.4) -- ++(150:1.4)
	(60:1.4) -- (0,0) -- (120:1.4)
	(-30:1.4) -- ++(0,-0.7)
	(210:1.4) -- ++(0,-0.7);

\fill
	(0,0) circle (0.1);

\foreach \x/\y/\u in {0/-0.25/0, 1.05/-0.8/6, -1.05/-0.8/1}
\node at (\x,\y) {$\u$};

\foreach \x/\y/\u in {0.55/-0.5/a, 1.05/-1.2/b, 1.85/-0.5/a, 1.2/-0.5/i, 0.25/0.05/j}
\foreach \v in {1,-1}
\node at (\x*\v,\y) {$\u$};

\foreach \x/\y/\a in {0/-1/1, -1.2/0.3/2, 1.2/0.3/3}
\node[inner sep=0.5, draw, shape=circle] at (\x,\y) {\a};
	
\end{scope}
			
\end{tikzpicture}
\caption{Vertices for tilings by non-regular $7^1_{5+6}$-heptagon.}
\label{H111X}
\end{figure}
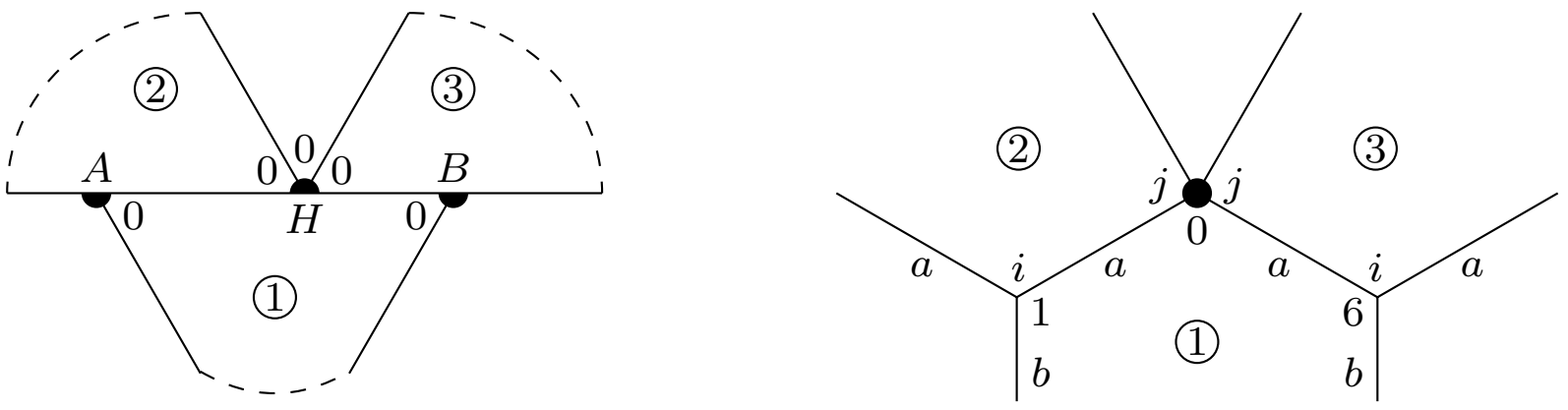

Next, we prove all vertices in the tiling has degree $3$. We note that all vertices of degree $>3$ in the complete corner system \eqref{h111v} contains corner 0. We consider corner 0 at a vertex $\bullet$ on the right of Figure \ref{H111X}. The corner 0 belongs to a tile \circled{1}, in which we also indicated corners 1 and 6, and some edge lengths. In \eqref{h111v}, a vertex $1\cdots=012,135$, and a vertex $6\cdots=066,236$. Moreover, both $1\cdots$ and $6\cdots$ consist of two $a$ and one $b$. Then we conclude the corner $i$ bounded by two $a$ at the two vertices are $i=0$ or 3. This implies the corners of \circled{2} and \circled{3} at $\bullet$ are $j=1,6,2,4$. However, none of the vertices of degree $>3$ in \eqref{h111v} contains two such $j$. The contradiction shows $\bullet$ cannot have degree $>3$.
\end{proof}

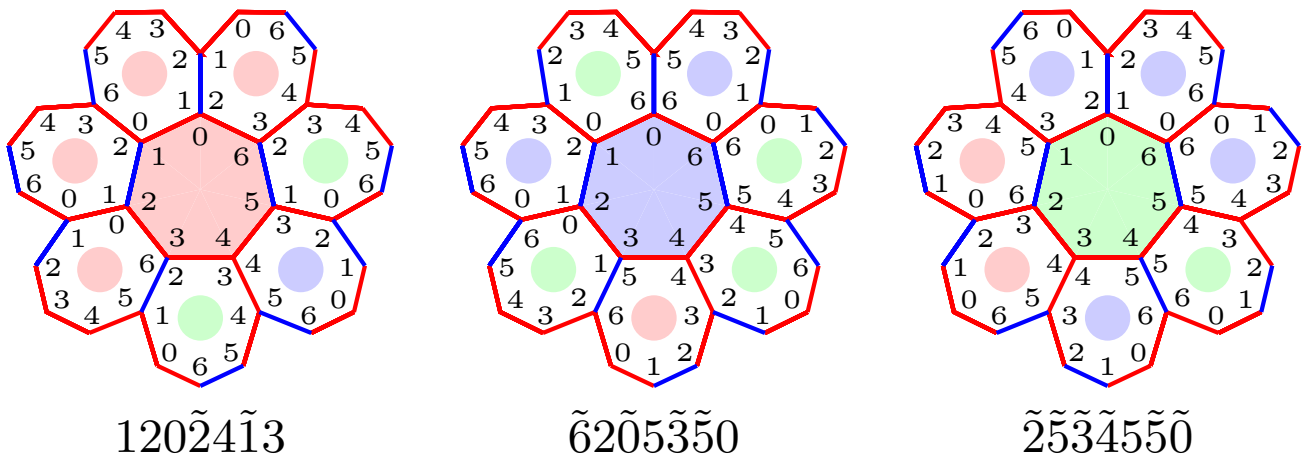
\begin{figure}[htp]
\centering
\begin{tikzpicture}[>=latex, scale=1, thick, font=\tiny]

\foreach \a in {0,...,6}
\foreach \b in {0,...,6}
\foreach \x in {0,1,2}
{
\pgfmathtruncatemacro{\patchx}{ifthenelse(\x==1,2,ifthenelse(\x==2,1,0))}
\begin{scope}[shift={(3*\patchx, 0)}, xscale=-1]

\coordinate (A\x\a) at (90+51.43*\a:0.5);
\coordinate (B\x\a 0) at (90+51.43*\a:0.5);
\coordinate (B\x\a 1) at (90+51.43*\a:0.9);
\coordinate (B\x\a 2) at (102+51.43*\a:1.2);
\coordinate (B\x\a 3) at (115.714+51.43*\a:1.3);
\coordinate (B\x\a 4) at (129.43+51.43*\a:1.2);
\coordinate (B\x\a 5) at (141.43+51.43*\a:0.9);
\coordinate (B\x\a 6) at (141.43+51.43*\a:0.5);

\coordinate (Z\x\a) at (115.7+51.43*\a:0.85);
\coordinate (X\x\a 0) at (100+51.43*\a:0.6);
\coordinate (X\x\a 1) at (99+51.43*\a:0.87);
\coordinate (X\x\a 2) at (105+51.43*\a:1.1);
\coordinate (X\x\a 3) at (115.714+51.43*\a:1.17);
\coordinate (X\x\a 4) at (126.43+51.43*\a:1.1);
\coordinate (X\x\a 5) at (132.43+51.43*\a:0.87);
\coordinate (X\x\a 6) at (131.43+51.43*\a:0.6);

\end{scope}
}

\foreach \a in {0,...,6}
{
\fill[red!20, xscale=-1, rotate=51.43*\a]
	(0,0) -- (90:0.5) -- (90+51.43:0.5);
\fill[blue!20, xshift=3cm, xscale=-1, rotate=51.43*\a]
	(0,0) -- (90:0.5) -- (90+51.43:0.5);
\fill[green!20, xshift=6cm, xscale=-1, rotate=51.43*\a]
	(0,0) -- (90:0.5) -- (90+51.43:0.5);
}

\foreach \a in {0,...,6}
\foreach \x in {0,1,2}	
{
\pgfmathtruncatemacro{\patchx}{ifthenelse(\x==1,2,ifthenelse(\x==2,1,0))}
\begin{scope}[shift={(3*\patchx,0)}, xscale=-1]

\pgfmathtruncatemacro{\cornerlabel}{mod(7-\a,7)}
\node at (90+51.43*\a:0.35) {\cornerlabel};

\draw[red]
	(B\x 12) -- (B\x 11) -- (A\x 1) -- (A\x 0) -- (A\x 6) -- (B\x 55) -- (B\x 54) -- (B\x 53)
	(B\x 14) -- (B\x 15) -- (A\x 2) -- (A\x 3) -- (A\x 4) -- (A\x 5) -- (B\x 51) -- (B\x 52)
	(B\x 03) -- (B\x 02) -- (B\x 01) -- (B\x 65) -- (B\x 64) -- (B\x 63)
	(B\x 23) -- (B\x 24)
	(B\x 31) -- (B\x 32)
	(B\x 34) -- (B\x 35)
	(B\x 42) -- (B\x 43) -- (B\x 44);

\draw[blue]
	(A\x 1) -- (A\x 2)
	(A\x 5) -- (A\x 6)
	(B\x 00) -- (B\x 01)
	(B\x 44) -- (B\x 45)
	(B\x 52) -- (B\x 53)
	;

\end{scope}
}

\foreach \x/\a in {0/0, 0/4, 0/5, 0/6, 2/3, 1/4, 1/5}
\fill[red!20] 
	(Z\x\a) circle (0.15);

\foreach \x/\a in {0/2, 2/0, 2/5, 1/0, 1/1, 1/3, 1/6}
\fill[blue!20] 
	(Z\x\a) circle (0.15);

\foreach \x/\a in {0/1, 0/3, 2/1, 2/2, 2/4, 2/6, 1/2}
\fill[green!20] 
	(Z\x\a) circle (0.15);
		
\foreach \x/\a/\u/\v in 
	{0/0/4/5, 0/1/2/3, 0/2/2/3, 0/2/5/6, 0/3/3/4, 0/4/1/2, 0/6/2/3,
	2/0/3/4, 2/1/3/4, 2/2/2/3, 2/2/5/6, 2/3/3/4, 2/4/1/2, 2/6/2/3,
	1/0/3/4, 1/1/3/4, 1/2/1/2, 1/2/4/5, 1/3/2/3, 1/3/5/6, 1/6/1/2}
\draw[red, thick]
	(B\x\a\u) -- (B\x\a\v);

\foreach \x/\a/\u/\v in 
	{0/0/3/4, 0/1/3/4, 0/2/1/2, 0/3/2/3, 0/2/4/5, 0/3/5/6, 0/6/1/2,
	2/0/4/5, 2/1/2/3, 2/2/1/2, 2/2/4/5, 2/3/2/3, 2/3/5/6, 2/6/1/2,
	1/0/4/5, 1/1/2/3, 1/2/2/3, 1/2/5/6, 1/3/3/4, 1/4/1/2, 1/6/2/3}
\draw[blue, thick]
	(B\x\a\u) -- (B\x\a\v);

\foreach \x/\a/\u/\v in {
	0/0/0/5, 0/0/1/6, 0/0/2/0, 0/0/3/1, 0/0/4/2, 0/0/5/3, 0/0/6/4,
	0/1/0/5, 0/1/1/4, 0/1/2/3, 0/1/3/2, 0/1/4/1, 0/1/5/0, 0/1/6/6,
	0/2/0/4, 0/2/1/5, 0/2/2/6, 0/2/3/0, 0/2/4/1, 0/2/5/2, 0/2/6/3,
	0/3/0/4, 0/3/1/3, 0/3/2/2, 0/3/3/1, 0/3/4/0, 0/3/5/6, 0/3/6/5,
	0/4/0/1, 0/4/1/2, 0/4/2/3, 0/4/3/4, 0/4/4/5, 0/4/5/6, 0/4/6/0,
	0/5/0/6, 0/5/1/0, 0/5/2/1, 0/5/3/2, 0/5/4/3, 0/5/5/4, 0/5/6/5,
	0/6/0/0, 0/6/1/1, 0/6/2/2, 0/6/3/3, 0/6/4/4, 0/6/5/5, 0/6/6/6,
	2/0/0/1, 2/0/1/2, 2/0/2/3, 2/0/3/4, 2/0/4/5, 2/0/5/6, 2/0/6/0,
	2/1/0/1, 2/1/1/0, 2/1/2/6, 2/1/3/5, 2/1/4/4, 2/1/5/3, 2/1/6/2,
	2/2/0/3, 2/2/1/2, 2/2/2/1, 2/2/3/0, 2/2/4/6, 2/2/5/5, 2/2/6/4,
	2/3/0/3, 2/3/1/4, 2/3/2/5, 2/3/3/6, 2/3/4/0, 2/3/5/1, 2/3/6/2,
	2/4/0/6, 2/4/1/5, 2/4/2/4, 2/4/3/3, 2/4/4/2, 2/4/5/1, 2/4/6/0,
	2/5/0/6, 2/5/1/0, 2/5/2/1, 2/5/3/2, 2/5/4/3, 2/5/5/4, 2/5/6/5,
	2/6/0/0, 2/6/1/6, 2/6/2/5, 2/6/3/4, 2/6/4/3, 2/6/5/2, 2/6/6/1,
	1/0/0/6, 1/0/1/5, 1/0/2/4, 1/0/3/3, 1/0/4/2, 1/0/5/1, 1/0/6/0,
	1/1/0/1, 1/1/1/0, 1/1/2/6, 1/1/3/5, 1/1/4/4, 1/1/5/3, 1/1/6/2,
	1/2/0/3, 1/2/1/4, 1/2/2/5, 1/2/3/6, 1/2/4/0, 1/2/5/1, 1/2/6/2,
	1/3/0/2, 1/3/1/1, 1/3/2/0, 1/3/3/6, 1/3/4/5, 1/3/5/4, 1/3/6/3,
	1/4/0/3, 1/4/1/2, 1/4/2/1, 1/4/3/0, 1/4/4/6, 1/4/5/5, 1/4/6/4,
	1/5/0/1, 1/5/1/0, 1/5/2/6, 1/5/3/5, 1/5/4/4, 1/5/5/3, 1/5/6/2,
	1/6/0/4, 1/6/1/3, 1/6/2/2, 1/6/3/1, 1/6/4/0, 1/6/5/6, 1/6/6/5}
{\pgfmathtruncatemacro{\cornerlabel}{mod(7-\v,7)}
\node at (X\x\a\u) {\cornerlabel};}

\begin{scope}[font=\small]

\node at (0,-1.6) {$120\tilde{2}4\tilde{1}3$};
\node at (3,-1.6) {$\tilde{6}2\tilde{0}5\tilde{3}\tilde{5}0$};
\node at (6,-1.6) {$\tilde{2}\tilde{5}\tilde{3}\tilde{4}5\tilde{5}\tilde{0}$};

\end{scope}

\end{tikzpicture}
\caption{All ${\mc D}_1$-patches for $7^1_{5+6}$-heptagon.}
\label{H111A}
\end{figure}

Proposition \ref{h111tiling} shows that tilings of ${\bb H}^2$ by non-regular $7^1_{5+6}$-heptagon is edge-to-edge and constructed according to the edge and corner system \eqref{h11ev1}. Guided by the system, we find there are exactly three ${\mc D}_1$-patches, in Figure \ref{H111A}.

For each ${\mc D}_1$-patch in Figure \ref{H111A}, we consider 7 neighboring tiles. Each neighboring tile has 3 neighboring tiles in the ${\mc D}_1$-patch. This gives partial information about the ${\mc D}_1$-patch around this neighboring tile, provided the ${\mc D}_1$-patch we start with can be extended to ${\mc D}_2$-patch. We list the information in Table \ref{H111tableA}.

\begin{table}[htp]
\centering
\renewcommand{\arraystretch}{1.2}
\begin{tabular}{|c||c|c|c|c|c|c|c|}
\hline
${\mc D}_1$-patch &
$t_0$ &
$t_1$ &
$t_2$ &
$t_3$ &
$t_4$ &
$t_5$ &
$t_6$  \\
\hline \hline
{\color{red} $120\tilde{2}4\tilde{1}3$} &
{\color{red} $12{*}{*}{*}{*}3$} &
{\color{red} $120{*}{*}{*}{*}$} &
{\color{red} $1{*}{*}{*}{*}\tilde{1}3$} &
{\color{green} ${*}\tilde{5}\tilde{3}\tilde{4}{*}{*}{*}$} &
{\color{blue} ${*}{*}\tilde{0}5\tilde{3}{*}{*}$} &
{\color{green} $\tilde{2}\tilde{5}\tilde{3}{*}{*}{*}{*}$} &
{\color{red} ${*}20\tilde{2}{*}{*}{*}$}  \\
\hline
{\color{blue} $\tilde{6}2\tilde{0}5\tilde{3}\tilde{5}0$} &
{\color{green} $\tilde{2}{*}{*}{*}{*}\tilde{5}\tilde{0}$} &
{\color{blue} $\tilde{6}2\tilde{0}{*}{*}{*}{*}$} &
{\color{green} $\tilde{2}\tilde{5}{*}{*}{*}{*}\tilde{0}$} &
{\color{red} ${*}{*}{*}\tilde{2}4\tilde{1}{*}$} &
{\color{green} ${*}{*}\tilde{3}\tilde{4}5{*}{*}$} &
{\color{green} ${*}{*}{*}{*}5\tilde{5}\tilde{0}$} &
{\color{blue} $\tilde{6}{*}{*}{*}{*}\tilde{5}0$}  \\
\hline
{\color{green} $\tilde{2}\tilde{5}\tilde{3}\tilde{4}5\tilde{5}\tilde{0}$} &
{\color{blue} ${*}2\tilde{0}5{*}{*}{*}$} &
{\color{red} ${*}{*}{*}{*}4\tilde{1}3$} &
{\color{red} ${*}{*}0\tilde{2}4{*}{*}$} &
{\color{blue} ${*}{*}{*}5\tilde{3}\tilde{5}{*}$} &
{\color{green} ${*}{*}{*}\tilde{4}5\tilde{5}{*}$} &
{\color{blue} ${*}{*}{*}{*}\tilde{3}\tilde{5}0$} &
{\color{blue} $\tilde{6}2{*}{*}{*}{*}0$}  \\
\hline
\end{tabular}
\caption{${\mc D}_1$-patches for neighboring tiles of a $7^1_{5+6}$-heptagon.}
\label{H111tableA}
\end{table}

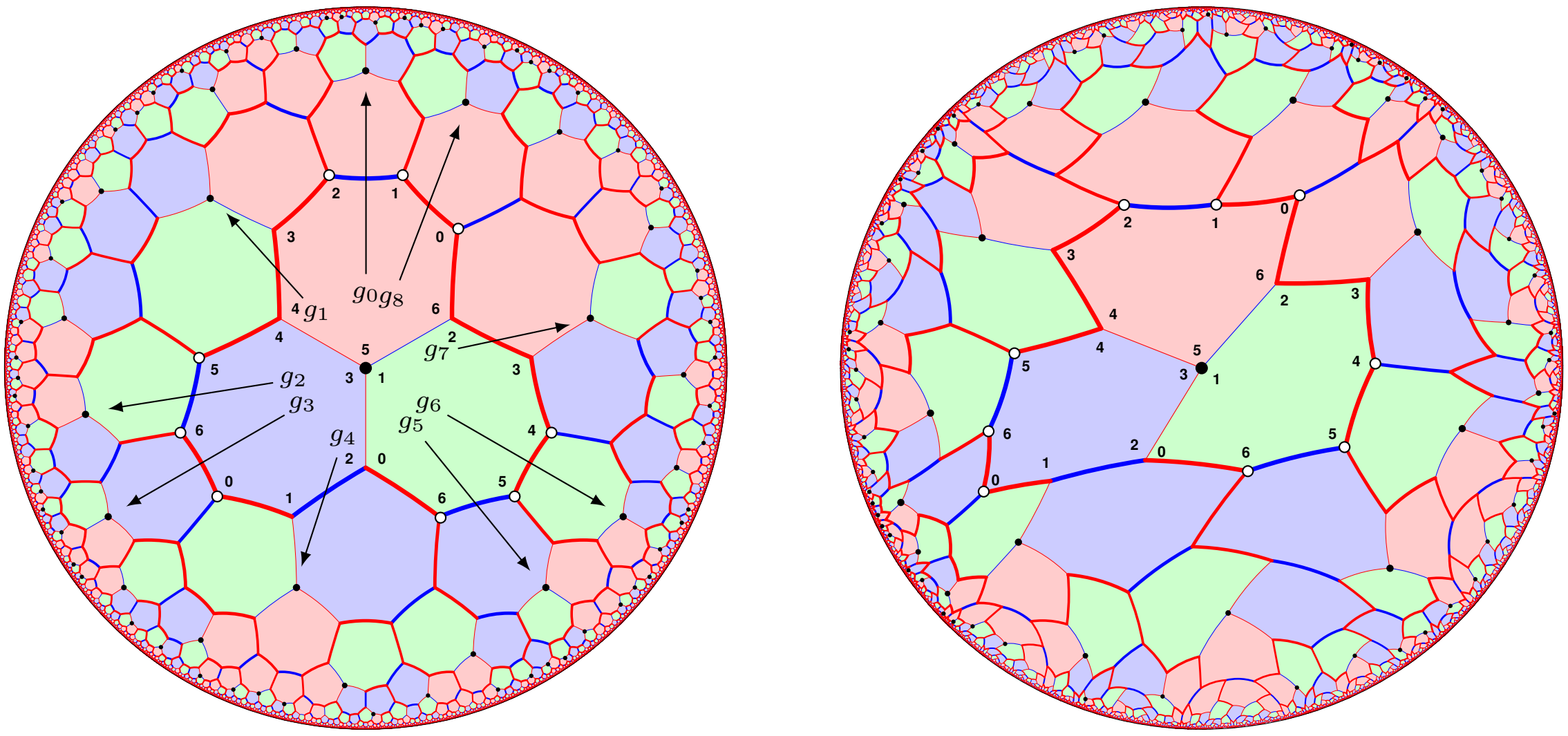
\begin{figure}[h!]
\centering
\begin{tikzpicture}[>=latex,scale=1]

\pgftext{
	\includegraphics[scale=0.4]{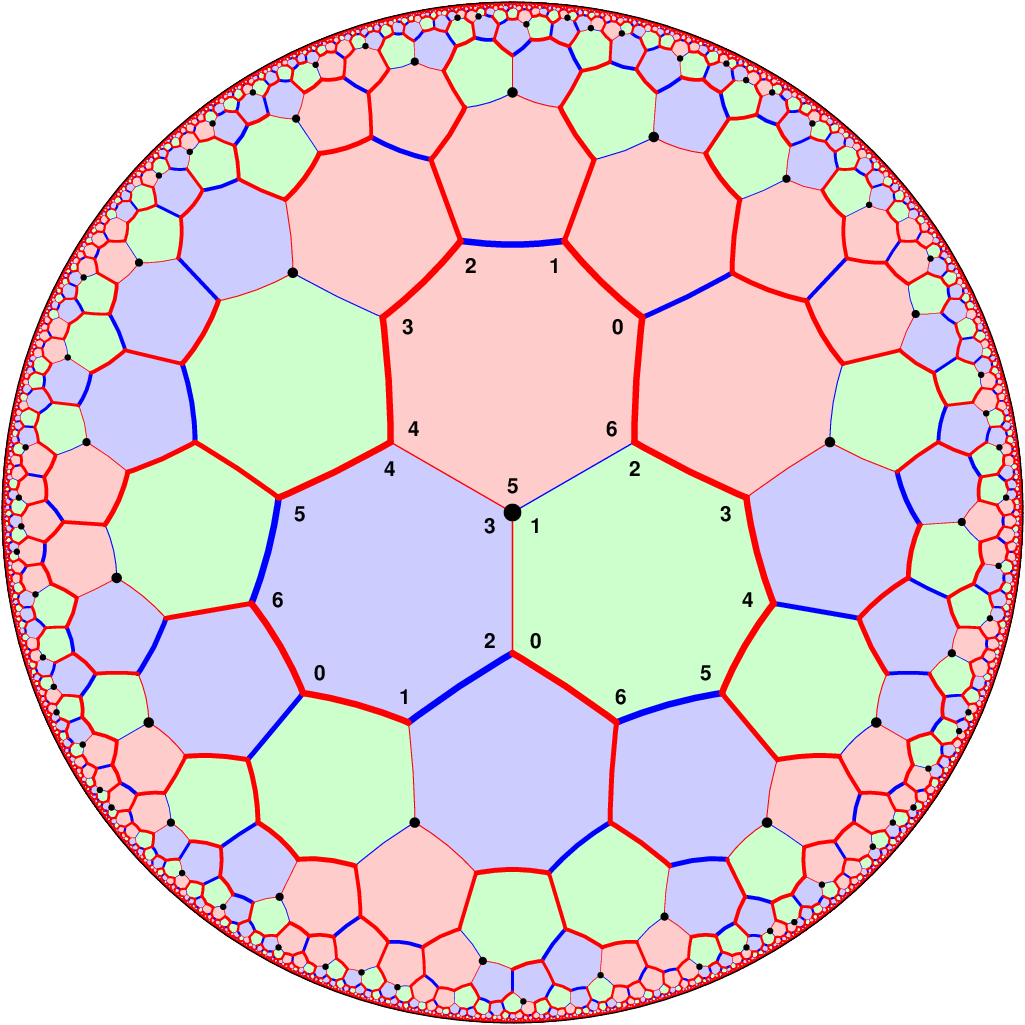}
	};
\begin{scope}[xscale=-1]

\foreach \a/\b in {90/2.65, 111/2.5, 168/2}
\foreach \u in {0,1,2}	
\draw[rotate=120*\u, ->]
	(\a:0.9) -- (\a:\b);

\foreach \u in {0,1,2}	
\filldraw[fill=white, rotate=120*\u]
	(79.2+120*\u:1.88) circle (0.05)
	(100.8+120*\u:1.88) circle (0.05)
	(3.5+120*\u:1.6) circle (0.05);

\foreach \a/\b in {90/0, 111/8, 168/7} 
\node at (\a:0.7) {\scriptsize $g_{\b}$}; 

\foreach \a/\b in {90/6, 111/5, 168/4} 
\node at (120+\a:0.7) {\scriptsize $g_{\b}$};

\foreach \a/\b in {90/3, 111/2, 168/1} 
\node at (-120+\a:0.7) {\scriptsize $g_{\b}$};

\end{scope}

\begin{scope}[xshift=8cm]

\pgftext{
	\includegraphics[scale=0.4]{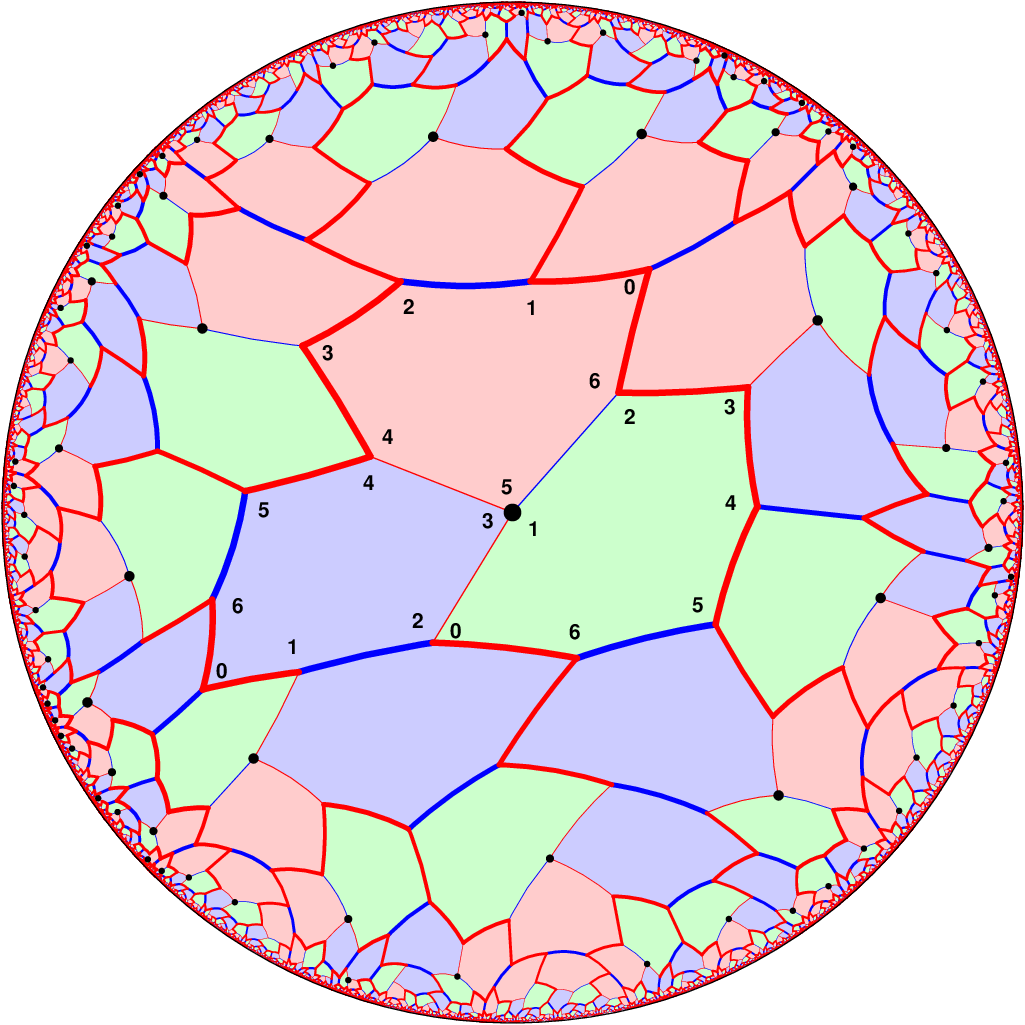}
	};

\foreach \a/\b in 
	{175.5/1.8, 115.5/1.73, 85/1.57, 60.5/1.9, 1.5/1.66, -29/1.56, -66/1.08, 209.5/2.4, 196.5/2.13}	
\filldraw[fill=white]
	(\a:\b) circle (0.05);
		
\end{scope}

\end{tikzpicture}
\caption{Schematic and real drawings of the unique tiling ${\mc T}(7^1_{5+6})$ of the Poincar\'e disk ${\bb H}^2$ by the non-regular $7^1_{5+6}$-heptagon.}
\label{H111B}
\end{figure}

Compared with the ${\mc D}_1$-patches in Figure \ref{H111A}, we may easily fill all $*$ and find the corresponding ${\mc D}_1$-patches, which are indicated by colored dots in Figure \ref{H111A}. Then we use the data to construct ${\mc D}_1$-patches around all the neighboring tiles in Figure \ref{H111A}. We find that all these ${\mc D}_1$-patches are compatible, and form unique extensions of the ${\mc D}_1$-patches in Figure \ref{H111A} to ${\mc D}_2$-patches. This implies that we may repeatedly use the three ${\mc D}_1$-patches to further construct the tilings, and get a unique tiling ${\mc T}(7^1_{5+6})$ of the Poincar\'e disk ${\bb H}^2$ in Figure \ref{H111B}. 

The tiling ${\mc T}(7^1_{5+6})$ is actually the tiling ${\mc T}(T)$ by the triple $T$ of red, blue and green tiles at the center of Figure \ref{H111B}. All the triple tiles are indicated by three tiles around $\bullet$ vertices. The symmetry $G(7^1_{5+6})$ of the tiling ${\mc T}(7^1_{5+6})$ is transitive on the triples, and the map
\[
g\in G(7^1_{5+6})\mapsto gT\in {\mc T}(T)
\]
is a one-to-one correspondence similar to \eqref{oneone}. The triple $T$ at the center is surrounded by nine triples $g_0T,g_1T,\dots,g_8T$ similar to $g_0t,g_1t,\dots,g_6t$ in Figure \ref{H5B}. We number these neighboring triples counterclockwise, starting with $g_0T$ at the top of the schematic drawing. Then we similarly get the presentation of the symmetry group
\begin{align*}
G(7^1_{5+6})
&=\langle g_0,g_1,g_2,g_3,g_4,g_5,g_6,g_7,g_8
\; | \\
&\qquad 
g_0^2=g_3^2=g_4^2=g_6^2=g_8^2
=g_1g_7=g_2g_5
=g_0g_1g_8=g_2g_6g_7=g_3g_4g_5=1\rangle \\
&=\langle g_0,g_4,g_6,g_8 \; | \;
g_0^2=g_4^2=g_6^2=g_8^2=(g_{4}g_{6}g_{8}g_{0})^2=1\rangle,
\end{align*}
with
\[
g_1=g_0g_8,\;
g_2=g_0g_8g_6,\;
g_3=g_4g_6g_8g_0,\;
g_5=g_6g_8g_0,\;
g_7=g_8g_0.
\]

Since $G(7^1_{5+6})$ has three orbits, tilings of a closed surface of Euler characteristic $\chi$ by the non-regular $7^1_{5+6}$-heptagon are the quotients of ${\bb H}^2$ by conjugate classes of torsionless subgroups of $G(7^1_{5+6})$ of finite index $\frac{f}{3}=-2\chi$. Since all triple tiles in ${\mc T}(T)$ have the same orientation, the surface must be orientable. We use standard computer package GAP to find the numbers of tilings of $g{\bb T}^2$, $g=2,3,4$, in Table \ref{table5}.

\begin{figure}[htp]
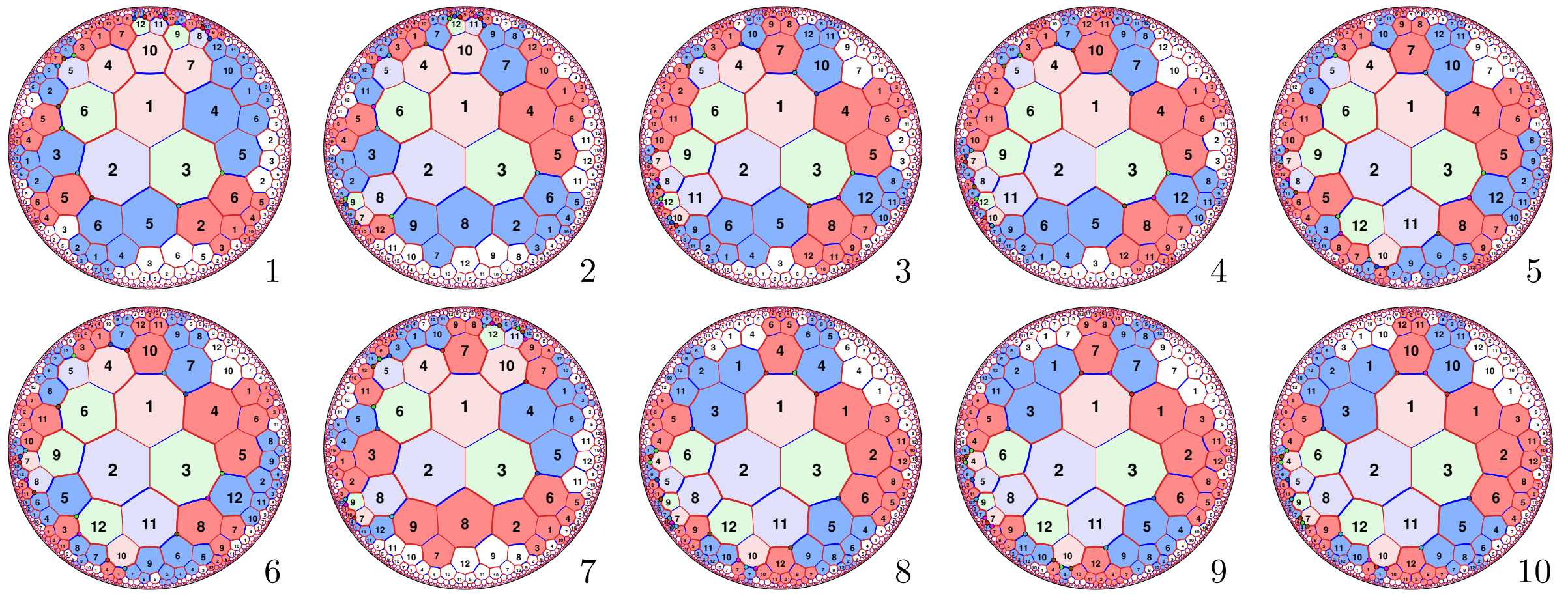

\centering
\begin{tikzpicture}[>=latex,scale=1]

%% 1-5

\foreach \a in {1,...,5}
{
\begin{scope}[shift={(-3.2*1+3.2*\a,0)}]

\pgftext{
	\includegraphics[scale=0.166]{H111domain/H111T\a}
	};

\node at (1.25,-1.25) {\a};

\end{scope}
}

%% 6-10

\foreach \a in {6,...,10}
{
\begin{scope}[shift={(-3.2*6+3.2*\a,-3.05)}]

\pgftext{
	\includegraphics[scale=0.166]{H111domain/H111T\a}
	};

\node at (1.25,-1.25) {\a};

\end{scope}
}

\end{tikzpicture}
\caption{Fundamental domains of ${\mc T}_k(7^1_{5+6},2{\bb T}^2)$, $k=1,2,\dots,10$, of tilings of $2{\bb T}^2$ by generic $7^1_{5+6}$-heptagon.}
\label{H111C}
\end{figure}

There are 10 tilings ${\mc T}_k(7^1_{5+6},2{\bb T}^2)$ of $2{\bb T}^2$ by the non-regular $7^1_{5+6}$-heptagon. Figure \ref{H111C} presents the lightly colored fundamental domains and the red/blue neighboring domains. We choose the more compact versions of the fundamental domains that are unions of four translations of the triple tile. 

\begin{figure}[htp]
\centering
\begin{tikzpicture}[>=latex,scale=1]

\pgftext{
	\includegraphics[scale=0.4]{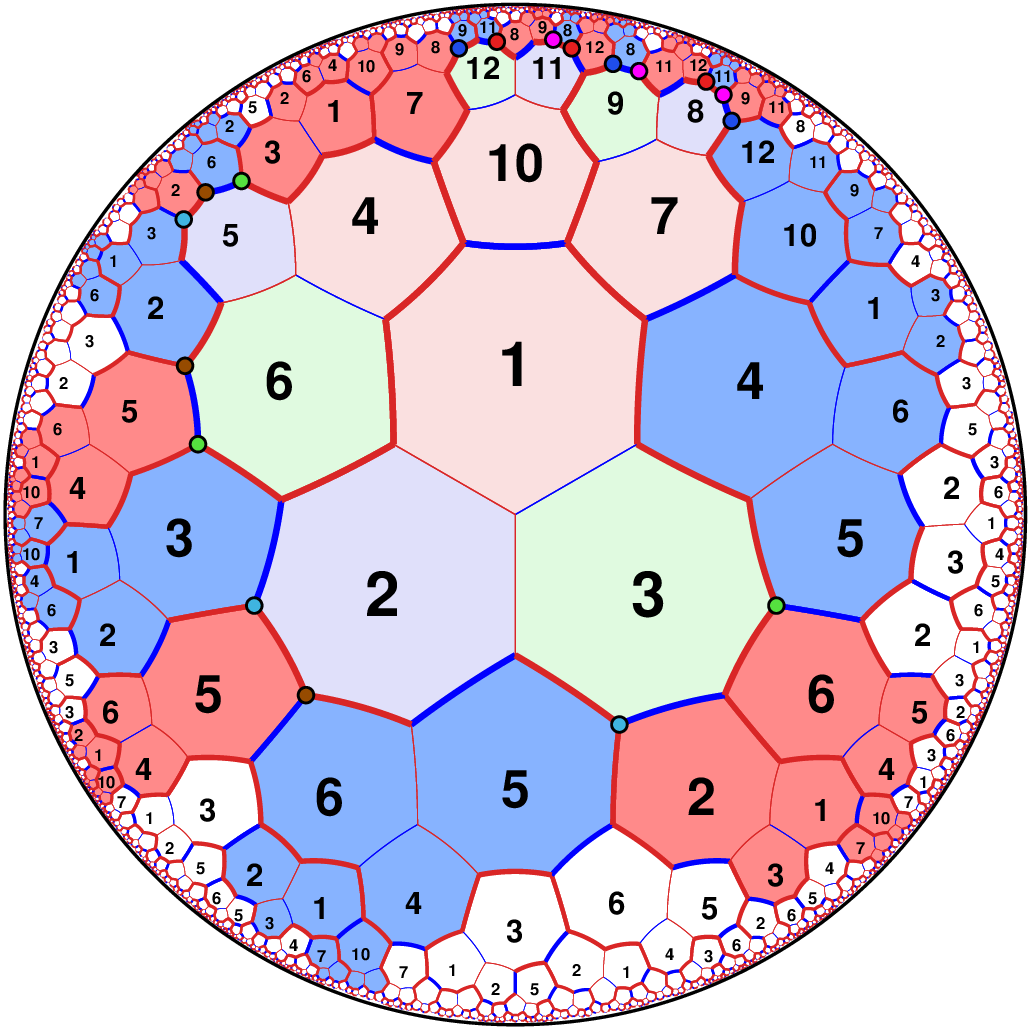}
	};

\foreach \a/\b in {
	59/1, 64/2, 70/3, 76/4, 80/5, 84/6, 90/7, 96/8, 112/9, 130/10, 138/11, 150/12, 170/13, 193/14, 210/15, 255/16, 320/17, 20/18}
\node at (\a:3.7) {\small \b};

\end{tikzpicture}
\caption{Fundamental domain of the first tiling (in Figure \ref{H111C}) ${\mc T}_1(7^1_{5+6},2{\bb T}^2)$ of $2{\bb T}^2$ by the non-regular $7^1_{5+6}$-heptagon, and neighboring translations of the fundamental domain.}
\label{H111E}
\end{figure}

Figure \ref{H111E} is the detail of the first tiling ${\mc T}_1(7^1_{5+6},2{\bb T}^2)$ of $2{\bb T}^2$ in Figure \ref{H111C}. The lightly colored fundamental domain $d$ consists of four triple tiles $T,g_4T,g_7T,g_8T$, and is surrounded by 18 neighboring domains $a_id$, numbered counterclockwise from $a_1d$ at the top, where
\begin{align*}
a_{1} &=g_8g_2g_8, &
a_{4} &=g_8g_5g_8, &
a_{7} &=g_0g_4g_8, &
a_{10} &=g_1g_2g_7, &
a_{13} &=g_1g_5g_7, &
a_{16} &=g_4g_7, \\
a_{2} &=g_8g_3g_0, &
a_{5} &=g_0g_2g_0, &
a_{8} &=g_0g_5g_0, &
a_{11} &=g_1g_3, &
a_{14} &=g_2, &
a_{17} &=g_5, \\
a_{3} &=g_8g_4g_0, &
a_{6} &=g_0g_3g_8, &
a_{9} &=g_1^2, &
a_{12} &=g_1g_4, &
a_{15} &=g_3g_7, &
a_{18} &=g_7^2.
\end{align*}
The torsionless subgroup $G_1(7^1_{5+6},2{\bb T}^2)$ is generated by the 18 isometries $a_i$, subject to the relations 
\[
a_1a_4=
a_2a_6=
a_3a_7=
a_5a_8=
a_9a_{18}=
a_{10}a_{13}=
a_{11}a_{15}=
a_{12}a_{16}=
a_{14}a_{17}=1,
\]
and 
\[
{\color{red} a_3a_8a_6}
={\color{blue} a_1a_5a_9}
={\color{green} a_{10}a_{14}a_{18}}
={\color{brown} a_{11}a_{16}a_{13}}
={\color{cyan} a_{12}a_{17}a_{15}}
={\color{magenta} a_2a_7a_4}
=1.
\]

\subsection{Tilings by $7^2_{5+6}$-Heptagon}
\label{H112}

By the explanation in Section \ref{calculation}, the 
non-regular $7^2_{5+6}$-heptagon has the following specific values
\begin{align}
7^2_{5+6} &\colon  
|\bar{0}|=|\bar{1}|=|\bar{2}|=|\bar{4}|=|\bar{5}|=|\bar{6}|=\text{arccosh}\tfrac{3t_0^3-15t_0^2+t_0+3}{3(1+t_0)^2(1-3t_0)}
\approx 0.62245, \nonumber \\
&\quad
|\bar{3}|=2\text{arccosh}\tfrac{2}{\sqrt{3}(1+t_0)}
\approx 0.51965;  \nonumber \\
&\quad
[2]=[4]=[6]=\tfrac{2\pi}{3}, \label{112heptagon} \\
&\quad
[0]=\tfrac{2\pi}{3}-2\theta_0\approx 0.41254\pi,\;
[1]=\tfrac{2\pi}{3}+2\theta_0\approx 0.92080\pi, \nonumber \\
&\quad
[3]=\tfrac{2\pi}{3}-\theta_0\approx 0.53960\pi,\;
[5]=\tfrac{2\pi}{3}+\theta_0\approx 0.79373\pi.  \nonumber 
\end{align}
Here $t_0\approx 0.11679$ is the unique root of \eqref{712solution} in the interval $(-\frac{1}{3},\frac{1}{3})$, and $\theta_0=2\text{arctan}\sqrt{3}t_0$. 

To better describe tiling by the heptagon, we denote the corners $2,4,6$ by $\alpha$. Then the complete vertex set is 
\begin{equation}\label{h112v}
\bar{V}(7^2_{5+6})=\{\alpha^3,\alpha 01, \alpha 35, 055,133\},\quad
\alpha=2,4,6. 
\end{equation}
We find 17 ${\mc D}_1$-patches, among which 13 pass $\partial {\mc D}_1$-verification. These 13 ${\mc D}_1$-patches are given by the first column of Table \ref{H112tableA}. Note that each ${\mc D}_1$-patch occupies two rows, with the upper row being the notation for (partial) ${\mc D}_1$-patches, and the lower row being the possible labels of the (partial) ${\mc D}_1$-patches. Figure \ref{H112A} shows the first three of these ${\mc D}_1$-patches.

\begin{figure}[htp]
\centering
\begin{tikzpicture}[>=latex, scale=1, thick, font=\tiny]

\foreach \a in {0,...,6}
{
\fill[red!20, rotate=51.43*\a]
	(0,0) -- (90:0.5) -- (90+51.43:0.5);
\fill[blue!20, xshift=3cm, rotate=51.43*\a]
	(0,0) -- (90:0.5) -- (90+51.43:0.5);
\fill[green!20, xshift=6cm, rotate=51.43*\a]
	(0,0) -- (90:0.5) -- (90+51.43:0.5);
}

\foreach \a in {0,...,6}
\foreach \x in {0,1,2}
\draw[red, xshift=3*\x cm]
	(38.57+51.43*\a:0.5) -- (90+51.43*\a:0.5) -- (90+51.43*\a:0.9) --  (102+51.43*\a:1.2) -- (115.714+51.43*\a:1.3) -- (129.43+51.43*\a:1.2) -- (141.43+51.43*\a:0.9);

\foreach \x in {0,1,2}	
{
\begin{scope}[xshift=3*\x cm]

\foreach \a in {0,1,3,5}
\node at (90+51.43*\a:0.35) {\a};

\foreach \a in {2,4,6}
\node at (90+51.43*\a:0.35) {$\alpha$};

\draw[blue]
	(90+51.43*3:0.5) -- (90+51.43*4:0.5);

\end{scope}
}

\foreach \a in {0,...,6}
\foreach \b in {0,...,6}
\foreach \x in {0,1,2}
{
\begin{scope}[shift={(3*\x, 0)}]

\coordinate (A\x\a) at (90+51.43*\a:0.5);
\coordinate (B\x\a 0) at (90+51.43*\a:0.5);
\coordinate (B\x\a 1) at (90+51.43*\a:0.9);
\coordinate (B\x\a 2) at (102+51.43*\a:1.2);
\coordinate (B\x\a 3) at (115.714+51.43*\a:1.3);
\coordinate (B\x\a 4) at (129.43+51.43*\a:1.2);
\coordinate (B\x\a 5) at (141.43+51.43*\a:0.9);
\coordinate (B\x\a 6) at (141.43+51.43*\a:0.5);

\coordinate (Z\x\a) at (115.7+51.43*\a:0.85);
\coordinate (X\x\a 0) at (100+51.43*\a:0.6);
\coordinate (X\x\a 1) at (99+51.43*\a:0.87);
\coordinate (X\x\a 2) at (105+51.43*\a:1.1);
\coordinate (X\x\a 3) at (115.714+51.43*\a:1.17);
\coordinate (X\x\a 4) at (126.43+51.43*\a:1.1);
\coordinate (X\x\a 5) at (132.43+51.43*\a:0.87);
\coordinate (X\x\a 6) at (131.43+51.43*\a:0.6);

\end{scope}
}
	
\foreach \x/\a in {0/2, 0/3, 0/4, 0/5, 1/1, 1/2, 2/2}
\fill[red!20] 
	(Z\x\a) circle (0.15);

\foreach \x/\a in {0/1, 0/6, 1/6, 2/0, 2/1, 2/3, 2/5}
\fill[blue!20] 
	(Z\x\a) circle (0.15);

\foreach \x/\a in {0/0, 1/0, 1/3, 1/4, 1/5, 2/4, 2/6}
\fill[green!20] 
	(Z\x\a) circle (0.15);
		
\foreach \x/\a/\u/\v in 
	{0/0/5/6, 0/2/4/5, 0/5/0/1, 0/6/2/3,
	1/0/4/5, 1/1/3/4, 1/2/1/2, 1/4/4/5, 1/5/3/4, 1/6/3/4,
	2/0/1/2, 2/1/3/4, 2/2/3/4, 2/4/2/3, 2/5/5/6}
\draw[blue]
	(B\x\a\u) -- (B\x\a\v);

\foreach \x/\a/\u/\v in {
	0/0/0/$\alpha$, 0/0/1/1, 0/0/2/0, 0/0/3/$\alpha$, 0/0/4/5, 0/0/5/$\alpha$, 0/0/6/3,
	0/1/0/3, 0/1/1/$\alpha$, 0/1/2/5, 0/1/3/$\alpha$, 0/1/4/0, 0/1/5/1, 0/1/6/$\alpha$,
	0/2/0/$\alpha$, 0/2/1/0, 0/2/2/1, 0/2/3/$\alpha$, 0/2/4/3, 0/2/5/$\alpha$, 0/2/6/5,
	0/3/0/$\alpha$, 0/3/1/5, 0/3/2/$\alpha$, 0/3/3/0, 0/3/4/1, 0/3/5/$\alpha$, 0/3/6/3,
	0/4/0/5, 0/4/1/$\alpha$, 0/4/2/0, 0/4/3/1, 0/4/4/$\alpha$, 0/4/5/3, 0/4/6/$\alpha$,
	0/5/0/3, 0/5/1/$\alpha$, 0/5/2/5, 0/5/3/$\alpha$, 0/5/4/0, 0/5/5/1, 0/5/6/$\alpha$,
	0/6/0/$\alpha$, 0/6/1/3, 0/6/2/$\alpha$, 0/6/3/5, 0/6/4/$\alpha$, 0/6/5/0, 0/6/6/1,
	1/0/0/1, 1/0/1/0, 1/0/2/$\alpha$, 1/0/3/5, 1/0/4/$\alpha$, 1/0/5/3, 1/0/6/$\alpha$,
	1/1/0/0, 1/1/1/1, 1/1/2/$\alpha$, 1/1/3/3, 1/1/4/$\alpha$, 1/1/5/5, 1/1/6/$\alpha$,
	1/2/0/$\alpha$, 1/2/1/3, 1/2/2/$\alpha$, 1/2/3/5, 1/2/4/$\alpha$, 1/2/5/0, 1/2/6/1,
	1/3/0/3, 1/3/1/$\alpha$, 1/3/2/1, 1/3/3/0, 1/3/4/$\alpha$, 1/3/5/5, 1/3/6/$\alpha$,
	1/4/0/$\alpha$, 1/4/1/0, 1/4/2/1, 1/4/3/$\alpha$, 1/4/4/3, 1/4/5/$\alpha$, 1/4/6/5,
	1/5/0/0, 1/5/1/$\alpha$, 1/5/2/5, 1/5/3/$\alpha$, 1/5/4/3, 1/5/5/$\alpha$, 1/5/6/1,
	1/6/0/0, 1/6/1/1, 1/6/2/$\alpha$, 1/6/3/3, 1/6/4/$\alpha$, 1/6/5/5, 1/6/6/$\alpha$,
	2/0/0/5, 2/0/1/$\alpha$, 2/0/2/3, 2/0/3/$\alpha$, 2/0/4/1, 2/0/5/0, 2/0/6/$\alpha$,
	2/1/0/0, 2/1/1/$\alpha$, 2/1/2/5, 2/1/3/$\alpha$, 2/1/4/3, 2/1/5/$\alpha$, 2/1/6/1,
	2/2/0/0, 2/2/1/$\alpha$, 2/2/2/5, 2/2/3/$\alpha$, 2/2/4/3, 2/2/5/$\alpha$, 2/2/6/1,
	2/3/0/3, 2/3/1/$\alpha$, 2/3/2/1, 2/3/3/0, 2/3/4/$\alpha$, 2/3/5/5, 2/3/6/$\alpha$,
	2/4/0/$\alpha$, 2/4/1/5, 2/4/2/$\alpha$, 2/4/3/3, 2/4/4/$\alpha$, 2/4/5/1, 2/4/6/0,
	2/5/0/5, 2/5/1/$\alpha$, 2/5/2/0, 2/5/3/1, 2/5/4/$\alpha$, 2/5/5/3, 2/5/6/$\alpha$,
	2/6/0/$\alpha$, 2/6/1/3, 2/6/2/$\alpha$, 2/6/3/1, 2/6/4/0, 2/6/5/$\alpha$, 2/6/6/5}
\node at (X\x\a\u) {\v};

\begin{scope}[font=\small]

\node at (0,-1.6) {{\color{red} 1}: $\tilde{\alpha} 3\alpha\alpha 53\alpha$};
\node at (3,-1.6) {{\color{blue} 2}: $\tilde{1}0\alpha \tilde{3}\alpha\tilde{0}0$};
\node at (6,-1.6) {{\color{green} 3}: $\tilde{5}\tilde{0}\tilde{0}\tilde{3}\tilde{\alpha}5\tilde{\alpha}$};

\end{scope}

\end{tikzpicture}
\caption{The first three of 13 ${\mc D}_1$-patches for non-regular $7^2_{5+6}$-heptagon in Table \ref{H112tableA}.}
\label{H112A}
\end{figure}

Similar to the discussion for non-regular $7^1_{5+6}$-heptagon, for each of 13 ${\mc D}_1$-patches, we consider 7 neighboring tiles. Each neighboring tile has 3 neighboring tiles in the ${\mc D}_1$-patch, which gives partial information about the ${\mc D}_1$-patch around this neighboring tile. We gather all these information in the upper rows of Table \ref{H112tableA} similar to Table \ref{H111tableA}. We also indicate which of the 13 ${\mc D}_1$-patches are compatible with the partial information, in the lower rows.

\begin{table}[htp]
\centering
\renewcommand{\arraystretch}{1.1}
\begin{tabular}{|c||c|c|c|c|c|c|c|} 
\hline
${\mc D}_1$-patch & 
$t_0$ & 
$t_1$ & 
$t_2$ & 
$t_3$ & 
$t_4$ & 
$t_5$ & 
$t_6$ \\
\hline
\hline
\textcolor{red}{$\tilde{\alpha}3\alpha\alpha53\alpha$} & 
${*}\tilde{0}\tilde{0}\tilde{3}{*}{*}{*}$ &
${*}0\alpha\tilde{3}{*}{*}{*}$ & 
${*}{*}{*}{*}53\alpha$ &  
${*}{*}\alpha\alpha5{*}{*}$ & 
${*}{*}{*}\alpha53{*}$ &  
${*}3\alpha\alpha{*}{*}{*}$ & 
$\tilde{1}0\alpha{*}{*}{*}{*}$  \\
\arrayrulecolor{gray!50}\hline\arrayrulecolor{black}
\textcolor{red}{1} & 
\textcolor{green}{3},10,11 & 
\textcolor{blue}{2},5,13 & 
\textcolor{red}{1},12 & 
\textcolor{red}{1},12 & 
\textcolor{red}{1},12 & 
\textcolor{red}{1} & 
\textcolor{blue}{2} \\
\hline
\textcolor{blue}{$\tilde{1}0\alpha\tilde{3}\alpha\tilde{0}0$} & 
$\tilde{5}\tilde{0}\tilde{0}{*}{*}{*}{*}$ & 
$\tilde{\alpha}{*}{*}{*}{*}3\alpha$ & 
$\tilde{\alpha}3\alpha{*}{*}{*}{*}$ & \
${*}{*}\tilde{0}\tilde{3}\tilde{\alpha}{*}{*}$ &
${*}{*}{*}{*}\tilde{\alpha}5\tilde{\alpha}$ &
$\tilde{5}\tilde{0}{*}{*}{*}{*}\tilde{\alpha}$ &  
$\tilde{1}{*}{*}{*}{*}\tilde{0}0$ \\
\arrayrulecolor{gray!50}\hline\arrayrulecolor{black}
\textcolor{blue}{2} & 
\textcolor{green}{3},11 & 
\textcolor{red}{1} & 
\textcolor{red}{1},4,6 & 
\textcolor{green}{3} & 
\textcolor{green}{3},5 & 
\textcolor{green}{3} & 
\textcolor{blue}{2},10 \\
\hline
\textcolor{green}{$\tilde{5}\tilde{0}\tilde{0}\tilde{3}\tilde{\alpha}5\tilde{\alpha}$} & 
${*}{*}{*}{*}\alpha\tilde{0}0$ &
$\tilde{1}0{*}{*}{*}{*}0$ &
$\tilde{\alpha}3{*}{*}{*}{*}\alpha$ &  
${*}{*}\alpha\tilde{3}\alpha{*}{*}$ & 
$\tilde{5}{*}{*}{*}{*}5\tilde{\alpha}$ & 
${*}{*}{*}\tilde{3}\alpha\tilde{0}{*}$ &  
${*}{*}{*}\tilde{3}\tilde{\alpha}5{*}$ \\
\arrayrulecolor{gray!50}\hline\arrayrulecolor{black}
\textcolor{green}{3} & 
\textcolor{blue}{2},7,10 & 
\textcolor{blue}{2} & 
\textcolor{red}{1},6 & 
\textcolor{blue}{2},4,6,13 & 
\textcolor{green}{3},5 & 
\textcolor{blue}{2},4,7,9,10 & 
\textcolor{green}{3},5 \\
\hline
$\tilde{\alpha}3\alpha\tilde{3}\alpha\tilde{0}\tilde{0}$ & 
${*}\alpha\tilde{0}\tilde{3}{*}{*}{*}$ & 
${*}3\alpha\tilde{3}{*}{*}{*}$ &
$\tilde{\alpha}3\alpha{*}{*}{*}{*}$ & 
${*}{*}\tilde{0}\tilde{3}\tilde{\alpha}{*}{*}$ & 
${*}{*}{*}{*}\tilde{\alpha}5\tilde{\alpha}$ & 
$\tilde{5}0{*}{*}{*}{*}\tilde{\alpha}$ & 
$\tilde{\alpha}\alpha{*}{*}{*}{*}\alpha$ \\
\arrayrulecolor{gray!50}\hline\arrayrulecolor{black}
4 & 
7,8,9 & 
4,6 & 
1,4,6 & 
3 & 
3,5 & 
5 & 
8 \\
\hline
$\tilde{5}0\alpha\tilde{3}\tilde{\alpha}5\tilde{\alpha}$ & 
${*}{*}{*}{*}\alpha\tilde{0}\tilde{0}$ & 
$\tilde{\alpha}{*}{*}{*}{*}3\alpha$ &  
$\tilde{\alpha}3\alpha{*}{*}{*}{*}$ & 
${*}{*}\tilde{0}\tilde{3}\alpha{*}{*}$ & 
$\tilde{5}{*}{*}{*}{*}5\tilde{\alpha}$ &  
${*}{*}{*}\tilde{3}\alpha\tilde{0}{*}$ & 
${*}{*}{*}\tilde{3}\tilde{\alpha}5{*}$  \\
\arrayrulecolor{gray!50}\hline\arrayrulecolor{black}
5 & 
4,9 & 
1 & 
1,4,6 & 
7,8,9,10,11 & 
3,5 & 
2,4,7,9,10 & 
3,5 \\
\hline
$\tilde{\alpha}3\alpha\tilde{3}\alpha0\alpha$ & 
${*}\tilde{0}\tilde{0}\tilde{3}{*}{*}{*}$ &
${*}3\alpha\tilde{3}{*}{*}{*}$ & 
$\tilde{\alpha}3\alpha{*}{*}{*}{*}$ & 
${*}{*}\tilde{0}\tilde{3}\tilde{\alpha}{*}{*}$ &
${*}{*}{*}{*}15\tilde{\alpha}$ &
$5{*}{*}{*}{*}3\alpha$ &
$\tilde{1}0\alpha{*}{*}{*}{*}$ \\
\arrayrulecolor{gray!50}\hline\arrayrulecolor{black}
6 & 
3,10,11 & 
4,6 & 
1,4,6 & 
3 & 
& 
& 
2 \\
\hline
$1\alpha\tilde{0}\tilde{3}\alpha\tilde{0}0$ & 
$1\alpha{*}{*}{*}{*}3$ &  
$\tilde{\alpha}\alpha0{*}{*}{*}{*}$ &  
$\tilde{\alpha}3{*}{*}{*}{*}\tilde{0}$ &
${*}{*}\alpha\tilde{3}\tilde{\alpha}{*}{*}$ & 
${*}{*}{*}{*}\tilde{\alpha}5\tilde{\alpha}$ &
$\tilde{5}\tilde{0}{*}{*}{*}{*}\tilde{\alpha}$ & 
$\tilde{1}{*}{*}{*}{*}\alpha0$  \\
\arrayrulecolor{gray!50}\hline\arrayrulecolor{black}
7 & 
& 
& 
4 & 
5 & 
3,5 & 
&  \\
\hline
$\tilde{\alpha}\alpha\tilde{0}\tilde{3}\alpha0\alpha$ & 
$\tilde{\alpha}{*}{*}{*}{*}\tilde{0}\tilde{0}$ & 
$\tilde{\alpha}\alpha\tilde{0}{*}{*}{*}{*}$ & 
$\tilde{\alpha}3{*}{*}{*}{*}\tilde{0}$ & 4 
${*}{*}\alpha\tilde{3}\tilde{\alpha}{*}{*}$ & 
${*}{*}{*}{*}15\tilde{\alpha}$ & 
$5{*}{*}{*}{*}3\alpha$ &
$\tilde{5}0\alpha{*}{*}{*}{*}$ \\
\arrayrulecolor{gray!50}\hline\arrayrulecolor{black}
8 &
4,9 & 
8,9 & 
4 & 
5 & 
& 
& 
5,12,13 \\
\hline
$\tilde{\alpha}\alpha\tilde{0}\tilde{3}\alpha\tilde{0}\tilde{0}$ & 
$\tilde{\alpha}{*}{*}{*}{*}\alpha\tilde{0}$ &  
$\tilde{\alpha}\alpha\tilde{0}{*}{*}{*}{*}$ & 
$\tilde{\alpha}3{*}{*}{*}{*}\tilde{0}$ &  
${*}{*}\alpha\tilde{3}\tilde{\alpha}{*}{*}$ & 
${*}{*}{*}{*}\tilde{\alpha}5\tilde{\alpha}$ &
$\tilde{5}0{*}{*}{*}{*}\tilde{\alpha}$ & 
$\tilde{\alpha}\alpha{*}{*}{*}{*}\alpha$  \\
\arrayrulecolor{gray!50}\hline\arrayrulecolor{black}
9 & 
& 
8,9 & 
4 & 
5 & 
3,5 & 
5 & 
8 \\
\hline
$\tilde{1}\tilde{0}\tilde{0}\tilde{3}\alpha\tilde{0}0$ & 
$\tilde{5}\tilde{0}0{*}{*}{*}{*}$ & 
$\tilde{1}0{*}{*}{*}{*}3$ &  
$\tilde{\alpha}3{*}{*}{*}{*}\alpha$ & 1
${*}{*}\alpha\tilde{3}\tilde{\alpha}{*}{*}$ &  
${*}{*}{*}{*}\tilde{\alpha}5\tilde{\alpha}$ & 
$\tilde{5}\tilde{0}{*}{*}{*}{*}\tilde{\alpha}$ & 
$\tilde{1}{*}{*}{*}{*}\tilde{0}0$  \\
\arrayrulecolor{gray!50}\hline\arrayrulecolor{black}
10 & 
& 
& 
1,6 & 
5 & 
3,5 & 
3 & 
2,10 \\
\hline
$\tilde{5}\tilde{0}\tilde{0}\tilde{3}\alpha0\alpha$ & 
${*}{*}{*}{*}\tilde{\alpha}\tilde{0}0$ &  
$\tilde{1}0{*}{*}{*}{*}0$ & 
$\tilde{\alpha}3{*}{*}{*}{*}\alpha$ &
${*}{*}\alpha\tilde{3}\tilde{\alpha}{*}{*}$ & 
${*}{*}{*}{*}15\tilde{\alpha}$ & 
$5{*}{*}{*}{*}0\alpha$ & 
${*}{*}{*}{*}\tilde{\alpha}0\alpha$  \\
\arrayrulecolor{gray!50}\hline\arrayrulecolor{black}
11 & 
& 
2 & 
1,6 & 
5 & 
& 
&   \\
\hline
$\tilde{5}0\alpha\alpha53\alpha$ & 
${*}{*}{*}{*}\tilde{\alpha}\tilde{0}\tilde{0}$ &  
$\tilde{\alpha}{*}{*}{*}{*}0\alpha$ &
${*}{*}{*}{*}53\alpha$ &  
${*}{*}\alpha\alpha5{*}{*}$ & 
${*}{*}{*}\alpha53{*}$ & 
${*}0\alpha\alpha{*}{*}{*}$ & 
${*}{*}{*}{*}\tilde{\alpha}0\alpha$  \\
\arrayrulecolor{gray!50}\hline\arrayrulecolor{black}
12 & 
& 
6,8 & 
1,12 & 
1,12 & 
1,12 & 
12 &  \\
\hline
$\tilde{5}0\alpha\tilde{3}\alpha0\alpha$ & 
${*}{*}{*}{*}\tilde{\alpha}\tilde{0}\tilde{0}$ &  
$\tilde{\alpha}{*}{*}{*}{*}3\alpha$ & 
$\tilde{\alpha}3\alpha{*}{*}{*}{*}$ &
${*}{*}\tilde{0}\tilde{3}\tilde{\alpha}{*}{*}$ & 3
${*}{*}{*}{*}15\tilde{\alpha}$ &  
$5{*}{*}{*}{*}0\alpha$ &  
${*}{*}{*}{*}\tilde{\alpha}0\alpha$   \\
\arrayrulecolor{gray!50}\hline\arrayrulecolor{black}
13 & 
& 
1 & 
1,4,6 & 
3 & 
& 
&  \\
\hline
\end{tabular}
\caption{${\mc D}_1$-patches for neighboring tiles of non-regular $7^2_{5+6}$-heptagon.}
\label{H112tableA}
\end{table}

We observe that some partial information is not compatible with any of the 13 ${\mc D}_1$-patchs. This is indicated by empty boxes in the lower rows. On the other hand, some partial information is compatible with several ${\mc D}_1$-patches. This is indicated by more than one labels in the boxes in the lower rows. 

\begin{figure}[h!]
\centering
\begin{tikzpicture}[>=latex,scale=1]

\pgftext{
	\includegraphics[scale=0.4]{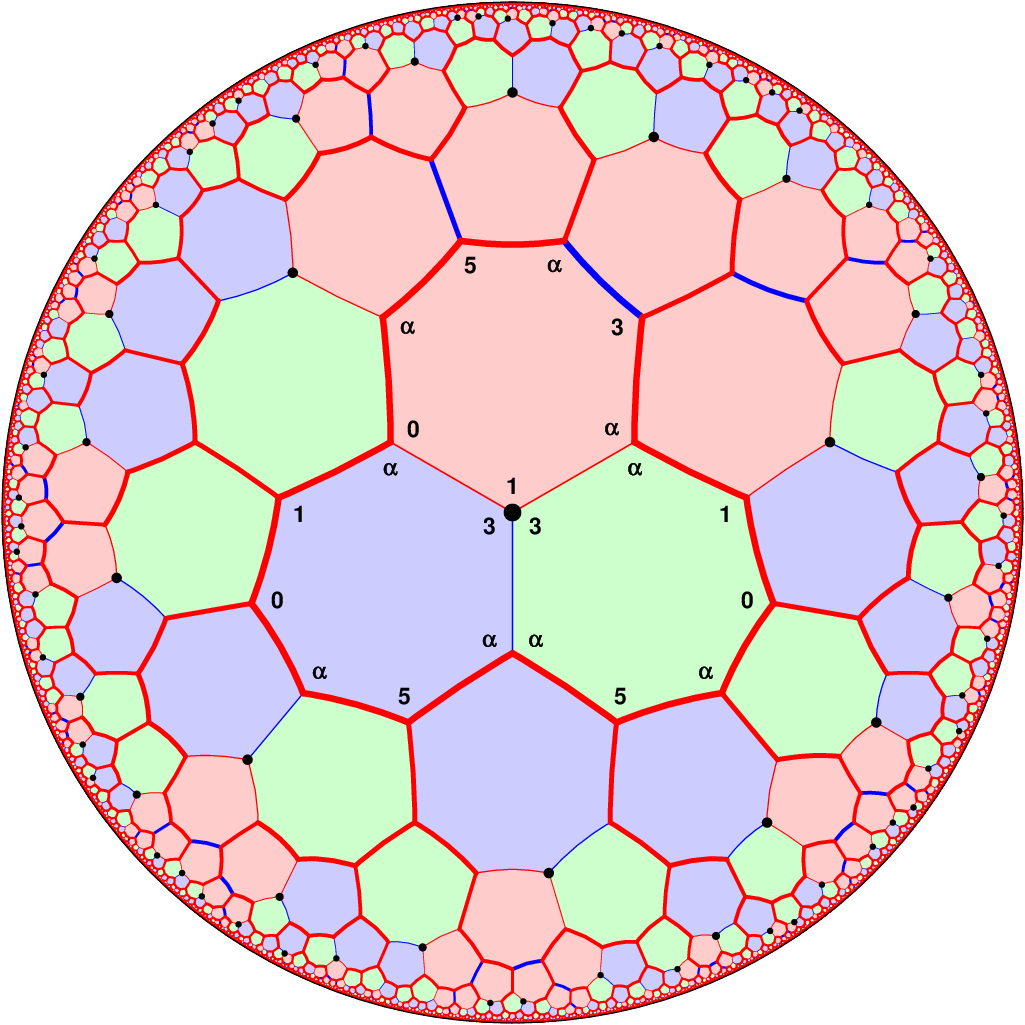}
	};

\foreach \a/\b in {12/2, 69/2.5, 90/2.65, 132/2, 190/2.5, 223/2.15, 276/2.2, -51/2.5, -30/2.65}	
\draw[->]
	(\a:0.9) -- (\a:\b);

\foreach \a in {79.2, 100.8, 199.2, -40.8, -19.2}	
\filldraw[fill=white]
	(\a:1.9) circle (0.05);
\foreach \a in {56.5, 176.5, -116.5, -63.5}	
\filldraw[fill=white]
	(\a:1.6) circle (0.05);
	
\foreach \a/\b in 
	{90/0, 132/1, 190/2, 223/3, 276/4, -51/5, -30/6, 12/7, 69/8} 
\node at (\a:0.7) {\scriptsize $g_{\b}$};

\begin{scope}[xshift=8cm]

\pgftext{
	\includegraphics[scale=0.4]{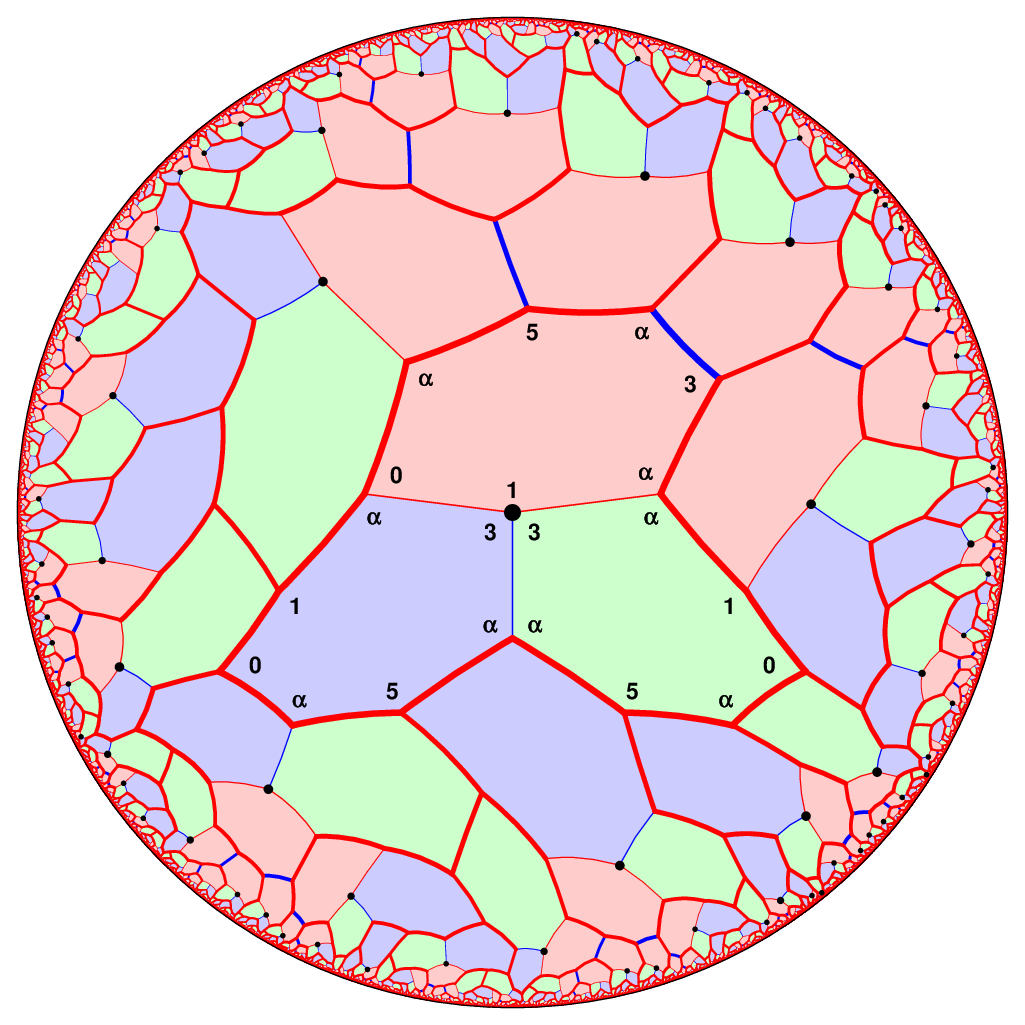}
	};

\foreach \a/\b in 
	{33/1.72, 56/1.72, 86.3/1.43, 198.5/1.72, 208.5/2.32, 241/1.62, -60.5/1.62, -43.5/2.13, -28.2/2.33}	
\filldraw[fill=white]
	(\a:\b) circle (0.05);
		
\end{scope}

\end{tikzpicture}
\caption{Schematic and real drawings of the unique tiling ${\mc T}(7^2_{5+6})$ of the Poincare disk ${\bb H}^2$ by the non-regular $7^2_{5+6}$-heptagon.}
\label{H112B}
\end{figure}
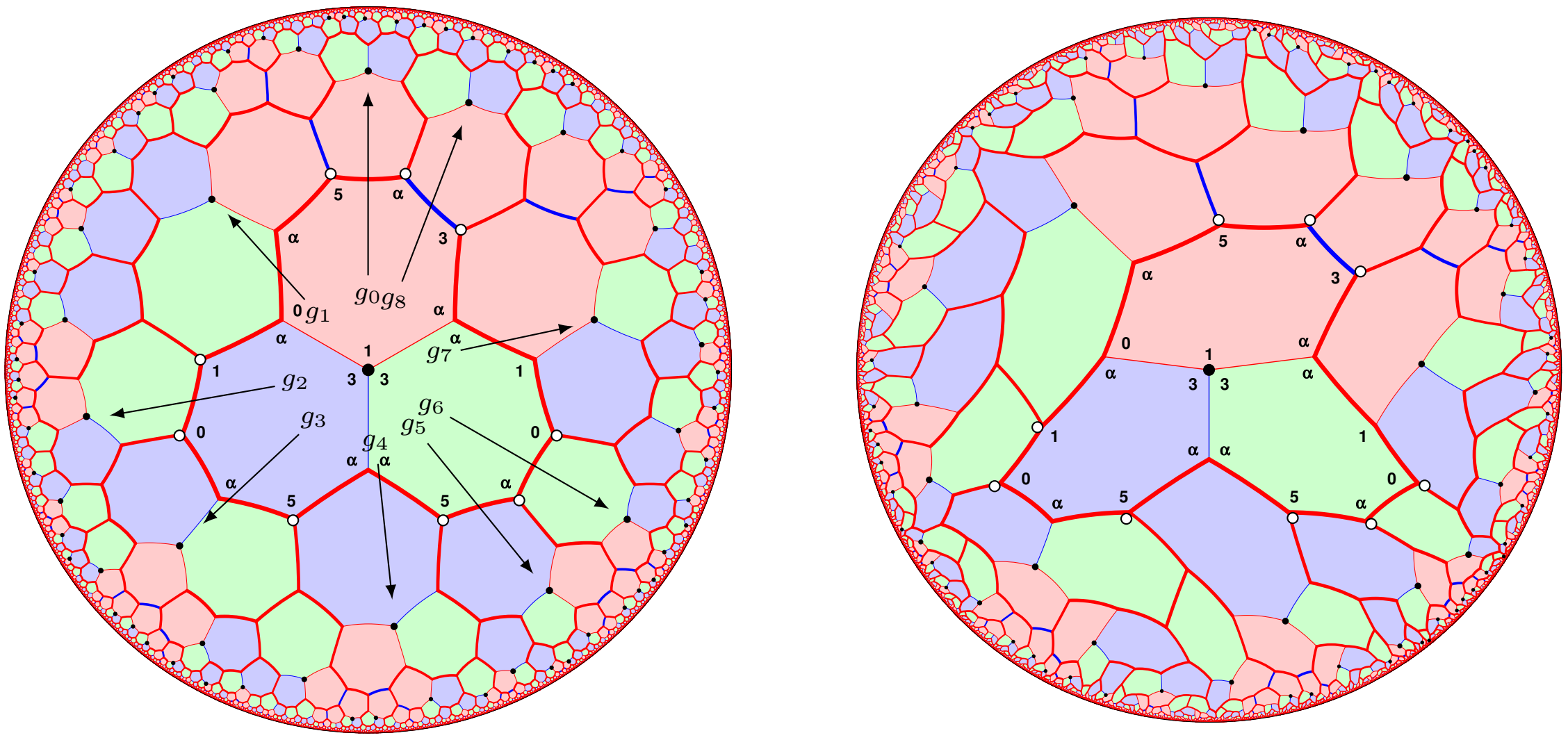

If some ${\mc D}_1$-patch in the first column has empty boxes on the right, then the ${\mc D}_1$-patch cannot be extended to ${\mc D}_2$-patch. Therefore the ${\mc D}_1$-patch cannot appear in the tiling. This means the ${\mc D}_1$-patches labeled by 6 though 13 are not allowed, and Table \ref{H112tableA} is reduced to a smaller table consisting of only ${\mc D}_1$-patches labeled by 1 through 5. The reduction of ${\mc D}_1$-patch labels to 1 through 5 also applies to all boxes in the lower rows. In particular, the 0th neighboring tile of the 4th ${\mc D}_1$-patch has partial information ${*}\alpha\tilde{0}\tilde{3}{*}{*}{*}$, which is only compatible with the 7th, 8th, and 9th ${\mc D}_1$-patches. Therefore this lower box becomes empty in the smaller table. We have the similar situation with the partial information ${*}{*}\tilde{0}\tilde{3}\alpha{*}{*}$ around the 3rd neighboring tile of the 5th ${\mc D}_1$-patch. Therefore the 4th and 5th ${\mc D}_1$-patches are also not allowed. Then the smaller table is further reduced to a table consisting of only the 1st, 2nd and 3rd ${\mc D}_1$-patches. We color the three ${\mc D}_1$-patches by red, blue and green, and Table \ref{H112tableA} is reduced to the colored part only. In the reduced table, each box contains only one ${\mc D}_1$-patch, i.e., one color. Therefore we are in the same unique situation similar to Table \ref{H111tableA}. The subsequent discussion is the same as the one for non-regular $7^1_{5+6}$-heptagon. We get a unique tiling of the Poincar\'e disk ${\bb H}^2$ by non-regular  $7^2_{5+6}$-heptagon in Figure \ref{H112B}. 

Similar to $7^1_{5+6}$, the triple $T$ consisting of tiles of three colors at the center is surrounded by nine triples $g_0T,g_1T,\dots,g_8T$, and we have the canonical presentation of the symmetry group
\begin{align*}
G(7^2_{5+6})
&=\langle g_0,g_1,g_2,g_3,g_4,g_5,g_6,g_7,g_8
\; | \\
&\qquad 
g_0^2=g_6^2=g_8^2=g_1g_7=g_2g_5=g_3g_4=g_0g_1g_8=g_2g_6g_7=g_3^2g_5=1\rangle \\
&=\langle g_0,g_3,g_6 \; | \;
g_0^2=g_6^2=(g_6g_3^{-2}g_0)^2=1\rangle,
\end{align*}
with
\[
g_1=g_3^2g_6,\;
g_2=g_3^2,\;
g_4=g_3^{-1},\;
g_5=g_3^{-2},\;
g_7=g_6g_3^{-2},\;
g_8=g_6g_3^{-2}g_0.
\]

Tilings of a closed surface of Euler characteristic $\chi$ by the non-regular  $7^2_{5+6}$-heptagon are the quotients of ${\bb H}^2$ by conjugate classes of torsionless subgroups of $G(7^2_{5+6})$ of finite index $\frac{f}{3}=-2\chi$. We know $g_3,g_4$ reverse the orientation, and all other generators preserve the orientation. The surface is orientable if and only if the numbers of $g_3,g_4$ (or only $g_3$ in the simplified presentation) in the generators of the corresponding torsionless subgroup are all even. We have a result similar to Proposition \ref{even6}.

\begin{proposition}\label{even112}
If the non-regular  $7^2_{5+6}$-heptagon tiles $g{\bb P}^2$, then $g$ is even.
\end{proposition}

\begin{proof}
Let $G$ be the torsionless subgroup of $G(7^2_{5+6})$ of index $-2\chi=-2(2-g)$, that gives a tiling ${\mc T}$ of $g{\bb P}^2$ by the non-regular $7^2_{5+6}$-heptagon. We have a natural homomorphism $x\in G\mapsto \bar{x}\in S_{-2\chi}$ to the permutation group of all triple tiles in ${\mc T}$. 

Since $g_0^2=g_6^2=g_8^2=1$ and $G$ is torsionless, similar to the proof of Proposition \ref{even6}, we know $\bar{g}_0,\bar{g}_6,\bar{g}_8$ act freely on the $-2\chi$ triple tiles of ${\mc T}$. Then by $\bar{g}_0^2=\bar{g}_6^2=\bar{g}_8^2=1$, we know each of $\bar{g}_0,\bar{g}_6,\bar{g}_8$ is a product of $\frac{-2\chi}{2}=-\chi$ disjoint transpositions. Therefore the signs of $\bar{g}_0,\bar{g}_6,\bar{g}_8$ are all $(-1)^{-\chi}$. Then the sign $(-1)^{-\chi}$ of $g_8=g_6g_3^{-2}g_0$ is the product $(-1)^{-\chi}(-1)^{-\chi}=1$ of the signs of $g_6$ and $g_0$. Therefore $(-1)^{-\chi}=1$, and $\chi$ is even, i.e., $g$ is even.
\end{proof}

\begin{figure}[htp]
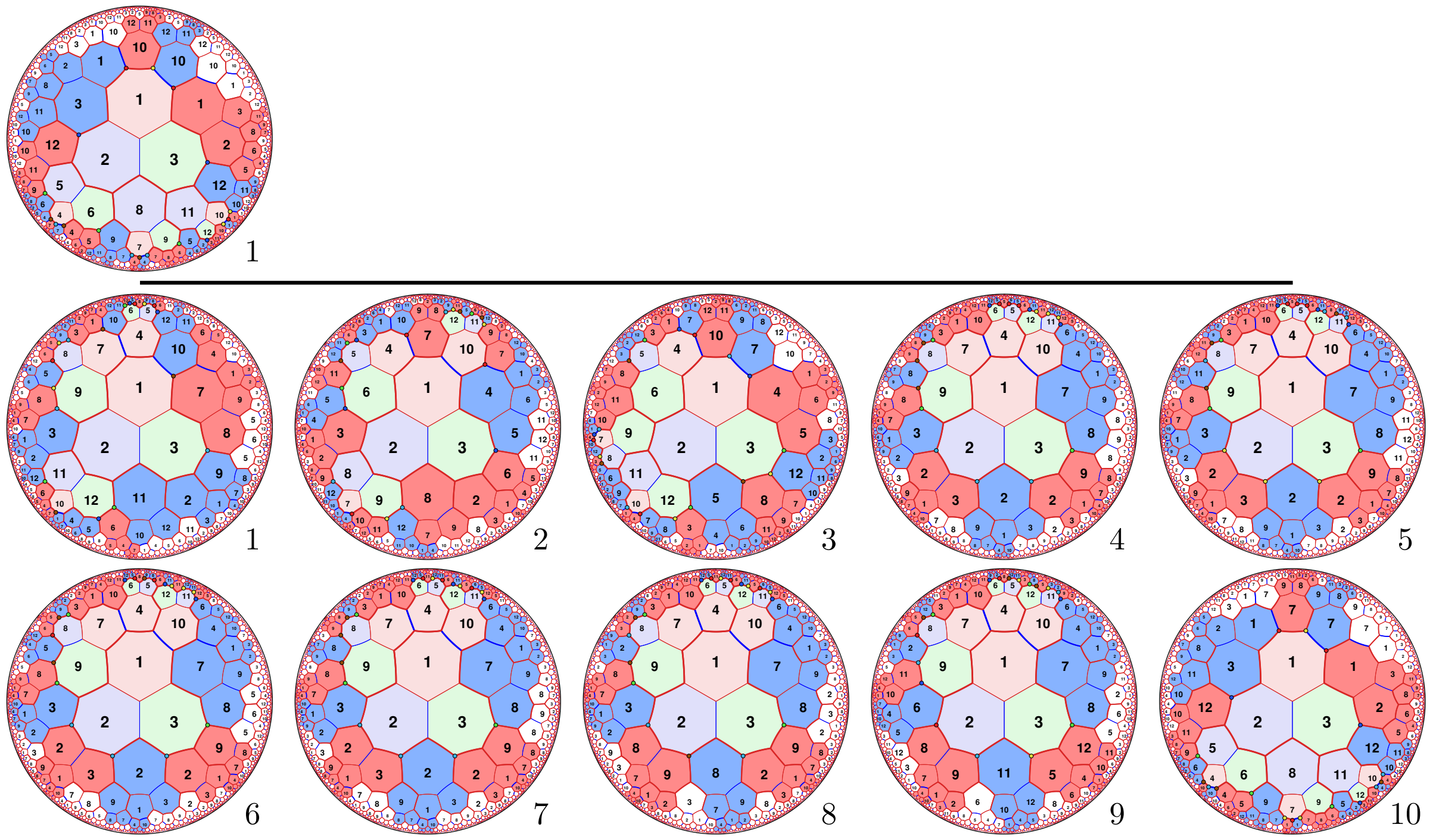

\centering
\begin{tikzpicture}[>=latex,scale=1]

%% 1 2T^2

\foreach \a in {1}
{
\begin{scope}[shift={(-3.2*1+3.2*\a,3.2)}]

\pgftext{
	\includegraphics[scale=0.166]{H112domainT/H112T\a}
	};

\node at (1.25,-1.25) {\a};

\end{scope}
}

%% 1-5 4T^2

\foreach \a in {1,...,5}
{
\begin{scope}[shift={(-3.2*1+3.2*\a,0)}]

\pgftext{
	\includegraphics[scale=0.166]{H112domainP/H112P\a}
	};

\node at (1.25,-1.25) {\a};

\end{scope}
}

%% 6-10

\foreach \a in {6,...,10}
{
\begin{scope}[shift={(-3.2*6+3.2*\a,-3.05*1)}]

\pgftext{
	\includegraphics[scale=0.166]{H112domainP/H112P\a}
	};

\node at (1.25,-1.25) {\a};

\end{scope}
}

\draw[very thick] 
	(0,1.6) -- ++(12.8,0);

\end{tikzpicture}
\caption{Fundamental domains of the only tiling ${\mc T}(7^2_{5+6},2{\bb T}^2)$ of $2{\bb T}^2$ by generic $7^2_{5+6}$-heptagon, and tilings ${\mc T}_k(7^2_{5+6},4{\bb P}^2)$, $k=1,2,\dots,10$, of $4{\bb P}^2$ by generic $7^2_{5+6}$-heptagon.}
\label{H112C}
\end{figure}

We use standard computer package GAP to find the numbers of tilings of $g{\bb T}^2$, $g=2,3,4$, in Table \ref{table5}. There is only one tiling ${\mc T}(7^2_{5+6},2{\bb T}^2)$ of $2{\bb T}^2$ by the generic $7^2_{5+6}$-heptagon. There are 10 tilings ${\mc T}_k(7^2_{5+6},4{\bb P}^2)$ of $4{\bb P}^2$ by the generic $7^2_{5+6}$-heptagon. Figure \ref{H112C} presents the lightly colored fundamental domains and the red/blue neighboring domains.  

\begin{figure}[htp]
\centering
\begin{tikzpicture}[>=latex,scale=1]

\pgftext{
	\includegraphics[scale=0.4]{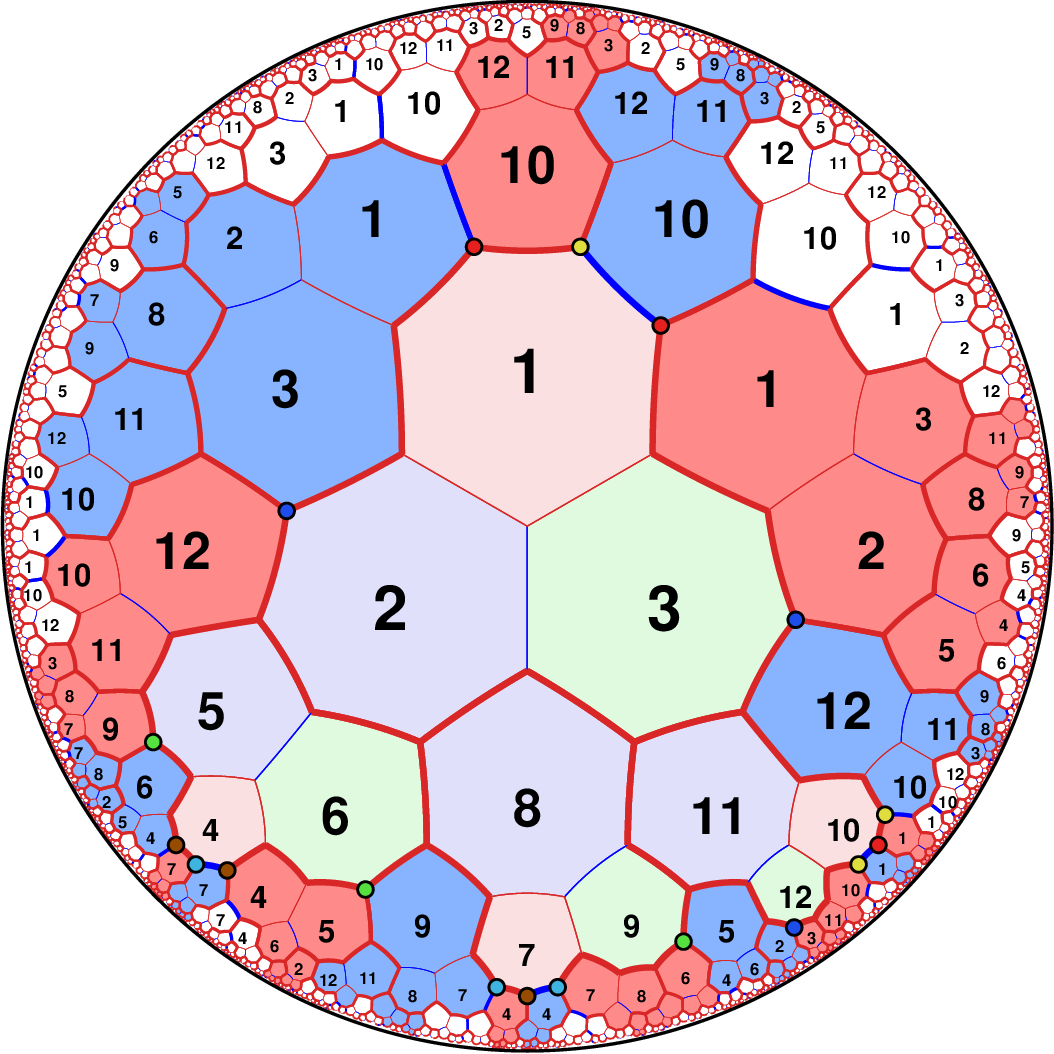}
	};

\begin{scope}[xshift=8cm]

\pgftext{
	\includegraphics[scale=0.4]{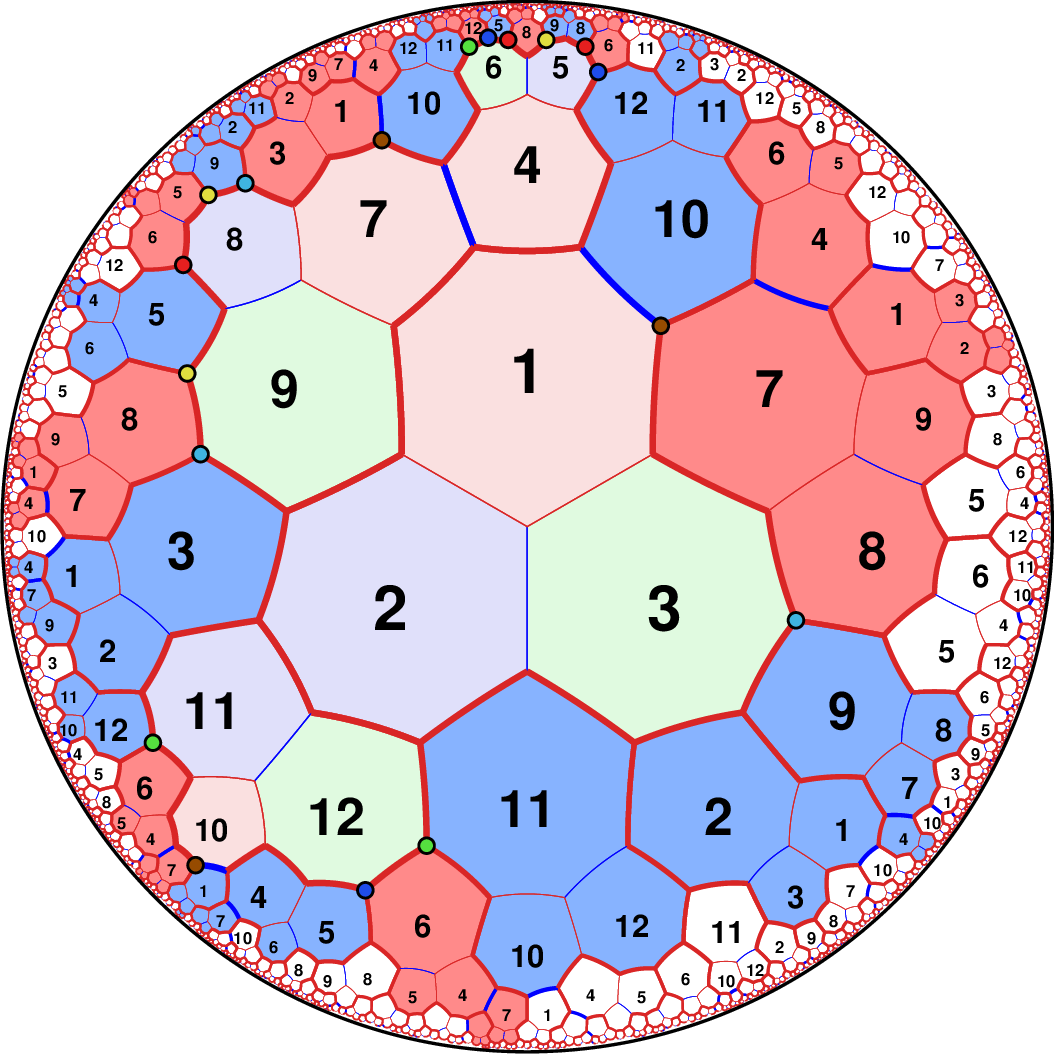}
	};

\foreach \a/\b in {
	80/1, 85/2, 90/3, 94/4, 98/5, 104/6, 115/7, 128/8, 142/9, 154/10, 172/11, 195/12, 220/13, 240/14, 260/15, 300/16, 30/17, 70/18}
\node at (\a:3.8) {\small \b};

\end{scope}

\end{tikzpicture}
\caption{The only tiling ${\mc T}(7^2_{5+6},2{\bb T}^2)$ of $2{\bb T}^2$, and first tiling ${\mc T}_1(7^2_{5+6},4{\bb P}^2)$ of $4{\bb P}^2$, in Figure \ref{H112C}, by the non-regular $7^2_{5+6}$-heptagon.}
\label{H112E}
\end{figure}

Figure \ref{H112E} is the detail of the only tiling ${\mc T}(7^2_{5+6},2{\bb T}^2)$ of $2{\bb T}^2$ and the first tiling ${\mc T}_1(7^2_{5+6},4{\bb P}^2)$ of $4{\bb P}^2$ in Figure \ref{H112C}. The lightly colored fundamental domain $d$ consists of four triple tiles $T,g_3T,g_4T,g_5T$ in ${\mc T}(7^2_{5+6},2{\bb T}^2)$, and consists of four triple tiles $T,g_0T,g_1T,g_3T$ in ${\mc T}_1(7^2_{5+6},4{\bb P}^2)$. The canonical presentation for the torsionless subgroup $G(7^2_{5+6},2{\bb T}^2)$ is similar to $G_1(7^1_{5+6},2{\bb T}^2)$ in Section \ref{H111}. We present here the canonical presentation of the group $G_1(7^2_{5+6},4{\bb P}^2)$ for tiling of non-orientable surface. Using the labels for the neighboring domains on the right of Figure \ref{H112E}, the canonical generators of $G_1(7^2_{5+6},4{\bb P}^2)$ are
\begin{align*}
a_1 &=g_0g_2g_0, & 
a_4 &=g_0g_5g_0, &
a_7 &=g_1^2, & 
a_{10} &=g_1g_4g_0, & 
a_{13} &=g_3g_1g_0, & 
a_{16} &=g_5, & \\
a_2 &=g_0g_3g_7, & 
a_5 &=g_0g_6g_4, &
a_8 &=g_1g_2g_7, &
a_{11} &=g_1g_5g_7, & 
a_{14} &=g_3g_7g_0, & 
a_{17} &=g_7^2,  \\
a_3 &=g_0g_4g_7, & 
a_6 &=g_0g_7g_4, &
a_9 &=g_1g_3g_0, &
a_{12} &=g_2, & 
a_{15} &=g_3g_6g_0, & 
a_{18} &=g_8g_4. 
\end{align*}
The torsionless subgroup $G_1(7^2_{5+6},4{\bb P}^2)$ is generated by the 18 isometries $a_i$, subject to the relations 
\[
a_1a_4=
a_2a_{10}=
a_3a_9=
a_5a_{15}=
a_6a_{13}=
a_7a_{17}=
a_8a_{11}=
a_{12}a_{16}=
a_{14}a_{18}=1,
\]
and 
\[
{\color{red} a_1a_3a_{10}}
={\color{blue} a_4a_18a_{15}}
={\color{green} a_5a_{16}a_{13}}
={\color{brown} a_6a_{14}a_{17}}
={\color{cyan} a_8a_{12}a_{17}}
={\color{magenta} a_2a_{11}a_9}
=1.
\]

\end{document}